\documentclass[11pt,twoside]{article}
\usepackage{latexsym, amssymb, amsthm, amsmath, epsfig, setspace, tikz, slashed}
\usepackage{mathrsfs}
\usepackage[toc,page]{appendix}
\usepackage[all]{xy}
\usepackage[retainorgcmds]{IEEEtrantools}
\usetikzlibrary{decorations.markings}
 \usetikzlibrary{arrows}
\usepackage{hyperref}
\usepackage[margin=10pt,font=small,labelfont=bf]{caption}

\theoremstyle{plain}
\newtheorem{theorem}{Theorem}[section]
\newtheorem{proposition}[theorem]{Proposition}

\theoremstyle{definition}
\newtheorem{definition}[theorem]{Definition}

\newtheorem{lemma}[theorem]{Lemma}
\newtheorem{corollary}[theorem]{Corollary}

\newtheorem{example}[theorem]{Example}

\newenvironment{renumerate}%
{%
\begin{enumerate}}%
{\end{enumerate}%
}%

\newenvironment{remark}%
{\vskip6pt%
\noindent%
{\it Remark.}}%
{\vskip6pt}

{\vskip6pt%
\noindent%
{\it Remarks}. %
\begin{renumerate}}%
{\end{renumerate}\vskip6pt}

\makeatletter
\def\Ddots{\mathinner{\mkern1mu\raise\p@
\vbox{\kern7\p@\hbox{.}}\mkern2mu
\raise4\p@\hbox{.}\mkern2mu\raise7\p@\hbox{.}\mkern1mu}}
\makeatother

\newcommand{\R}{\text{${\mathbb R}$}}
\newcommand{\C}{\text{$\mathbb C$}}

\newcommand{\KK}{\text{$\mathbb K$}}
\newcommand{\Z}{\text{$\mathbb Z$}}
\newcommand{\GG}{\text{$\mathbb G$}}

\newcommand{\D}{\slashed{D}}
\newcommand{\T}{\text{$\mathbb{T}$}}
\newcommand{\TM}{\text{$\mathbb{T}M$}}
\renewcommand{\frak}[1]{\text{$\mathfrak{#1}$}}
\renewcommand{\tilde}{\widetilde}

\renewcommand{\L}{\text{$\mathcal{L}$}}

\renewcommand{\AA}{\text{$\mathbb{A}$}}
\newcommand{\V}{\text{$\mathbb{V}$}}
\newcommand{\J}{\text{$\mathcal{J}$}}
\newcommand{\I}{\text{$\mathcal{I}$}}
\newcommand{\W}{\text{$\mathcal{W}$}}
\newcommand{\JJ}{\text{$\mathbb{J}$}}
\newcommand{\M}{\text{$\mathcal{M}$}}

\newcommand{\N}{\text{$\mathcal{N}$}}
\newcommand{\G}{\text{$\mathcal{G}$}}
\newcommand{\U}{\text{$\mathcal{U}$}}
\newcommand{\Gau}{\mathscr{G}}
\newcommand{\Clif}{\mathrm{Clif}}

\newcommand{\Spin}[1]{\mathrm{Spin}({#1})}

\newcommand{\e}{\text{$\varepsilon$}}

\renewcommand{\bar}{\overline}

\newcommand{\gf}{\text{$\varphi$}}

\newcommand{\Id}{\mathrm{Id}}

\newcommand{\del}{\text{$\partial$}}

\newcommand{\delbar}{\text{$\overline{\partial}$}}
\newcommand{\deltabar}{\text{$\overline{\delta}$}}
\newcommand{\tensor}{\otimes}

\newcommand{\Kperp}{\mathbb{K}^{\perp}}
\newcommand{\mc}[1]{\text{$\mathcal{#1}$}}

\newcommand{\into}{\longrightarrow}
\newcommand{\noqed}{\let\qed\relax}

\renewcommand{\Im}{\mathrm{Im}\,}

\newcommand{\Gg}{\mathfrak{g}}

\newcommand{\IP}[1]{\langle #1 \rangle}

\newcommand{\nhood}{neighbourhood}

\newcommand{\gcs}{generalized complex structure}
\newcommand{\gacs}{generalized almost complex structure}

\newcommand{\gcss}{generalized complex structures}

\newcommand{\gc}{generalized complex}

\newcommand{\gk}{generalized K\"ahler}
\newcommand{\gks}{generalized K\"ahler structure}
\newcommand{\gkss}{generalized K\"ahler structures}
\newcommand{\gkm}{generalized K\"ahler manifold}

\newcommand{\wrt}{with respect to}
\renewcommand{\iff}{if and only if}

\newcommand{\Cour}[1]{[\![#1]\!]}
\newcommand{\gsd}{\text{$\mathrm{SD}$}}
\newcommand{\gasd}{\text{$\mathrm{ASD}$}}

\date{} \usepackage{color} \definecolor{tocolor}{rgb}{.1,.1,.5}
\definecolor{urlcolor}{rgb}{.2,.2,.6}
\definecolor{linkcolor}{rgb}{.1,.1,.6}
\definecolor{citecolor}{rgb}{.6,.2,.1}
\hypersetup{backref=true, colorlinks=true, urlcolor=urlcolor,
  linkcolor=linkcolor, citecolor=citecolor}

\begingroup
\hypersetup{linkcolor=blue}
\endgroup

\numberwithin{equation}{section}

\makeatletter
\newcommand{\subjclass}[2][2010]{%
  \let\@oldtitle\@title%
  \gdef\@title{\@oldtitle\footnotetext{#1 \emph{Mathematics subject classification.} #2}}%
}
\newcommand{\keywords}[1]{%
  \let\@@oldtitle\@title%
  \gdef\@title{\@@oldtitle\footnotetext{\emph{Keywords.} #1.}}%
}
\makeatother

\begin{document}
\title{Hodge theory and deformations of SKT manifolds}
\author{Gil R. Cavalcanti\thanks{{\tt gil.cavalcanti@gmail.com}} \\
       Department of Mathematics\\
Utrecht University\\
}

\maketitle

\abstract{We use tools from generalized complex geometry to develop the theory of SKT (a.k.a. pluriclosed Hermitian) manifolds and more generally manifolds with special holonomy \wrt\ a metric connection with closed skew-symmetric torsion. We develop Hodge theory on such manifolds showing how the reduction of the holonomy group causes a decomposition of the twisted cohomology. For SKT manifolds this decomposition is accompanied by an identity between different Laplacian operators and forces the collapse of a spectral sequence at the first page. Further we study the deformation theory of  SKT structures, identifying the space where the obstructions live. We illustrate our theory with examples based on Calabi--Eckmann manifolds, instantons, Hopf surfaces and Lie groups.
}

\tableofcontents

\section*{Introduction}

Looking beyond the Levi--Civita connection in Riemannian geometry, one finds a number  of other metric connections with interesting properties. Normally these families of connections are defined by some characteristic of the  torsion tensor and recurrent themes of research are connections with parallel torsion or connections with skew symmetric torsion. The ``strong" torsion condition refers to the latter: a {\it strong torsion} connection on a Riemannian manifold is a metric connection whose torsion is skew symmetric and closed. In this setting, a {\it K\"ahler structure with strong torsion}, or {\it SKT structure} is a Hermitian structure $(g,I)$ together with a strong torsion connection for which $I$ is parallel. A weaker notion is that of a (strong) parallel Hermitian structure which for us means  a connection with closed, skew symmetric torsion for which the holonomy is $U(n)$. The difference between a parallel Hermitian structure and an SKT structure being that in the former case integrability of the complex structure is not required.

A reason to study of such objects comes from  string theory, where closed 3-forms arise naturally as fields in their sigma models \cite{Neveu197186,MR0342064}. Once a 3-form is added to the sigma model, if one still requires a nontrivial amount of supersymmetry, the type of geometry of the target space has to move away from the usual K\"ahler geometry. It was precisely following this path that Gates, Hull and Ro\v cek \cite{MR776369} discovered the bi-Hermitian geometry that nowadays also goes by the name of generalized K\"ahler geometry \cite{gualtieri-2010} as the solutions to the $(2,2)$-supersymmetric sigma model. Requiring less supersymmetry without giving up on the idea altogether  leads one to consider models where there is more left than right supersymmetry or models where the right side is simply absent. These conditions lead to  $(2,1)$ or $(2,0)$ supersymmetric sigma models and supersymmetry holds  \iff\  the target space  has an SKT structure \cite{MR818681}. This point of view also leads one to consider parallel Hermitian and bi-Hermitian structures as these are geometric structures imposed by a sigma model with  an extended supersymmetry algebra \cite{MR981114,MR1025428}.

Mathematical properties of SKT structures have been subject of study of several papers since the 90's, including those of Bismut \cite{MR1006380}, Grantcharov {\it et al} \cite{MR1908699},  Fino and collaborators \cite{MR2059435,MR2511043} and more recently Streets and Tian \cite{MR2673720,MR2755684}. Yet, until now there was no framework in which basic K\"ahler results, such as Hodge theory and deformation theory, obtained meaningful counterparts in the SKT world. Indeed the opposite seemed to be the case: Since SKT manifolds are, in particular, complex, their space of forms inherits a natural bi-grading, but a simple check in concrete examples shows that there is no corresponding  decomposition of their cohomology. Further, by studying these structures on six dimensional nilmanifolds Fino, Parton and Salamon \cite{MR2059435} produced examples showing that their Fr\"olicher spectral sequence does not necessarily degenerate at the first page, they do not satisfy the $dd^c$-lemma, manifolds carrying these structures may not be formal and  that these structures are not stable under deformations of the complex structure. In short, several K\"ahler properties seem to have been lost once the torsion was included. Put another way, given a complex manifold, we were left with no tools to decide whether it admited an SKT metric or not.

Here we tackle SKT structures from a new point of view. The key observation is that SKT structures have yet another description, this time, as a `generalized structure', i.e., a geometric structure on $\TM = TM\oplus T^*M$. In fact, in this paper we show that, in a very precise way, SKT structures lie half way between generalized Hermitian and \gkss.  Using this approach, we show that some of the negative results mentioned earlier have a positive counterpart involving the torsion while those with negative answers obtain a conceptual explanation for their failure to hold.

Indeed, the first observation is that, as structures defined on $\TM$ with the $H$-Courant bracket, the natural differential operator to consider is $d^H = d + H \wedge$, where $H = d^c \omega$, and $\omega$ is the Hermitian form. Hence questions about $d$, and its decomposition as $\del + \delbar$ miss an important ingredient and were doomed from the start. The cohomology of $d^H$ is only $\Z_2$-graded, yet we show that a parallel Hermitian structure induces a $\Z\times \Z_2$ grading on the space of forms which itself induces a $\Z\times \Z_2$ grading on the $d^H$-cohomology. This is achieved by introducing the {\it intrinsic torsion} of a generalized almost Hermitian structure and using it to show that the $d^H$-Laplacian preserves a $\Z \times \Z_2$-graded decomposition of the space of forms. For  SKT manifolds one can go further and prove an identity of Laplacians, extending Gualtieri's work on  generalized K\"ahler geometry \cite{gualtieri-2004}.  This way we relate a cohomology naturally defined in terms of the SKT data with the $d^H$-cohomology. These identities imply further that a variation of the Fr\"olicher spectral sequence degenerates at the first page. As an application, we return to the moduli space of instantons on a bundle over a compact complex surface and show that the existence of an SKT structure in this space (obtained by L\"ubke and Teleman \cite{MR1370660} by an {\it ad hoc} method) can be seen as a consequence of the Hodge theory developed for SKT manifolds.

With this new framework at hand,  we study of the problem of stability of SKT structures. Precisely, we use the operators introduced by the generalized complex setup to determine the obstruction space to the stability problem. In its classical version, we prove that if the Dolbeault cohomology group $H^{1,2}(M)$ vanishes, then any deformation of the complex structure can be accompanied by a deformation of the metric so that the pair remains SKT.

This paper is organized as follows. In Section \ref{sec:linear algebra} we develop the linear algebra pertinent to generalized complex, generalized Hermitian and SKT structures. In particular we show that an SKT structure gives rise to a $\Z\times \Z_2$-grading on the space of forms. In Section  \ref{sec:intrinsic torsion} we introduce the intrinsic torsion of a generalized Hermitian structure and, in Sections \ref{sec:parallel} and \ref{subsec:integrability} , we relate SKT structures and parallel Hermitian structures to the vanishing of certain components of the intrinsic torsion. In Section \ref{sec:Hodge} we study Hodge theory for SKT  and parallel Hermitian structures and prove that in both cases the $d^H$-cohomology decomposes according to the decomposition of forms induced by the structure. As an application, in Section \ref{subsec:instantons} we recover the result that the moduli space of instantons over a complex surface has an SKT structure. In  Section \ref{sec:deformations} we use our framework  to study the  problem of deformations of SKT and  \gkss.

{\bf Acknowledgements}: This research was supported by a Marie Curie Intra European Fellowship and a VIDI grant from NWO, the Dutch science foundation. The author is thankful to Anna Fino, S\"onke Rollenske,  Ulf Lindstrom, Martin Ro\v cek, Stefan Vandoren and  Maxim Zabzine for useful conversations.

\section{Linear algebra}\label{sec:linear algebra}

Given a vector space $V^m$ we let $\V = V \oplus V^*$ be its ``double".  $\V$ is endowed with a natural symmetric pairing:
$$\IP{X+\xi,Y+\eta} =\frac{1}{2}(\eta(X) + \xi(Y)),\qquad X,Y \in V~\xi,\eta \in V^*.$$

Elements of $\V$ act on $\wedge^{\bullet}V^*$ via
$$(X +\xi) \cdot \gf = i_X\gf + \xi\wedge \gf.$$
One can easily check that for $v \in \V$
$$v\cdot(v\cdot \gf) = \IP{v,v}\,\gf,$$
hence $\wedge^{\bullet} V^*$ is naturally a module for the Clifford algebra of $\V$. In fact, it is the space of spinors for $Spin(\V)$ and hence comes equipped with a spin invariant pairing, the {\it Chevalley pairing}:
$$(\gf,\psi)_{Ch} = -(\gf \wedge \psi^t)_{top}, $$
where $\cdot^t$ indicates transposition, an $\R$-linear operator defined on decomposable forms by
$$(\alpha_1\wedge \cdots \wedge \alpha_k)^t = \alpha_k \wedge \cdots \wedge \alpha_1,$$
and $top$ means taking the degree $m$ component.

The spin group, $\Spin\V$, acts on both $\V$ and on spinors in a compatible manner, namely, its action on $\V$ is by conjugation using Clifford multiplication
$$g_* v = gvg^{-1}\qquad \mbox{ for all } g \in \Spin\V, v \in \V$$
and on $\wedge^{\bullet} V^*$ by the Clifford action described above, so we have
\begin{equation}\label{eq:spin action}
g\cdot (v\cdot \gf) = (gvg^{-1}) \cdot (g\gf) = (g_* v) \cdot(g\cdot \gf)
\end{equation}
for all $g\in \Spin\V$, $v \in \V$ and $\gf \in \wedge^{\bullet}V^*$.


\begin{example}[$B$-field transform]\label{ex:B field}
Given $B \in \wedge^2 V^*\subset \frak{spin}(\V) $, it acts on $\V$. Namely, we can regard $B$ as a map from $V$ into $V^*$
$$B(X+\xi)= -i_XB = -B(X).$$
Exponentiating this map we have
$$e^B_*:\V \into \V \qquad e^B_*(X+\xi) = X + \xi -B(X).$$ 
And this action is compatible with the action of $e^B \in \mathrm{Spin}(\V)$ on forms:
$$e^B:\gf \into e^B \wedge \gf.$$
So we have $e^B_*:\wedge^{k} \V \into \wedge^{k}\V$ and for $\alpha \in \wedge^{\bullet}\V$ and $\gf \in \wedge^{\bullet}V^*$
$$(e^B_*\alpha) \cdot e^B \gf = e^B(\alpha \cdot \gf).$$
\hfill$\blacksquare$
\end{example}

We will be interested in introducing geometric structures on $\V$. The first we consider is a generalized metric, as introduced by Gualtieri \cite{gualtieri-2010}.

\begin{definition}
A {\it generalized metric} on $V$ is an automorphism $\G:\V \into \V$ which is orthogonal and self-adjoint \wrt\ the natural pairing and for which the bilinear tensor
$$\IP{\G v,w}, \qquad v,w \in \V$$
is positive definite.
\end{definition}

Since $\G$ is orthogonal and self-adjoint, we have $\G^{-1} = \G^t =\G$, hence $\G^{2}= \Id$. Therefore $\G$ splits $\V$ into its $\pm 1$-eigenspaces: $\V = V_+ \oplus V_-$ and the projection $\pi_V:\V \into V$ gives isomorphisms $\pi:V_{\pm} \into V$. Further, given a generalized metric $\G$ we can write $\V = \G V^* \oplus V^*$ and $\G V^*$ is isomorphic to $V$ via the projection $\pi_V:\V \into V$. Since both $V$ and $\G V$ are isotropic subspaces of $\V$ which project isomorphically onto $V$, we can describe $\G V$ as the graph of a linear map $b: V\into V^*$, that is $b \in \tensor^2 V^*$. Isotropy means that $b \in \wedge^2 V^*$ and hence gives rise to an orthogonal transformation of the natural pairing, $e^b_*$. This map has the property that $e^b_*:\G V^* \into V$, hence, after an orthogonal transformation of $\V$, we can assume that $\G V^* = V$. For this splitting, 
\begin{equation}\label{eq:generalized metric}
\G = \begin{pmatrix}
0& g\\
g^{-1} & 0
\end{pmatrix}\end{equation}
where $g$ is an ordinary metric on $V$. The splitting of $\V$ determined by a generalized metric is the {\it metric splitting}.

If $V$ is endowed with an orientation, we can define a generalized {\it Hodge star operator} \cite{gualtieri-2004} as follows. Since $\pi_V:V_+ \into V$ is an isomorphism, $V_+$  inherits an orientation. Then we let $\{e_1, e_2,\cdots, e_m\}$ be a positive orthonormal basis of $V_+$ and let $\star = -e_m\cdot\cdots e_2\cdot e_1 \in \mathrm{Clif}(\V)$. Then
$$\star \gf := \star \cdot \gf,$$
where $\cdot$ denotes Clifford action.

With this definition, we have
\begin{equation}\label{eq:positive pairing}
(\gf,\star \gf)_{Ch} >0 \qquad \mbox{if}~ \gf \neq 0
\end{equation}
If the splitting of $\V$ is the metric splitting, we have
$$(\gf,\star \psi)_{Ch} = \gf \wedge * \psi,$$
where $*$ is the usual Hodge star, hence, in this splitting, $\star$ is the usual Hodge star except for a change in signs given by the Chevalley pairing. Since $\star =  -e_m\cdot\ldots \cdot e_2\cdot e_1$, we have that
$$\star^2 = (-1)^{\frac{m(m-1)}{2}}$$
and hence it splits the space of forms into its eigenspaces, namely, into its $\pm 1$-eigenspaces $\wedge_{\pm}^{\bullet}V^*$ if $m$ is zero or one modulo four or its $\pm i$-eigenspaces if $n$ is 2 or 3 modulo 4. This allows us to define {\it self-dual} and {\it anti self-dual forms} in all dimensions:  
\begin{definition}
We say that a form $\gf \in \wedge^{\bullet}V^*_\C$  is {\it \gsd} if $\star \gf = -i^{\frac{m(m-1)}{2}}\gf$ and  that it is {\it \gasd} if $\star\gf = i^{\frac{m(m-1)}{2}}\gf$. We denote the space of \gsd-forms by $\wedge_{+}^{\bullet}V^*$ and the space of \gasd-forms by $\wedge_{-}^{\bullet}V^*$.\footnote{The choice of signs on the Chevalley pairing and of $\star$ were made so that on a four manifold, the notion of \gsd\ and \gasd\ agrees with the usual notions of self-dual and anti self-dual on 2-forms.}
\end{definition}

The Clifford action of elements in $V_{\pm}$ either preserves or switches the eigenspaces of $\star$:

\begin{lemma}\label{lem:V+ V-}
Let $v_\pm \in V_{\pm}$. Then acting via Clifford action on forms we have
$$v_+ \star = (-1)^{m-1}\star v_+\qquad\mbox{ and }\qquad v_-\star =  (-1)^m\star v_-$$
hence for $m$ even, 
$$V_+:\wedge_{\pm}^\bullet V^* \into \wedge_{\mp}^\bullet V^*\qquad \mbox{ and }\qquad V_-:\wedge_{\pm}^\bullet V^* \into \wedge_{\pm}^\bullet V^*$$
and for $m$ odd
$$V_+:\wedge_{\pm}^\bullet V^* \into \wedge_{\pm}^\bullet V^*\qquad \mbox{ and }\qquad V_-:\wedge_{\pm}^\bullet V^* \into \wedge_{\mp}^\bullet V^*$$
\end{lemma}
\begin{proof}
For $v_- \in V_-$, we have that $\IP{v_-,V_+} =0$ hence $v_-$ graded commutes with $\star$. For $v_+\in V_+$, we can choose an orthonormal  basis for $V_+$ for which $e_1 = v/\|v\|$, so that $v$ anti commutes with all the remaining elements of the basis and commutes with $e_1$. 
\end{proof}

For the rest of this section we will introduce structures on $V$ which force its dimension to be even so we let $m =2n$.

\begin{definition}\label{def:linear gcs}
A {\it \gcs} on $V$ is a complex structure on $\V$ which is orthogonal \wrt\ the natural pairing. A {\it generalized Hermitian structure}  or a {\it $U(n)\times U(n)$ structure} on $V$ is a \gcs\ $\J_1$ on $V$ and a generalized metric $\G$ such that $\J_1$ and $\G$ commute.
\end{definition}

Given a generalized complex structure $\J$ on $V$, we can split $\V_\C$, the complexification of $\V$, into the $\pm i$-eigenspaces of $\J$: $\V_\C = L \oplus \bar{L}$. The spaces $L$ and $\bar{L}$ are maximal isotropic subspaces of $\V_\C$ such that $L \cap \bar{L} =\{0\}$. Since the natural pairing is nondegenerate, we can use it to identify $\bar{L} = L^*$. Precisely, we identify  
\begin{equation}\label{eq:Lbar and L*}
v\mapsto  2\IP{v,\cdot} \in L^*,\qquad \mbox{ for } v \in \bar{L}.
\end{equation}

Given a generalized Hermitian structure $(\J_1,\G)$ on $V$, $\J_2 = \G \J_1$ is orthogonal \wrt\ the natural pairing and squares to $-\Id$, hence it is also a \gcs. Since $\pi_V:V_{\pm}\into V$ are isomorphisms, and $\mc{J}_1|_{V_{\pm}}$ is a  complex structure on $V_{\pm}$ orthogonal \wrt\ the natural pairing, it  induces complex structures $I_{\pm}$ on $V$ compatible with the metric $g$ induced by $\G$ making $V$ into a bi-Hermitian vector space. We can further form the corresponding Hermitian forms $\omega_{\pm} = g\circ I_{\pm}$.

Given any generalized Hermitian structure, $\wedge^{\bullet}V^*_\C$ splits as the intersections of the eigenspaces of $\J_1$ and $\J_2$: $U^{p,q} = U^p_{\J_1}\cap U^q_{\J_2}$, where $U^p_{\J_i}$ is the $ip$-eigenspace of $\J_i$ in $\wedge^\bullet V^*_{\C}$.  In this context, the generalized Hodge star is related to the action of $\JJ_i = e^{\frac{\pi \J_i}{2}}$, namely:

\begin{lemma}\label{lem:star} { (Gualtieri \cite{gualtieri-2004})} In a generalized Hermitian vector space one has
$$\star = -\JJ_1\JJ_2.$$
\end{lemma}

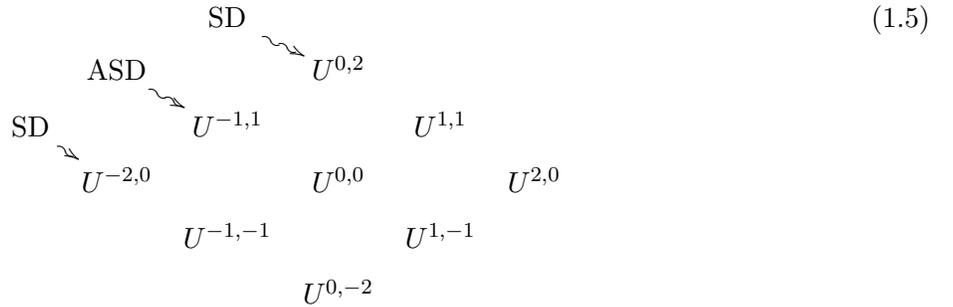
\begin{figure}[h!!]
\begin{center}
\begin{equation}
\xymatrix@R=6pt@C=6pt{& & \gsd \ar@{~>}[rd]&&&\\
&          \gasd \ar@{~>}[rd]      &               &U^{0,2}& &\\
\gsd\ar@{~>}[rd] &               &U^{-1,1}&              & U^{1,1}&\\
&U^{-2,0} &               &U^{0,0} &              & U^{2,0}\\
&              & U^{-1,-1} &            & U^{1,-1}&\\
&              &                  &U^{0,-2} &            &
}
\end{equation}
\caption{Representation of the spaces of \gsd\ and \gasd\ in terms of the $(p,q)$-decomposition of forms on a 4-dimensional generalized Hermitian structure.}\label{fig:figure1}
\end{center}
\end{figure}

This means that we can read the decomposition of forms into \gsd\ and \gasd\ from the $U^{p,q}$ decomposition, namely $\star|_{U^{p,q}} = i^{p+q}$. If we plot the (nontrivial) spaces $U^{p,q}$ in a lattice, each diagonal is made either of \gsd- or \gasd-forms, with $U^{n,0}$ made of \gsd-forms (see Figure \ref{fig:figure1}).

\begin{definition}\label{def:linear SKT}
A {\it positive $U(n)$ structure} or a {\it positive Hermitian structure} on $V$ is a generalized metric $\G$ and a complex structure $\I_+$ on $V_+$, the $+1$-eigenspace of $\G$, orthogonal \wrt\ the natural pairing. A {\it negative $U(n)$ structure}  or {\it negative Hermitian structure} on $V$ is a generalized metric with an orthogonal complex structure $\I_-$ on its $-1$-eigenspace. We say that a \gcs\ $\J$ {\it extends} a positive/negative $U(n)$ structure $(\G,\I)$  if $\I$ is the restriction of $\J$ to the appropriate space and $(\G,\J)$ is a generalized Hermitian structure.
\end{definition}


Given a generalized Hermitian structure, since $\J_1$ and $\G$ commute, $\J_1$ preserves the eigenspaces of $\G$ and hence, upon restriction to $V_{\pm}$, one obtains a positive and a negative Hermitian structure. Conversely, a positive (resp. negative) Hermitian structure can be extended to $\V$ by declaring that it vanishes on $V_-$ (resp. $V_+$). Then a pair of positive and negative Hermitian structures, $\I_+$, $\I_-$  gives rise to a generalized Hermitian structure by declaring that $\J_1 = \I_+  + \I_-$. 

Given a positive $U(n)$ structure on $V$, we can use the isomorphism $V_+ \cong V$ to transport the metric and the complex structure $\I_+$ from $V_+$ to $V$, making it into a Hermitian vector space $(V,g,I)$. Further, we can use $\I_+$ to define a complex structure $\I_-$  on $V_-$ using the isomorphisms  $V_+\cong V \cong V_-$ and this way we have an extension of the $U(n)$ structure to a generalized Hermitian structure: namely we declare that $\J_1$ is $\I_+$ on $V_+$ and $\I_-$ on $V_-$, hence, in the metric splitting of $\V$, $\J_1$ is the \gcs\ associated to  the complex structure\footnote{Note that usually, the \gcs\ associated to a given a complex structure $I$ on a vector space differs from $\J_1$ by a sign, or said another way, by overall conjugation. Hence, all the usual concepts get conjugated with this definition, for example, later we will have $\del_{\J_1} = \delbar$ and the generalized canonical bundle is $\wedge^{0,n}T^*M$}  $I$ and consequently $\J_2$ is the \gcs\ associated to the Hermitian form $\omega = g \circ I$:
\begin{equation}\label{eq:J1 and J2}
\J_1 = \begin{pmatrix}
I & 0\\
0 & -I^*
\end{pmatrix} \qquad
\J_2 = \begin{pmatrix}
0 & -\omega^{-1}\\
\omega & 0
\end{pmatrix}.
\end{equation}

For this set of choices, there is a relation between the $(p,q)$-decomposition of forms determined by the generalized Hermitian structure and the usual $(p,q)$-decomposition of forms determined by the complex structure $I$ on $V$. 

\begin{proposition}\label{prop:kahlerupq}
{\em (Cavalcanti \cite{MR2265463})} Let $\J_1$ and $\J_2$ be as above and let  $L_2$ be the $+i$-eingenspace of $\J_2$ . The map
$$\psi:\V_\C \into \V_\C;\qquad \psi(X +\xi) = e^{i\omega}_*e^{\frac{i}{2}\omega^{-1}}_*(X+ \xi)$$
preserves $\J_1$ and
$$\psi(V_\C)  = L_2,\qquad \psi(V^*_\C) = \bar{L_2}.$$
Therefore the corresponding action on spinors, 
$$\Psi:\wedge^{\bullet}V^*_{\C} \into \wedge^{\bullet}V^*_{\C};\qquad \Psi(\gf)= e^{i\omega}e^{\frac{i}{2}\omega^{-1}}(\gf),
$$
preserves the eigenspaces of $\J_1$ and maps $\wedge^kV^*_\C$ into  $U^{n-k}_{\J_2}$ that is
$$
\Psi(\wedge^{p,q}V^*) = U^{q-p,n-p-q}_{\J_1,\J_2}.$$
\end{proposition}

A positive Hermitian structure is fully determined by  the $+i$-eigenspace of $\I_+$, that is, an isotropic $n$-dimensional subspace $V^{1,0}_+\subset \T_\C M$ for which $V^{1,0}_+\cap \overline{V^{1,0}_+} =\{0\}$ and  $\IP{v,\overline{v}} > 0$ for all $v\in V^{1,0}_+\backslash\{0\}$.

If we are given a generalized complex  extension $\J$ of a positive Hermitian structure $(\G,\I_+)$ we obtain a bigrading of forms into $U^{p,q}$ as explained earlier. However, from the point of view of the $U(n)$ structure, the natural spaces to consider are  $W_+^k = \sum_{p+q=k}U^{p,q}$. Indeed $W_+^{n}$ corresponds to the vector space of all forms  that are annihilated by the Clifford action of $V_+^{1,0}$ and $W_+^{n-2k} = \wedge^{k} V^{0,1}_+ \cdot W_+^n$, that is, the spaces $W_+^k$ are solely determined by the complex structure $\I_+$. Notice that $W_+^k$ is only nontrivial if $-n\leq k\leq n$ and $n-k =0$ mod $2$. Further, since the spaces $W_+^k$ are the diagonals of the $U^{p,q}$ decomposition, each $W_+^k$ is made of either \gsd- or \gasd-forms, with $W_+^n$ being \gsd.

Another description of the spaces $W_+^\bullet$ is obtained by extending $\I_+$ to an endomorphism of $\TM$  by declaring that $\I_+|_{V_-}$ vanishes, so that   $\I_+$ is a skew-symmetric operator on $\TM$ \wrt\ the natural pairing, that is, $\I_+ \in \frak{spin}(\TM)$. Similarly to  a  generalized complex structure, $\J$, for which the space $U^k$ is the $ik$-eigenspace of the action of $\J$,  letting $\I_+$ act  on forms one sees that the space $W_+^k$ is the $i\frac{k}{2}$-eigenspace of $\I_+$. Hence, Lemma \ref{lem:star} takes the following form for positive $U(n)$ structures:

\begin{lemma}\label{lem:starSKT}
Let $(\G,\I_+)$ be a positive Hermitian structure on $V$, and let $\mathbb{I_+} = e^{\frac{\pi \I_+}{2}}$. Then
$$\mathbb{I}_+^2 = -\star,$$
that is $\star|_{W_+^k} = -i^k$. 
\end{lemma}

Similarly, for a negative Hermitian structure $(\G,\I_-)$, we can extend $\I_-$ to $\V$ by declaring that it vanishes on $V_+$. If $\J_1$ is a \gc\ extension of $\I_-$, then the eigenspaces of $\I_-$ correspond to the anti-diagonals of the $U^{p,q}$ decomposition: $W_-^k = \sum_{p-q=k} U^{p,q}$. 

Finally we observe that the spaces $W_\pm^k$  have both even and odd forms, so one can refine this grading to a $\Z \times \Z_2$-grading:
\begin{equation}\label{eq:W double graded}
W_\pm^{k,0} = W_\pm^k \cap (\wedge^{ev}V^*\tensor \C),\qquad W_\pm^{k,1} = W_\pm^k \cap (\wedge^{od}V^*\tensor \C).
\end{equation}

In what follows we will refer to both spaces $W_\pm^k$ and $W_\pm^{k,l}$, with the understanding that if the $\Z_2$-grading is not particularly important, we will simply omit it.

\section{Intrinsic torsion of generalized Hermitian structures}\label{sec:intrinsic torsion}

Except for a generalized metric, each of the structures introduced in Section \ref{sec:linear algebra} has an appropriate integrability condition. We let $(M^{2n},H)$ be a manifold  with a real closed 3-form $H$ and consider the Courant  bracket on sections of $\T M = TM \oplus T^*M$:
$$\Cour{X+\xi,Y +\eta}_H = [X,Y]+ \mc{L}_X \eta -i_Y d\xi - i_Y i_XH.$$
We will omit the 3-form from the bracket if it is clear from the context.

The Courant bracket is the derived bracket associated to the operator $d^H = d + H\wedge$, i.e., the following identity holds for all $v_1,v_2  \in \Gamma(\T M)$ and  $\gf \in \Omega^{\bullet}(M)$:
\begin{equation}
\begin{aligned}\label{eq:derived bracket}
\Cour{v_1,v_2}_H \cdot \gf& = \{\{v_1,d^H\},v_2\}\gf\\
& = d^H (v_1\cdot v_2 \cdot \gf) + v_1\cdot d^H(v_2\cdot \gf) - v_2\cdot d^H(v_1 \cdot \gf) - v_2\cdot v_1 \cdot d^H\gf,
\end{aligned}
\end{equation}
where $\cdot$ denotes the Clifford action of $\mathrm{Clif}(\T M)$ on $\wedge^{\bullet}T^*M$ and $\{\cdot,\cdot\}$ denotes the graded commutator of operators.

The orthogonal action of a 2-form $B \in \Omega^2(M)$ on $\T M$ relates different Courant brackets:
$$\Cour{e^B_*v_1,e^B_* v_2}_H = e^B_* \Cour{v_1,v_2}_{H+dB}.$$

\begin{definition}
 For each of the structures introduced in Definitions \ref{def:linear gcs} and \ref{def:linear SKT}, we refer to  the smooth assignment of such structure to $T_xM$ for each $x\in M$ by including the adjective {\it almost} in  the name of the structure.
\end{definition}

\begin{definition}[Integrability conditions]\label{def:integrability}~
\begin{itemize}
\item A {\it  (integrable) \gcs}  is a \gacs\ $\J$ whose  $+i$-eigenspace is involutive \wrt\ the Courant bracket.
\item A {\it generalized Hermitian structure} is a pair $(\G,\J_1)$ of generalized metric and compatible integrable generalized complex structure.
\item A {\it \gks} is a generalized Hermitian structure $(\G,\J_1)$ for which $\J_1$  and $\J_2 = \J_1\G$ are integrable.
\end{itemize}
\end{definition}

\subsection{The Nijenhuis tensor}\label{subsec:nijenhuis}

Let us spend some time to understand the {\it Nijenhuis tensor} of an almost  \gcs\ $\J$. 
This tensor is defined in the usual way, namely if $L$ is the $+i$-eigenspace of $\J$
\begin{equation}\label{eq:Nij}
\mathrm{Nij}:\Gamma(\bar{L}) \times \Gamma(\bar{L}) \into \Gamma(L);\qquad \mathrm{Nij}(X,Y) = -\Cour{X,Y}^{L},
\end{equation}
where $\cdot^{L}$ indicates projection onto $L$. We can alternatively use the identification $L = \bar{L}^*$ from \eqref{eq:Lbar and L*} to consider the operator
\begin{equation}\label{eq:def nijenhuis}
\begin{aligned}N: \Gamma(\bar{L}) \times \Gamma(\bar{L}) \times \Gamma(\bar{L}) &\into \Omega^0(M;\C);\\
N(X,Y,Z) = -2\IP{\Cour{X,Y},Z} &= 2\IP{\mathrm{Nij}(X,Y),Z}.
\end{aligned}
\end{equation}
As usual, $\mathrm{Nij}$  is a tensor, indeed, for $f \in \C^{\infty}(M;\C)$ we have
$$\mathrm{Nij}(X,fY) = -\Cour{X,fY}^{L} = -(f\Cour{X,Y} + (\mc{L}_{\pi_T(X)}f )Y)^{L} = -f\Cour{X,Y}^{L} = f \mathrm{Nij}(X,Y).$$
Further, the tensor $N \in \wedge^2L\tensor L$ defined above actually lies in $\wedge^3 L$. Indeed, for $X,Y,Z\in \Gamma(\bar{L})$ we have
$$0= \mc{L}_{\pi_T(X)}\IP{Y,Z} = -(\IP{\Cour{X,Y},Z} +\IP{Y,\Cour{X,Z}}) = \tfrac{1}{2}(N(X,Y,Z)  + N(X,Z,Y)),$$
which shows that $N$ is fully skew.

A different way to understand $N$ arises by using the $U^k$ decomposition of forms determined by $\J$. Namely, letting $\U^k = \Gamma(U^k)$ (throughout the paper we denote the sheaf of sections of $U^k$ by $\U^k$ and the sheaf of sections of $U^{p,q}$ by $\U^{p,q}$), one has:
\begin{lemma}\label{lem:nijenhuis}
 In an almost  generalized complex manifold
\begin{equation}\label{eq:deldelbar}
d^H:\U^k\into \U^{k-3} \oplus \U^{k-1} \oplus \U^{k+1} \oplus \U^{k+3}.
\end{equation}
Further the map $\pi_{k+3} \circ d^H:\U^k \into \U^{k+3}$ corresponds to the Clifford action of $N$, the Nijenhuis tensor defined in \eqref{eq:def nijenhuis}.
\end{lemma}
\begin{proof}
We prove first that
\begin{equation}\label{eq:Nijenhuis step 1}
d^H:\U^k \into \sum_{j\geq k-3} \U^j
\end{equation}
and will do so by induction, starting at $k=n+1$ and working our way down. For $k = n+1$ we have that  $U^{n+1}=\{0\}$ and the claim follows trivially.

Next we assume the result to be true for all $j >k$ and let $\rho \in \U^k$. For $v_1,v_2 \in  \Gamma(L)$, using \eqref{eq:derived bracket} we have
$$v_2\cdot v_1 \cdot d^H\rho  = -\Cour{v_1,v_2}\cdot \rho + d^H (v_1\cdot v_2 \cdot \rho) + v_1 d^H (v_2 \cdot \rho) -v_2 d^H (v_1 \cdot \rho).$$
Since the Clifford action of $v_i$ sends $\U^j$ to $\U^{j+1}$, the inductive hypothesis implies that the last three terms lie in $\oplus_{j\geq k-1}\U^{j}$, while  the first term, being the action of an element of $L \oplus \bar{L}$ on $\rho$ lies in $\U^{k+1}\oplus \U^{k-1}$. Therefore we conclude that 
$$v_2\cdot v_1 \cdot d^H\rho \in \oplus_{j\geq k-1} \U^{j}\qquad \mbox{for all } v_1,v_2 \in \Gamma(L)$$
and \eqref{eq:Nijenhuis step 1} follows.

As $d^H$ is a real operator and $U^{-k} = \bar{U^k}$ conjugating \eqref{eq:Nijenhuis step 1} we have that
$$d:\U^k \into \sum_{j\leq k-3} \U^j.$$
Furthermore, if $U^k$ is made of even forms then $U^{k+1}$ is made of odd forms and vice versa, we have that
$$d^H:\U^k \into \sum_{\substack{|j-k |\leq3\\ j-k =1\mbox{ mod } 2}} \U^j,$$
therefore proving \eqref{eq:deldelbar}.

Now we prove that $\pi_{k+3} \circ d^H$ corresponds to the Clifford action of $N$. Once again we use induction, this time starting at $k=-n-1$ and moving upwards. Since $\U^{-n-1} = \{0\}$, the claim is trivial there. Assume now that for $j < k$ we have proved that $\pi_{j+3} \circ d^H$ is the Clifford action of $N$. Let $\rho \in \U^k$ and $v_1,v_2 \in \Gamma(\bar{L})$. Then using \eqref{eq:derived bracket} we have
\begin{align*}
v_2 v_1&  \pi_{k+3}\circ d^H \rho =  \pi_{k+1}(v_2 v_1  d^H \rho)\\
&=\pi_{k+1}(-\Cour{v_1,v_2} \rho + d^H v_1  v_2  \rho + v_1  d^H v_2 \rho - v_2  d^H v_1 \rho)\\
&= N(v_1,v_2) \rho + \pi_{k+1}\circ d^H  (v_1  v_2  \rho ) + v_1\pi_{k+2}\circ d^H v_2 \rho -v_2 \pi_{k+2}\circ d^H v_1  \rho\\
&= \{v_2,\{v_1,N\}\} \rho +N v_1  v_2 \rho  + v_1  N  v_2  \rho - v_2  N  v_1  \rho\\
&=v_2 v_1 N  \rho,
\end{align*}
where in third equality we have used that the component of $\Cour{v_1,v_2}\rho$ in $\U^{k+1}$ is given by the Clifford action of the $L$ component of  $\Cour{v_1,v_2}$, that is $\iota_{v_2} \iota_{v_1}N$ and in the fourth equality we used the inductive hypothesis as well as the fact that when acting on forms, the interior product of $v \in \Gamma(\bar{L})$ with $\gf \in \Gamma(\wedge^{\bullet} L)$ is given by $\gf(v) \rho = \{v,\gf \}\rho$. The last equality follows by expanding the graded commutator and canceling out similar terms. 
 \end{proof}

If we compose $d^H|_{\U^k}$ with projection onto $\U^{k+1}$ and $\U^{k-1}$ we get operators
$$\del_{\J} = \pi_{k+1}\circ d^H :\U^k \into \U^{k+1}\qquad\mbox{and}\qquad \delbar_{\J} =\pi_{k-1}\circ d^H:\U^k\into \U^{k-1},$$
and if the \gcs\ $\J$ is clear from the context, we denote these operators simply  by $\del$ and $\delbar$. As proved by Gualtieri in \cite{MR2811595}, integrability is equivalent to the requirement
$$d^H:\U^k\into  \U^{k-1} \oplus \U^{k+1};\qquad d^H = \del+\delbar,$$
so we can also see from this point of view how the vanishing of the Nijenhuis tensor implies integrability.


\subsection{The intrinsic torsion and the road to integrability}\label{subsec:intrinsic}

\begin{figure}[h!!]
\begin{center}
$$\xymatrix@R=18pt@C=12pt{
\U^{p-3,q+3}& &\U^{p-1,q+3} &  &\U^{p+1,q+3} &&\U^{p+3,q+3}\\
\U^{p-3,q+1}& &\U^{p-1,q+1} &  &\U^{p+1,q+1} &&\U^{p+3,q+1}\\
&&                &\ar@/ ^0.5pc/[llluu]_(0.7){\bar{N_-}} \ar@/ _0.5pc/[luu]_{N_1} \ar[lllu]^{\bar{N_4}} \ar[lu]_{\deltabar_-}  \ar[llld]_{\bar{N_3}} \ar[ld]^{\deltabar_+}\ar@/ _0.5pc/[llldd]^(0.7){\bar{N_+}} \ar@/ ^0.5pc/[ldd]^{\bar{N_2}} \U^{p,q} \ar@/ _0.5pc/[rdd]_{\bar{N_1}} \ar@/ ^0.5pc/[rrrdd]_(0.7){N_-}\ar[rd]_{\delta_-} \ar[rrrd]^{N_4} \ar[ru]^{\delta_+}\ar[rrru]_{N_3}\ar@/ ^0.5pc/[ruu]^{N_2} \ar@/ _0.5pc/[rrruu]^(0.7){N_+}             &&& \\
\U^{p-3,q-1}& &\U^{p-1,q-1} &  &\U^{p+1,q-1} &&\U^{p+3,q-1}\\
\U^{p-3,q-3}& &\U^{p-1,q-3} &  &\U^{p+1,q-3} &&\U^{p+3,q-3}\\
}
$$
\caption{Representation of the nontrivial components of $d^H$ when restricted to $\U^{p,q}$ for an generalized almost Hermitian structure.}\label{fig:nonintegrable}
\end{center}
\end{figure}
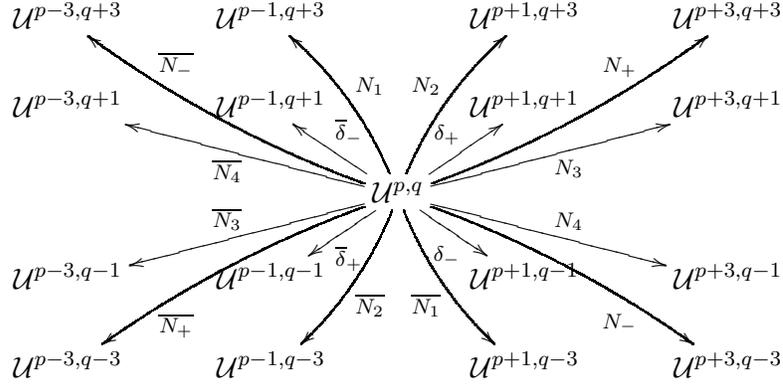

With this understanding of the Nijenhuis tensor, we can give a pictorial description of  the long road to integrability from almost generalized Hermitian to \gk. Indeed, given an almost generalized Hermitian structure, we get a splitting of forms into spaces $U^{p,q}$. According to Lemma \ref{lem:nijenhuis},  $d^H$ can not change either the `$p$' or  the `$q$' grading by more than three and it must switch parity. Hence $d^H$ decomposes as a sum of  eight operators and their complex conjugates
$$d^H =  \delta_++\deltabar_+ + \delta_-+\deltabar_- + N_+ +\bar{N_+}+ N_- + \bar{N_-} + N_1 +\bar{N_1}+ N_2 + \bar{N_2} + N_3 +\bar{N_3}+  N_4 +\bar{N_4};$$
\begin{equation}\label{eq:operators}
\begin{aligned}
\delta_+ :\U^{p,q} \into \U^{p+1,q+1},&\qquad \delta_-:\U^{p,q} \into \U^{p+1,q-1},\\
N_+ :\U^{p,q} \into \U^{p+3,q+3},&\qquad N_-:\U^{p,q} \into \U^{p+3,q-3},\\
N_1 :\U^{p,q} \into \U^{p-1,q+3},&\qquad N_2:\U^{p,q} \into \U^{p+1,q+3},\\
N_3 :\U^{p,q} \into \U^{p+3,q+1},&\qquad N_4:\U^{p,q} \into \U^{p+3,q-1}.
\end{aligned}
\end{equation}
and we can draw in a diagram all the possible nontrivial components of $d^H|_{\U^{p,q}}$ as arrows (see Figure \ref{fig:nonintegrable}).

\begin{figure}[h!!]
\begin{center}
$$\xymatrix@R=12pt@C=12pt{
& &\U^{p-1,q+3} &  &\U^{p+1,q+3} &&\\
& &\U^{p-1,q+1} &  &\U^{p+1,q+1} &&\\
&&                & \ar[luu]_(0.6){N_1} \ar[lu]^{\deltabar_-}   \ar[ld]_{\deltabar_+} \ar[ldd]^(0.6){\bar{N_2}} \U^{p,q} \ar[rdd]_(0.6){\bar{N_1}} \ar[rd]^{\delta_-}\ar[ru]_{\delta_+}\ar[ruu]^(0.6){N_2}              &&& \\
& &\U^{p-1,q-1} &  &\U^{p+1,q-1} &&\\
& &\U^{p-1,q-3} &  &\U^{p+1,q-3} &&\\
}$$
\caption{Representation of the nontrivial components of $d^H$ for a generalized Hermitian structure.}\label{fig:genhermitian}
\end{center}
\end{figure}
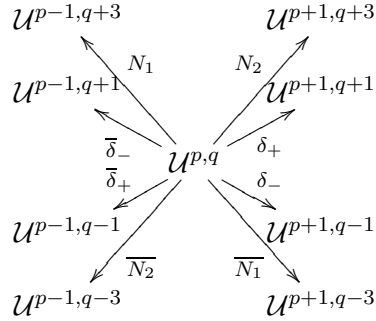

Once we require that $\J_1$ is integrable, i.e. we are dealing in fact with a generalized Hermitian structure, then $d^H$ only changes the `$p$' degree by $\pm 1$ and several components of $d^H$ present in the nonintegrable case, now vanish and the diagram from Figure \ref{fig:nonintegrable} clears up to the one presented in Figure \ref{fig:genhermitian}.

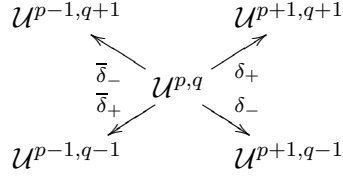
\begin{figure}[h!!]
\begin{center}
$$\xymatrix@R=12pt@C=6pt{
& &\U^{p-1,q+1} &  &\U^{p+1,q+1} &&\\
&&                & \ar[lu]^{\deltabar_-}   \ar[ld]_{\deltabar_+}  \U^{p,q} \ar[rd]^{\delta_-}\ar[ru]_{\delta_+}              &&& \\
& &\U^{p-1,q-1} &  &\U^{p+1,q-1} &&\\
}
$$
\caption{Representation of the nontrivial components of $d^H$ for a \gks.}\label{fig:gen kahler}
\end{center}
\end{figure}

Finally, if we require that $\J_2$ is also integrable, and hence we are in fact dealing with a \gks,  the last two components of the Nijenhuis tensor, labeled $N_1$ and $N_2$ above, vanish and $d^H$ decomposes as a sum of four operators, as pictured in Figure \ref{fig:gen kahler}.

This shows that the obstruction for a generalized almost Hermitian structure to be a \gks\ is given by the tensors $N_{i}$, $i=1,2,3,4$ and $N_{\pm}$

\begin{definition}
The {\it intrinsic torsion} of a generalized Hermitian manifold are the tensors $N_i$, $i=1,2,3,4$ and $N_\pm$. 
\end{definition}

Next sections we will introduce geometric structures on $M$ which are weaker than \gkss\ and show how these structures can be phrased in terms of the vanishing of certain components of the intrinsic torsion.

\subsection{The operators $\delta_\pm$ and $\deltabar_\pm$}
Not requiring integrability, $d^H$  restricted to $\U^{p,q}$ has sixteen components but only four are not tensorial, namely $\delta_\pm$ and $\deltabar_\pm$. Hence we have
\begin{equation}\label{eq:decomposition of dH}
d^H \sim \delta_+ + \deltabar_+ + \delta_-+\deltabar_-,
\end{equation}
where $\sim$ indicates that these operators agree up to lower order terms, i.e., they have the same symbol. Using the decomposition 
\begin{equation}\label{eq:decomposition of TM}
\T_\C M = V_+^{1,0}\oplus V_+^{0,1} \oplus V_-^{1,0} \oplus V_-^{0,1},
\end{equation}
one can easily see that the decomposition \eqref{eq:decomposition of dH} in terms of symbols corresponds to the decomposition of a 1-form $\xi \in \TM$ into its four components according to \eqref{eq:decomposition of TM} and one can also check that the symbol sequence for each of the operators $\delta_\pm$ and $\deltabar_\pm$ is exact, e.g., for $\delta_+$, the symbol sequence associated to sequence of operators
$$
\cdots \stackrel{\delta_+}{\into}\U^{p-1,q-1}\stackrel{\delta_+}{\into}\U^{p,q}\stackrel{\delta_+}{\into}\U^{p+1,q+1}\stackrel{\delta_+}{\into}\cdots
$$
is an exact sequence.
Adding over $p+q =k$ and letting $\W^k_+ = \oplus_{p+q =k} \U^{p,q} $ we have also that the symbol sequence associated to
\begin{equation}\label{eq:wk seq}
\cdots \stackrel{\delta_+}{\into}\W_+^{k-2}\stackrel{\delta_+}{\into}\W_+^{k}\stackrel{\delta_+}{\into}\W_+^{k+2}\stackrel{\delta_+}{\into}\cdots
\end{equation}
is exact.
And similarly, adding over $q$ we get that the symbol sequence of 
\begin{equation}\label{eq:up seq}
\cdots \stackrel{\delta_+}{\into}\U^{p-1}\stackrel{\delta_+}{\into}\U^{p}\stackrel{\delta_+}{\into}\U^{p+1}\stackrel{\delta_+}{\into}\cdots
\end{equation}
is exact.

\section{Parallel Hermitian and bi-Hermitian structures}\label{sec:parallel}

The first type of structures that we will relate to the intrinsic torsion are Hermitian structures which are parallel for a connection with closed skew torsion. The existence of a relationship between connections with closed, skew symmetric torsion and the Courant bracket was made evident by Hitchin in \cite{MR2253158}. Precisely, given a generalized metric on $\TM$, we let $H$ be the 3-form corresponding to the metric splitting and $g$ be the induced metric on $M$. Also, for $X \in \Gamma(TM)$ we let $X_\pm \in \Gamma(V_{\pm})$ be unique lifts of $X$ to $V_\pm$, we let $\pi_\pm:\TM \into V_{\pm}$ be the orthogonal projections onto $V_\pm$ and $\pi_T:\TM \into TM$ be the natural projection.

\begin{proposition}[Hitchin \cite{MR2253158}]\label{prop:connection}
Let $\tilde{\nabla}^{\pm}$ be the unique metric connection whose torsion is skew symmetric and equal to $\mp H$. Then
$$\tilde{\nabla}^+_X Y = \pi_T\circ  \pi_+\Cour{X_-,Y_+},$$
$$\tilde{\nabla}^-_X Y = \pi_T \circ \pi_-\Cour{X_+,Y_-}.$$
\end{proposition}

From this proposition, we see that the isomorphisms $\pi_T: V_{\pm} \into TM$ relate the connections with torsion $\mp H$ to the operators
\begin{IEEEeqnarray*} {l l l l}
\nabla^+&:\Gamma(V_-)\times \Gamma(V_+)\into \Gamma(V_+);&\qquad  \nabla^+_v w &= \pi_+\Cour{v,w}, \quad v\in \Gamma(V_+), w \in \Gamma(V_-),\IEEEyesnumber\label{eq:nabla+}\\
\nabla^-&:\Gamma(V_+)\times \Gamma(V_-)\into \Gamma(V_-);&\qquad  \nabla^-_w v &= \pi_-\Cour{w,v}, \quad v\in \Gamma(V_+), w \in \Gamma(V_-).\IEEEyesnumber\label{eq:nabla-}
\end{IEEEeqnarray*}
As we will work with the spaces $V_\pm$ directly, we will use  $\nabla^\pm$ instead of the connections they 
induce on $TM$, with the understanding that these are equivalent operators, so, for example,  if $\tilde\nabla^+$ has holonomy in $U(n)$, then $M$ has an almost Hermitian structure $(g,I)$ which is parallel \wrt\ $\tilde{\nabla}^+$ and via the isomorphism $\pi_T:V_+ \into TM$, $V_+$ gets a positive almost Hermitian structure  which is parallel for $\nabla^+$ so that $\nabla^+$ has holonomy in $U(n)$ and the converse statement also holds.

\begin{definition}
A {\it parallel positive (resp. negative) Hermitian structure} is a positive (resp. negative)  almost Hermitian structure which is parallel \wrt\ $\nabla^+$, (resp. $\nabla^-$). A {\it parallel bi-Hermitian structure} or {\it parallel $U(n)\times U(n)$-structure} is a triple $(g,I_+,I_-)$ such that $(g,I_+)$ is a parallel positive Hermitian structure and $(g,I_-)$ is a parallel negative Hermitian structure. 
\end{definition}

\begin{proposition}\label{prop:U(n) holonomy}
If the connection $\nabla^+$ has holonomy in $U(n)$, the corresponding almost Hermitian structure satisfies
\begin{equation}\label{eq:positive}
\Cour{\Gamma(V_+^{1,0}),\Gamma(V_+^{1,0})} \subset \Gamma(V_+\tensor\C).
\end{equation}
And conversely, a Hermitian structure satisfying \eqref{eq:positive} is parallel \wrt\ $\nabla^+$.
\end{proposition}
\begin{proof}
Let $\I_+$ be a complex structure on $V_+$ orthogonal \wrt\ the natural pairing and let  $V_+^{1,0}$ be its $+i$-eigenspace. Then, for $v_1, v_2 \in \Gamma(V_+^{1,0})$ and $w \in \Gamma(V_-)$, we have that $\IP{v_i,w} =0$, as $V_+$ and $V_-$ are orthogonal \wrt\ the natural pairing. Hence
\begin{equation}\label{eq:reduced holonomy computation}
\begin{aligned}
0&= \L_{\pi_T v_1}\IP{w,v_2}\\
&= \IP{\Cour{v_1,w},v_2} + \IP{w,\Cour{v_1,v_2}}\\
&= -\IP{\Cour{w,v_1},v_2} + \IP{w,\Cour{v_1,v_2}}\\
& =  - \IP{\nabla^+_w v_1,v_2} + \IP{w,\Cour{v_1,v_2}}.
\end{aligned}
\end{equation}
where in the third equality we used again that $V_+$ and $V_-$ are orthogonal \wrt\ the natural pairing, and hence $\Cour{v_1,w} = -\Cour{w,v_1}$ and in the fourth equality we used that $v_2 \in V_+$ and hence the $V_-$ component of $\Cour{w,v_1}$ is annihilated by the natural pairing.

If $\I_+$ is parallel \wrt\ $\nabla^+$, then $\IP{\nabla^+_w v_1,v_2} $ vanishes for all $v_1, v_2 \in \Gamma(V_+^{1,0})$ and $w \in \Gamma(V_-)$ and hence, according to \eqref{eq:reduced holonomy computation}, so does $\IP{w,\Cour{v_1,v_2}}$, showing that  $\Cour{v_1,v_2}$ must be orthogonal to $V_-$ and hence, a section of $V_+$. Conversely, if $,\Cour{v_1,v_2} \in \Gamma(V_+)$, then  for all $w \in \Gamma(V_-)$, $\IP{w,\Cour{v_1,v_2}}=0$ hence \eqref{eq:reduced holonomy computation} implies that $\IP{\nabla^+_w v_1,v_2}=0$ showing that $\nabla^+_w v_1$ is orthogonal to $V_+^{1,0}$. Since $V_+^{1,0}$ is a maximal isotropic of $V_+\tensor\C$, we conclude that $\nabla^+_w v_1 \in \Gamma(V_+^{1,0})$ for all $w \in V_-$ and $v_1 \in \Gamma(V_+^{1,0})$, hence $\I_+$ is parallel \wrt\ $\nabla^+$.
\end{proof}

Obviously the same result holds exchanging $\nabla^+$ and $V_+^{1,0}$ by $\nabla^-$ and $V_-^{1,0}$.

Next, and throughout these notes, we denote by $\mc{W}^k$ and $\W^{k,l}$ the sheaf of sections of the bundles $W^k$ and $W^{k,l}$ respectively. 

\begin{theorem}\label{prop:parallel via nijenhuis}
Let $(M,H)$ be a manifold with 3-form, $\G$ be a generalized metric on $M$ and $\I_+$, $\I_-$ be a positive and a negative almost Hermitian structure on $M$
\begin{itemize}
\item $\I_+$ is parallel \wrt\ $\nabla^+$ \iff\ $N_2 = N_3 =0$, i.e.,
$$d^H:\W^k_+ \into \W_+^{k-6} \oplus \W_+^{k-2} \oplus \W_+^k \oplus \W_+^{k+2}\oplus \W_+^{k+6};$$
\item $\I_-$ is parallel \wrt\ $\nabla^-$ \iff\ $N_1 = N_4 =0$, i.e.,
$$d^H:\W^k_- \into \W_-^{k-6} \oplus \W_-^{k-2} \oplus \W_-^k \oplus \W_-^{k+2}\oplus \W_-^{k+6};$$
\item $(\G,\I_+,\I_-)$ is a parallel bi-Hermitian structure \iff\ $N_i=0$ for $i=1,2,3,4$.
\end{itemize}
\end{theorem}
\begin{proof}
Since $V_-$ is the orthogonal complement of $V_+$ \wrt\ the natural pairing, it is clear that \eqref{eq:positive} is equivalent to the following two conditions
$$\IP{\Cour{v_1,v_2},w} = 0 \mbox{ for all }  v_1,v_2 \in \Gamma(V_+^{1,0}), w \in \Gamma(V_-^{0,1}).$$
$$\IP{\Cour{v_1,v_2},w} = 0 \mbox{ for all }  v_1,v_2 \in \Gamma(V_+^{1,0}), w \in \Gamma(V_-^{1,0});$$
but the first condition is equivalent to the vanishing of $\bar{N_2}$ and the second, to the vanishing of $\bar{N_3}$. Finally, since the component of $d^H$ mapping $\W^k$ to $\W^{k+4}$ is given by the sum $N_2 + N_3$, we see that the vanishing of this component is also equivalent to the parallel condition.

The remaining claims are proved similarly.
\end{proof}

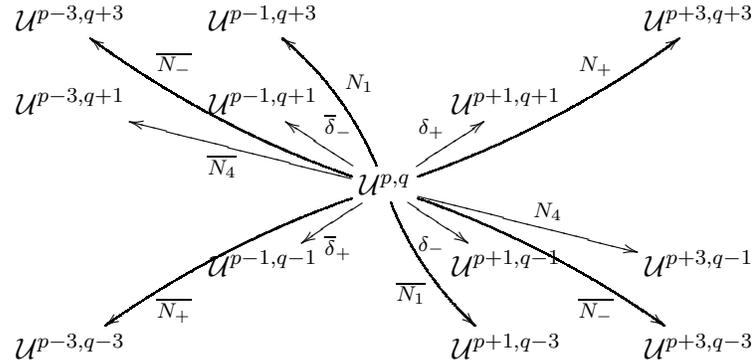
\begin{figure}[h!!]
\begin{center}
$$\xymatrix@R=16pt@C=10pt{
\U^{p-3,q+3}& &\U^{p-1,q+3} &  & &&\U^{p+3,q+3}\\
\U^{p-3,q+1}& &\U^{p-1,q+1} &  &\U^{p+1,q+1} &&\\
&&                &\ar@/ ^0.5pc/[llluu]_(0.7){\bar{N_-}} \ar@/ _0.5pc/[luu]_{N_1} \ar[lllu]^{\bar{N_4}} \ar[lu]_{\deltabar_-}   \ar[ld]^{\deltabar_+}\ar@/ _0.5pc/[llldd]^(0.7){\bar{N_+}}  \U^{p,q} \ar@/ _0.5pc/[rdd]_{\bar{N_1}} \ar@/ ^0.5pc/[rrrdd]_(0.7){\bar{N_-}}\ar[rd]_{\delta_-} \ar[rrrd]^{N_4} \ar[ru]^{\delta_+} \ar@/ _0.5pc/[rrruu]^(0.7){N_+}             &&& \\
& &\U^{p-1,q-1} &  &\U^{p+1,q-1} &&\U^{p+3,q-1}\\
\U^{p-3,q-3}& & &  &\U^{p+1,q-3} &&\U^{p+3,q-3}\\
}
$$
\caption{Representation of the nontrivial components of $d^H$ for an almost  generalized Hermitian structure extending a parallel positive Hermitian structure}\label{fig:type1}
\end{center}
\end{figure}

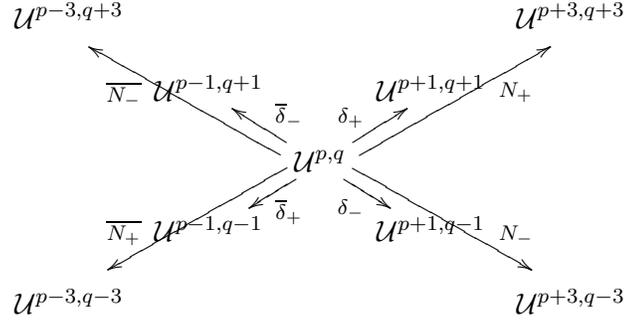
\begin{figure}[h!!]
\begin{center}
$$\xymatrix@R=12pt@C=6pt{
\U^{p-3,q+3} &  & &&\U^{p+3,q+3}\\
& \U^{p-1,q+1} &  &\U^{p+1,q+1} &\\
&                &\ar@<1ex>[lluu]^(0.7){\bar{N_-}}  \ar[lu]_(0.4){\deltabar_-}   \ar[ld]^(0.4){\deltabar_+}\ar@<-1ex>[lldd]_(0.7){\bar{N_+}}  \U^{p,q} \ar@<1ex>[rrdd]^(0.7){N_-}\ar[rd]_(0.4){\delta_-} \ar[ru]^(0.4){\delta_+} \ar@<-1ex>[rruu]_(0.7){N_+}             && \\
 &\U^{p-1,q-1} &  &\U^{p+1,q-1} &\\
\U^{p-3,q-3}& & &  & \U^{p+3,q-3}\\
}
$$\caption{Representation of the nontrivial components of $d^H$  for a parallel bi-Hermitian structure.}\label{fig:structureX}
\end{center}
\end{figure}

The decomposition of $d^H$ for an almost generalized Hermitian structure extending a parallel positive Hermitian structure is depicted in Figure \ref{fig:type1} and the decomposition of $d^H$ for a parallel bi-Hermitian structure is depicted in Figure \ref{fig:structureX}.

\section{SKT structures}\label{subsec:integrability}

Classically, an SKT structure is a Hermitian structure $(M,g,I)$  for which $dd^c \omega =0$, where $\omega$ is the Hermitian 2-form. In this case, $d^c \omega =H$ is a closed 3-form and the complex structure is in fact parallel \wrt\ the metric connection with torsion $-H$, hence SKT structures are a particular case of positive parallel Hermitian structures, the only difference being that now one requires $I$ to be integrable.  In this section we phrase the SKT condition in terms of generalized geometry and relate it with the intrinsic torsion. We will see that an SKT structure lies precisely half way between a generalized Hermitian and a \gks.

Given an SKT structure $(g,I)$ we let $H$ be the background 3-form and consider the generalized metric $\G$  as in \eqref{eq:generalized metric}. Using the isomorphism $V_+ \cong TM$ we use $I$ to induce a complex structure  $\I_+$ on $V_+$ and hence split $V_+$ into eigenspaces of $\I_+$: $V_+\tensor \C = V_+^{1,0} \oplus V_+^{0,1}$.

\begin{proposition}\label{prop:SKT to gc}{\em (Cavalcanti \cite{MR2314216})}
For an SKT structure  $(g,I)$ on $M$ with $d^c\omega =H$, $V_+^{1,0}$ is involutive \wrt\ the $H$-Courant bracket.

Conversely, given a positive almost Hermitian structure $(\G,\I_+)$ on $(M,H)$, where $H$ is the 3-form associated of the metric splitting of $\TM$, if $V_+^{1,0}$ is involutive, then the induced Hermitian structure $(g,I)$ on $M$ is an SKT structure with $d^c \omega= H$ .
\end{proposition}

\begin{remark}
Here we obtain, in a new light, a well known contrast between connections with torsion and the Levi--Civita connection regarding integrability. Indeed, for the Levi--Civita connection, reduction of the holonomy group to $U(n)$ implies integrability of the complex structure, but  that is known not to be the case for  connections with torsion. From our point of view, this is the difference between the reduced holonomy condition
$$\Cour{\Gamma(V_+^{1,0}),\Gamma(V_+^{1,0})} \subset \Gamma(V_+)$$
and the SKT condition
$$\Cour{\Gamma(V_+^{1,0}),\Gamma(V_+^{1,0})} \subset \Gamma(V_+^{1,0}).$$
\end{remark}

Of course, if we have $d^c \omega = -H$, we can lift the SKT data to $V_-$ to obtain a negative Hermitian structure and an analogous version of Proposition \ref{prop:SKT to gc} holds.

This result motivates our formulation of the SKT condition.
\begin{definition}
A {\it positive (resp. negative)  SKT structure} is a positive (resp. negative) almost $U(n)$ structure $(\G,\I)$ for which the $+i$-eigenspace of $\I$ is involutive.
A generalized (almost) complex structure $\J$ {\it extends} an SKT structure  if $\J$ is fiberwise an extension of $\I$, in which case we say that $\J$ is a {\it generalized (almost) complex/Hermitian extension} of the SKT structure.
\end{definition}

Firstly, we observe that the SKT condition can also be phrased in terms of the vanishing of components of the intrinsic torsion. Indeed, in the presence of an extension of, say, a positive $U(n)$ structure $(\G,\I_+)$ to a \gacs\ $\J_1$ if we  let $V_+^{1,0}$ be the $+i$-eigenspace of $\I_+$ then involutivity  is equivalent to
$$\IP{\Cour{v_1,v_2},w_1} = 0 \qquad\mbox{ for  all }v_1,v_2 \in \Gamma(V_+^{1,0}), w \in \Gamma(V_+^{1,0}\oplus V_-),$$
that is,  the components $N_2, N_3$ and $N_+$ of the intrinsic torsion vanish. A similar argument gives a characterization of negative SKT structures.

\begin{proposition}\label{prop:SKT from nijenhuis}
Let $\G$ be a generalized metric on  a manifold with 3-form $(M,H)$ and let $\I_\pm$ be a pair of positive and negative Hermitian structures compatible with $\G$. If  $\J_1 = \I_+ + \I_-$ then $(\G,\I_+)$ is an SKT structure \iff\  $N_2, N_3$ and $N_+$ vanish and $(\G,\I_-)$ is a negative SKT structure \iff\ $N_1,N_4$ and $N_-$ vanish.
\end{proposition}

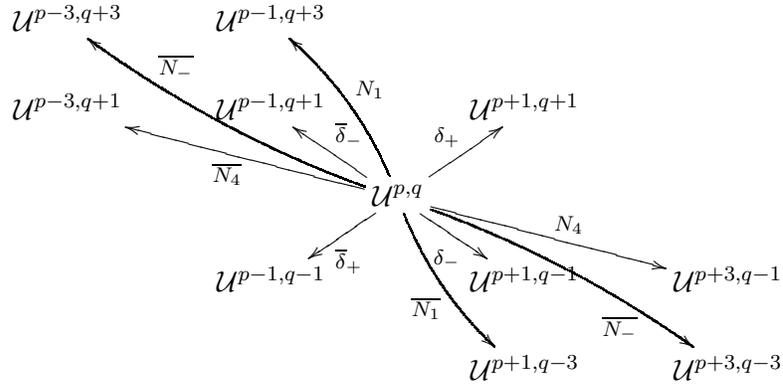
\begin{figure}[h!!]
\begin{center}
$$\xymatrix@R=18pt@C=12pt{
\U^{p-3,q+3}& &\U^{p-1,q+3} &  & &&\\
\U^{p-3,q+1}& &\U^{p-1,q+1} &  &\U^{p+1,q+1} &&\\
&&                &\ar@/ ^0.5pc/[llluu]_(0.7){\bar{N_-}} \ar@/ _0.5pc/[luu]_{N_1} \ar[lllu]^{\bar{N_4}} \ar[lu]_{\deltabar_-}   \ar[ld]^{\deltabar_+}\U^{p,q} \ar@/ _0.5pc/[rdd]_{\bar{N_1}} \ar@/ ^0.5pc/[rrrdd]_(0.7){\bar{N_-}}\ar[rd]_{\delta_-} \ar[rrrd]^{N_4} \ar[ru]^{\delta_+}             &&& \\
& &\U^{p-1,q-1} &  &\U^{p+1,q-1} &&\U^{p+3,q-1}\\
& & &  &\U^{p+1,q-3} &&\U^{p+3,q-3}\\
}
$$
\caption{Representation of the nontrivial components of $d^H$ for a positive SKT structure}\label{fig:+SKT}
\end{center}
\end{figure}

\begin{corollary}
Let $(M,H)$ be a manifold with 3-form.
\begin{itemize}
\item A parallel positive (resp. negative) almost  Hermitian structure is a positive (resp. negative) SKT structure \iff\  $N_+ =0$ (resp. $N_-=0$);
\item  A parallel bi-Hermitian structure is a \gks\ \iff\ $N_+ =N_-=0$.
\end{itemize}   
\end{corollary}

\begin{corollary}
Let $M$ be a four dimensional manifold.
\begin{itemize}
\item A parallel positive/negative  Hermitian structure  is a positive/negative SKT structure;
\item A parallel bi-Hermitian structure is a \gks.
\end{itemize}
\end{corollary}
\begin{proof}
Since $M$ is four dimensional, $V_\pm^{1,0}$ are two dimensional complex vector spaces, so  $\wedge^3 V_\pm^{1,0} = \{0\}$ and hence $N_\pm \in \wedge^3 V^{1,0}_\pm$ are the trivial tensors.
\end{proof}

\begin{corollary}[Gualtieri \cite{gualtieri-2010}]
Let $\G$ be a generalized metric on  a manifold with 3-form $(M,H)$, $\I_\pm$ be a pair of positive and negative SKT structures compatible with $\G$ and  $\J_1 = \I_+ + \I_-$. Then $(\G,\J_1)$ is a \gks.
\end{corollary}
\begin{proof}
According to Proposition \ref{prop:SKT from nijenhuis}, under the hypothesis, all components of the intrinsic torsion vanish hence $(\G,\J_1)$ is a \gks.
\end{proof}

Next, we describe the integrability condition for an SKT structure in terms of the decomposition of forms into $W_\pm^k$ and $U^{p,q}$ (for a fixed generalized complex extension $\J_1$)  described in the previous section.

\begin{theorem}\label{theo:skt}
Let $(\G,\I_+)$ be a positive almost Hermitian structure on a manifold with 3-form $(M^{2n},H)$. Then the following are equivalent:
\begin{enumerate}
\item $(\G,\I_+)$ is a positive SKT structure;
\item $d^H:\W_+^k \into \W_+^{k-2} \oplus  \W_+^{k} \oplus  \W_+^{k+2}$ for all $k$;
\item $d^H:\W_+^n \into \W_+^{n-2} \oplus  \W_+^{n}$.
\end{enumerate}
\end{theorem}
\begin{proof}
We will first prove that 1) implies 2). Let $\I_-$ be any complex structure on $V_-$ orthogonal \wrt\ the natural pairing so that $\G$ and  $\J_1 = \I_+ + \I_-$ form a generalized almost Hermitian structure. Then according to Proposition \ref{prop:SKT from nijenhuis}, $N_2,N_3$ and $N_+$ vanish. Then, we see that $d^H$ splits into three components:
\begin{equation}\label{eq:skt dh decomposition}
d^H = \delta_+^N + \bar{\delta_+^N} + \slashed \delta_-
\end{equation}
\vskip-18pt
\begin{IEEEeqnarray*}{l l l l}
\slashed\delta_-&:\W_+^k \into \W_+^k;&\qquad \slashed\delta_- &= \delta_-+ \deltabar_- + N_- + \bar{N_-},\\
\delta_+^N&:\W_+^k \into \W_+^{k+2};&\qquad \delta_+ ^N&= \delta_++N_1+ N_4,\\
\bar{\delta_+^N}&:\W_+^k \into \W_+^{k-2};&\qquad \bar{\delta_+^N} &= \deltabar_++\bar{N_1}+ \bar{N_4}.
\end{IEEEeqnarray*} 

Condition 2) clearly implies condition 3). Finally to prove that 3) implies 1) we once again choose a complex structure $\I_-$ on $V_-$ and observe that since $\bar{N_2}, \bar{N_3}$ and $\bar{N_+}$ are tensors, it is enough to check that they vanish when applied to spaces where their action is effective. But for $\gf \in \U^{n,0} \subset \W_+^n$, the components of $d^H$ landing on $\W_+^{n-6}$ and $\W_+^{n-4}$ are precisely $\bar{N_+}\cdot \gf$ and $\bar{N_3}\cdot\gf$ and the vanishing of these forms for $\gf \neq 0$ implies that $\bar{N_+} = \bar{N_3} =0$. Similarly, if $\psi \in \U^{0,n}$, the $\W_+^{n-4}$ component of $d^H$ is $\bar{N_2}\cdot\psi$, hence the vanishing of this component implies that $\bar{N_2}=0$ and Proposition \ref{prop:SKT from nijenhuis} implies that $\I_+$ is integrable.
\end{proof}

While the different  $U^{p,q}$ components of $d^H$ obtained in terms of the generalized almost Hermitian extension of the SKT structure depend on the particular extension chosen, the operators $\delta_+^N$, $\bar{\delta_+^N}$ and $\slashed\delta_-$ depend only on the SKT data, since they are the decomposition of $d^H$ obtained from the eigenspaces of  $\I_+$. Since $(d^H)^2 =0$ these operators satisfy some relations.

\begin{corollary}\label{cor:delta_+^N commutations}
The following hold
$$(\delta_+^N)^2 = (\overline{\delta_+^N})^2 = 0;\qquad \{\delta_+^N,\slashed \delta_-\} =0; \qquad  \{\bar{\delta_+^N},\slashed \delta_-\} =0;\qquad \{\delta_+^N,\overline{\delta_+^N}\}+  (\slashed\delta_-^H)^2 =0.$$
\end{corollary}


In the arguments up to now, given a positive SKT structure, we have chosen  $\I_-$ rather freely, but the complex structure $\I$ corresponding to the SKT data (see Proposition \ref{prop:SKT to gc}) is integrable and together with $\G$ forms a generalized Hermitian extension of $(\G,\I_+)$, that is, an SKT always has a generalized Hermitian extension.

\begin{theorem}\label{theo:genhermitian}
Let $(M,H)$ be a manifold with 3-form and $(\G,\J_1)$ be an generalized almost Hermitian structure on $M$. Then the following are equivalent:
\begin{enumerate}
\item $(\G,\J_1)$ is an integrable \gc\ extension of a positive SKT structure;
\item $d^H: \U^{p,q} \into \U^{p+1,q+1} \oplus \U^{p+1,q-1} \oplus\U^{p+1,q-3} \oplus\U^{p-1,q+3} \oplus\U^{p-1,q+1} \oplus\U^{p-1,q-1}$ for all $p,q \in \Z$;
\item $d^H: \U^{p,n-p} \into \U^{p+1,n-p-1} \oplus \U^{p+1,n-p-3} \oplus \U^{p-1,n-p+1} \oplus \U^{p-1,n-p-1}$ for all $p\in \Z$;
\item $d^H: \U^{0,n} \into \U^{-1,n-1} \oplus \U^{1,n-1} \oplus \U^{1,n-3}$ and $d^H: \U^{n,0} \into \U^{n-1,1} \oplus \U^{n-1,-1}$. 
\end{enumerate}
\end{theorem}

\begin{proof}
The fact that 1) implies 2) follows from Proposition \ref{prop:SKT from nijenhuis} and integrability of $\J_1$. The implications $2)\Rightarrow 3) \Rightarrow 4)$ are immediate. Finally,  similarly to the proof of Theorem \ref{theo:skt}, the values of the different components of the intrinsic torsion are fully determined by their action on $\U^{n,0}$ and $\U^{0,n}$ and 4) implies the vanishing of all components of the intrinsic torsion except from $N_1$.
\end{proof}

This theorem allows us to define six differential operators on a \gc\ extension of a positive SKT  structure:
\begin{IEEEeqnarray*}{c}
d^H = \delta_- + \deltabar_-+ N + \overline{N} + \delta_+ + \deltabar_+\IEEEyesnumber\label{eq:gen hermitian decomposition}\\ 
\delta_-:\U^{p,q}\into \U^{p+1,q-1};\qquad \delta_+:\U^{p,q}\into \U^{p+1,q+1};\qquad N:\U^{p,q}\into \U^{p-1,q+3}.
\end{IEEEeqnarray*}
\begin{figure}[h!!]
\begin{center}
$$\xymatrix@R=12pt@C=12pt{
\U^{p-1,q+3}&               && &\\
\U^{p-1,q+1} &               &\U^{p+1,q+1} &              &\\
               &\ar[rdd]_{\overline{N}}\ar[ld]_{\deltabar_+}\ar[lu]^{\deltabar_-} \U^{p,q} \ar[rd]^{\delta_-} \ar[ru]_{\delta_+}\ar[luu]_{N}         &   & &\\
\U^{p-1,q-1} &               &\U^{p+1,q-1} &            &  \\
&               & \U^{p+1,q-3} &              &\\
}
$$
\caption{Decomposition of $d^H$ for a \gc\ extension of a positive  SKT structure.}\label{fig:skt}
\end{center}
\end{figure}
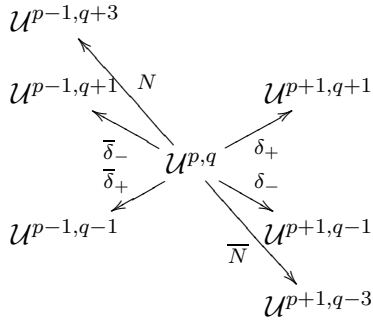

We see that $N = N_1$ is the Nijenhuis tensor of $\J_2$ and this theorem shows that if $(\G,\J_1)$ is a generalized Hermitian extension of an SKT structure then half of the Nijenhuis tensor of $\J_2$ vanishes, namely the component $N_2$ from \eqref{eq:operators}. Also, this decomposition allows us to express the operators $\delta^N_+$ and $\slashed\delta_-$ from \eqref{eq:skt dh decomposition} in terms of $\delta_+$, $\delta_-$, $N$ and their conjugates:
\begin{equation}\label{eq:sktoperators}
\delta_+^N = \delta_+ + N;\qquad \slashed\delta_- = \delta_- + \deltabar_-.
\end{equation} 

Since $(d^H)^2 =0$ these operators satisfy a number of relations.

\begin{corollary}
The following relations and their complex conjugates hold:
$$  \delta_+^2 =N^2 = 0;\qquad  \{\delta_+, N\}=0;$$
$$ \delta_-^2 = - \{\overline{N},\delta_+\};\qquad \{\overline{N},\delta_-\} =0;  \qquad \{\delta_-,\delta_+\}=0;\qquad \{\delta_+,\deltabar_-\}  =- \{\delta_-,N\};$$
$$ \{\delta_+, \deltabar_+\} + \{\delta_-,\deltabar_-\} + \{N,\overline{N}\}=0.$$
\end{corollary}

\begin{corollary}
Let $(\G,\J_1)$ be a generalized Hermitian structure on a manifold  with 3-form $(M,H)$ and let $\J_2 = \G \J_1$ be the associated \gacs. If  the canonical bundle of $\J_2$ admits a $\delbar_{\J_1}$-closed trivialization, then $(\G,\J)$ is an extension of an SKT structure.
\end{corollary}
\begin{proof}
Since $\J_1$ is integrable, $d^H:\mc{U}^{n,0}\into \U^{n-1,1} \oplus \U^{n-1,-1}$. Due to Lemma \ref{lem:nijenhuis}, we have that
$$d^H:\U^{0,n}\into \mc{U}^{-1,n-3} \oplus \U^{-1,n-1}\oplus \U^{1,n-1} \oplus \U^{1,n-3}.$$
But according to the hypothesis there is a trivialization $\rho$ of $U^{0,n}$ such that $\delbar_{\J_1} \rho =0$, so, in particular, $N_2\cdot \rho$, the component  of $\delbar_{\J_1}\rho$ in $\mc{U}^{-1,n-3}$, must vanish and since it vanishes on $\rho$, $N_2$ is the zero tensor. So we have in fact 
$$d^H:\U^{0,n}\into \U^{-1,n-1}\oplus \U^{1,n-1} \oplus \U^{1,n-3}.$$
and the last condition of Theorem \eqref{theo:genhermitian} holds.
\end{proof}

The same results hold for {\it negative SKT} structures. One should bear in mind, however that the relevant spaces  for a negative structure are given by $W_-^k =\sum_{p-q=k} U^{p,q}$, that is, each $W_-^{k}$ is an antidiagonal and integrability is equivalent to
$$d^H:\W_-^k \into \W_-^{k+2}\oplus \W_-^k \oplus \W_-^{k-2},$$
or, in terms of the $U^{p,q}$ decomposition obtained by choosing a \gc\ extension, the only nontrivial component the Nijenhuis tensor is $N_2$.

So an SKT structure (positive or negative) corresponds to a generalized Hermitian structure in which half of the Nijenhuis tensor of $\J_2$ vanishes.

\section{Hodge theory}\label{sec:Hodge}

In this section we develop Hodge theory for manifolds with parallel Hermitian, bi-Hermitian and SKT structures. Our main result is that for  a parallel positive Hermitian structure the Laplacian preserves the spaces $\W_+^{k,l}$ and hence  induces a decomposition of the $d^H$-cohomology accordingly. This is in contrast with the fact that for usual manifolds the $d^H$-cohomology has only a $\Z_2$-grading. For positive SKT manifolds, not only does the Laplacian preserve the spaces  $\W_+^{k,l}$, but in fact there is an identity between the Laplacians for the operators $d^H$, $\delta_+^N$ and $\overline{\delta_+^N}$. We start with real Hodge theory and some operators of interest.

\subsection{Differential operators, their adjoints and Laplacians}\label{subsec:operators and adjoints}
Given a generalized metric and orientation on a compact manifold $M^m$, we can form the Hodge star operator which gives us a positive definite inner product on forms:
\begin{equation}\label{eq:metric on forms}
\GG(\gf,\psi) = \int_M (\gf,\star\psi)_{Ch} \qquad \mbox{for all } \gf,\psi \in \Omega^{\bullet}(M;\R).
\end{equation}

Two basic results about $d^H$ are:
\begin{lemma} If $M^m$ is compact
\begin{equation}\tag{Integration by parts}
\int_M (d^H \gf,\psi)_{Ch} = (-1)^m\int_M(\gf,d^H\psi)_{Ch}.
\end{equation}
\begin{equation}\tag{Formal adjoint}
\GG(d^H\gf,\psi) = \GG(\gf,(-1)^m\star^{-1}d^H\star\psi)
\end{equation}
hence the formal adjoint of $d^H$ is given by $d^{H*} = (-1)^{\frac{m(m+1)}{2}}\star d^H \star$.
\end{lemma}
We will be mostly interested in the Dirac operators:
$$\begin{aligned}
\D^H_+ &= \tfrac{1}{2}(d^H - (-1)^{\frac{m(m-1)}{2}}\star d^{H}\star) = \tfrac{1}{2}(d^H + (-1)^{m+1}d^{H*})\\
\D^H_- &= \tfrac{1}{2}(d^H +(-1)^{\frac{m(m-1)}{2}} \star d^{H}\star )= \tfrac{1}{2}(d^H + (-1)^{m}d^{H*}).
\end{aligned}
$$

The operators $\D^H_\pm$ relate to the projections of $d^H$ onto \gsd- and \gasd- forms
$$d^H_{\pm}: \Omega^{\bullet}(M)\into\Omega^{\bullet}_{\pm}(M),$$
$$d^H_{\pm}\gf  = \frac{1}{2}(1 \mp i^{-\frac{m(m-1)}{2}}\star) d^H \gf.$$
\begin{lemma}\label{lem:alternative dH}
Given $\gf \in \Omega^{\bullet}(M)$ let  $\gf_\pm$ be its \gsd\ and \gasd\ components, then
$$\D^H_{\pm}\gf = d^H_{\mp}\gf_+ +d^H_{\pm}\gf_-.$$
In particular we see that $\D_-^H$ preserves the spaces $\Omega_{\pm}^{\bullet}(M)$ while $\D^H_+$ maps them to each other.
\end{lemma}
\begin{proof}
\begin{align*}
\D^H_{\pm} \gf& =\tfrac{1}{2}( (d^H \mp (-1)^{\frac{m(m-1)}{2}}\star d^H\star)(\gf_+ + \gf_-)\\
&= \tfrac{1}{2}((d^H \pm (-1)^{\frac{m(m-1)}{2}}i^{\frac{m(m-1)}{2}}\star d^H)\gf_+ +  \tfrac{1}{2}((d^H \mp (-1)^{\frac{m(m-1)}{2}}i^{\frac{m(m-1)}{2}}\star d^H)\gf_-\\
&= \tfrac{1}{2}((d^H \pm i^{3\frac{m(m-1)}{2}}\star d^H)\gf_+ + \tfrac{1}{2}( (d^H \mp i^{3\frac{m(m-1)}{2}}\star d^H)\gf_-\\
&=\tfrac{1}{2}( (d^H \pm i^{-\frac{m(m-1)}{2}}\star d^H)\gf_+ + \tfrac{1}{2}( (d^H \mp i^{-\frac{m(m-1)}{2}}\star d^H)\gf_-\\
&=d^H_{\mp}\gf_+ + d^H_{\pm}\gf_-.
\end{align*}
\end{proof}

\begin{lemma}\label{prop:Laplacian1}
Let $\triangle^H = d^Hd^{H*} + d^{H*}d^H$ be the $d^H$-Laplacian. Then 
\begin{enumerate}
\item $\star \triangle^H = \triangle^H\star$ and 
\item $(-1)^{m+1}(\D^H_+)^2 = (-1)^m (\D^H_-)^2 = \tfrac{1}{4}\triangle^H.$
\end{enumerate}
Therefore $\triangle^H$ preserves the decomposition of forms into $\Omega^{\bullet}_\pm(M)$ and hence the $d^H$-cohomology of $M$ splits as \gsd- and \gasd-cohomology: $H_{d^H}(M) = \mc{H}^H_{+}(M)\oplus \mc{H}^H_{-}(M)$.
\end{lemma}

Not only does $\D^H_+$ reverses the parity of forms but, according to Lemma \ref{lem:alternative dH}, it swaps $\Omega^{\bullet}_+(M)$ with $\Omega^{\bullet}_-(M)$ and hence it is a generalization of the signature operator. Since $\D_+^H$ has the same symbol as  $d +(-1)^{m+1} d^*$ the indices of $\D_+^H:\Omega^\bullet_+(M) \into \Omega^\bullet_-(M)$ and $\D_+^H:\Omega^{ev}(M) \into \Omega^{od}(M)$ are just the usual signature and Euler characteristic, so we have

\begin{lemma}\label{lem:chi and sigma}
The Euler characteristic, $\chi$, and the  signature, $\sigma$, of a compact manifold with closed 3-form $(M,H)$  are given by
$$\chi =  \dim(\mc{H}^{H,ev}(M)) - \dim(\mc{H}^{H,od}(M));\qquad\sigma =  \dim(\mc{H}^H_+(M)) - \dim(\mc{H}^H_-(M)).$$
\end{lemma}

\subsection{Signature and Euler characteristic of almost Hermitian manifolds}\label{subsec:signature}

In a compact generalized almost Hermitian manifold we have sixteen operators induced by the $U^{p,q}$ decomposition of forms. It turns out that for these operators taking complex conjugates or adjoints are nearly the same thing. Firstly, we extend  the real inner product on forms \eqref{eq:metric on forms} to complex valued forms by requiring it  to be Hermitian. If we let $\overline{\star}$ denote the operator given by $\overline{\star} \gf = \star \overline{\gf}$, we have
$$\GG(\gf,\psi) = (\gf,\overline{\star}\psi) \qquad \mbox{ for all } \gf, \psi \in \Omega^{\bullet}(M;\C). $$

\begin{proposition}\label{prop:adjoint} Let $M^{2n}$ be a compact  generalized almost Hermitian manifold and  let $\delta$ denote any of the operators $\delta_+,\deltabar_+,\delta_-, \deltabar_-, N_\alpha$ or $\overline{N_\alpha}$, where $\alpha = 1,2,3,4, \pm$. For $\gf,\psi\in \Omega^{\bullet}(M;\C)$ we have
\begin{align}
\int_M(\delta\gf,\psi)_{Ch} &= \int_M (\gf,\delta \psi)_{Ch}\tag{Integration by parts}\\
\delta^* &= \pm \deltabar,\tag{Formal adjoint}
\end{align}
where the sign for the adjoint is positive if $\delta:\Omega_{\pm}(M;\C)\into \Omega_{\pm}(M;\C)$ and negative otherwise.
\end{proposition}
\begin{proof}
The proof that integration by parts holds is the same for all of these operators, so we consider only $\delta_+$. It is enough to consider the case when $\gf \in \U^{p,q}$ and hence $\delta_+ \gf \in \U^{p+1,q+1}$ and it pairs trivially with $\U^{k,l}$, unless $k= -p-1$ and $l = -q-1$. Hence we may assume that $\psi \in \U^{-p-1,-q-1}$ and compute 
$$\int_M(\delta_+ \gf,\psi)_{Ch} = \int_M(d^H \gf,\psi)_{Ch}= \int_M(\gf,d^H\psi)_{Ch}=\int_M( \gf,\delta_+\psi)_{Ch},$$
where in the first equality we have used that the remaining components of $d^H\gf$  do not lie in $ \U^{p+1,q+1}$, hence they pair trivially with $\psi$, in the second equality we  integrated by parts and then  reversed the argument.

For the formal adjoint, again taking $\delta= \delta_+$, $\gf \in \U^{p,q}$ and $\psi \in \U^{p+1,q+1}$ we compute:
\begin{align*}
\GG(\delta_+\gf,\psi) &= \int_M(\delta_+ \gf,\overline{\star}\psi)_{Ch}  = i^{-p-q-2} \int_M(\gf,\delta_+\bar{ \psi})_{Ch}  \\
&= i^{-p-q-2} \int_M(\gf,(-1)^{n}\overline{\star}\,\overline{\star}\delta_+\overline{\psi})_{Ch}  =i^{-p-q-2}(-1)^{n} \int_M(\gf,\overline{\star}\,\star\overline{\delta}_+\psi)_{Ch} \\
& = i^{-2p-2q-2}(-1)^{n}\int_M(\gf,\overline{\star}\overline{\delta}_+\psi)_{Ch}= - \GG(\gf,\deltabar_+ \psi),
\end{align*}
where we have used several times that on $\U^{p,q}$ $\star$ is multiplication by $i^{p+q}$ and in the last equality we used that $p+q = n$ mod $2$.
\end{proof}

Therefore we can form the  Dirac operator corresponding to, say, $\delta_+$:
$$\slashed\delta_+ = \delta_+ -\delta_+^* = \delta_+ +\deltabar_+;$$

Since the symbol sequence of $\delta_+$ is exact, we get elliptic operators for the sequences \eqref{eq:wk seq} and \eqref{eq:up seq}:
$$\begin{aligned}
\slashed\delta_+:\Omega^{ev}(M)& \into \Omega^{od}(M);\\
\slashed\delta_+:\Omega_{+}(M)&  \into \Omega_{-}(M);\\
\end{aligned}
$$

The indices of these operators are just the Euler characteristic  and the signature of $M$:

\begin{theorem}\label{theo:indices}
Let $(M,H)$ be a compact manifold with a 3-form $H$ and let $(\G,\J_1)$ be an almost  generalized Hermitian structure on $M$. Let
$$h^{ev/od} = \dim(\ker(\slashed \delta_+:\Omega^{ev/od} \into \Omega^{od/ev})),$$
$$w^{\pm} = \dim(\ker(\slashed \delta_+:\Omega_{\pm}(M) \into \Omega_{\mp}(M)));$$
Then
$$\chi(M) = h^{ev} - h^{od}  \qquad \mbox{and}\qquad \sigma(M) =w^{+}-w^{-}.$$
\end{theorem}
\begin{proof}
To prove the claim about the Euler characteristic  we only have to show that  the symbols of $\slashed \delta_+$  and  $\D_+^H: \Omega^{ev}\into \Omega^{od}$  since due to Lemma \ref{lem:chi and sigma} the index of $\D_+^H$ is the Euler characteristic.
To compare the symbols we have
$$\D_+^H = d^H - (d^H)^* \sim \delta_+ + \deltabar_+ + \delta_- +\deltabar_- -(-\deltabar_+ - \delta_+ + \deltabar_- +\delta_-) = 2 \delta_+ + 2\deltabar_+ = 2\slashed\delta_+,
$$
where $\sim$ means that the operators have the same symbol and in the second passage we used Proposition \ref{prop:adjoint}.

The proof of the claim regarding the signature is done in the same lines.
\end{proof}

\subsection{Hodge theory on parallel Hermitian manifolds}\label{subsec:hodge theory I}

\begin{theorem}\label{theo:hodge decomposition}
Let $(M,H)$ be a compact manifold with $3$-form, let  $(\G,\I_+,\I_-)$ be a generalized metric with a pair of almost Hermitian structures and $\J_1 = \I_+ + \I_-$. Then
\begin{enumerate}
\item If $\I_{+}$ (resp. $\I_-$) is a parallel positive  (resp. negative) Hermitian structure, the $d^H$-Laplacian preserves the spaces $\W^{k,l}_+$ (resp. $\W_-^{k,l}$) and hence the $d^H$-cohomology of $M$ inherits a corresponding $\Z\times \Z_2$-grading;
\item If $(\G,\J_1)$ is a parallel bi-Hermitian structure, the $d^H$-Laplacian preserves the spaces $\U^{p,q}$ and hence the $d^H$-cohomology of $M$ inherits a corresponding $\Z^2$-grading; 
\end{enumerate}
\end{theorem}
\begin{proof}
The claim for positive and negative structures are analogous and together they imply the last claim, so it is enough to prove the first claim.

If $(M,H)$ has a generalized Hermitian structure $(\G,\J_1)$, then, since each diagonal of the $\U^{p,q}$ decomposition lies in either $\Omega_+^{\bullet}(M;\C)$ or in $\Omega^{\bullet}_-(M;\C)$, so Lemma \ref{lem:alternative dH} implies that  $\D_+^H$ and $\D_-^H$ are decomposed into a sum of operators:
\begin{equation}\label{eq:D+D-}
\begin{aligned}
\D^H_+ &= \delta_+ + N_1 + N_4 + N_+ + \deltabar_+ + \bar{N_1} + \bar{N_4} +\bar{N_+};\\
\D^H_- &= \delta_- + N_2 + N_3 + N_- + \deltabar_- + \bar{N_2} + \bar{N_3}+ \bar{N_-}.
\end{aligned}
\end{equation}

If $\I_+$ is parallel, then, according to Theorem \ref{prop:parallel via nijenhuis}, $N_2=N_3=0$ and from \eqref{eq:D+D-} we have that $\D^H_- = \delta_- + \deltabar_- + N_- + \overline{N_-}$, which clearly preserves the $\W^k_+$ decomposition, therefore $\triangle^H = 4(\D^H_-)^2$ also preserves the $\W^k$ decomposition of forms. Since the Laplacian has even parity, it also preserves $\Omega^{ev}(M)$ and $\Omega^{od}(M)$, hence the Laplacian preserves the spaces $\W^{k,l}_+$.
\end{proof}

\begin{corollary} Let $M$ be a manifold with 3-form.
\begin{itemize}
\item An SKT structure on $M$ induces a $\Z\times \Z_2$-grading on the  $d^H$-cohomology;
\item (Gualtieri \cite{gualtieri-2004}) A generalized K\"ahler structure induces a $\Z^2$-grading on the $d^H$-cohomology.
\end{itemize}
\end{corollary}

For each of the cases covered in the previous theorem, we denote by $H^{\bullet,\bullet}_{d^H}(M)$ the corresponding decomposition of the $d^H$-cohomology.

\begin{corollary}(Riemann bilinear relations)
In a  compact parallel positive Hermitian manifold
$$i^{-k}(a,\overline{a})_{Ch} >0 \qquad \mbox{ for all } a \in H^{k,l}_{d^H}(M)\backslash\{0\}.$$
 \end{corollary}
 \begin{proof}
 Indeed, let $\alpha$ be the harmonic representative for the class $a\in H^{k,l}_{d^H}(M)\backslash\{0\}$. Since $\alpha \in \W_+^{k,l}$, $\star \alpha = i^k\alpha$ and hence $\overline{\star}\alpha = i^{-k}\overline{\alpha}$ and 
 $$i^{-k}(\alpha,\overline{\alpha})_{Ch}  = (\alpha,\bar{\star}\alpha)_{Ch} > 0.$$
 \end{proof}
Finally in this context we get a clearer version of Theorem \ref{theo:indices}:
\begin{corollary}
Letting $w^{k,l}$ be the dimension of the space of harmonic sections of $\W^{k,l}_+$  in a parallel Hermitian manifold we have
$$\chi(M) = \sum (-1)^l{w^{k,l}}  \qquad \mbox{and}\qquad \sigma(M) =\sum (-1)^{\frac{n-k}{2}}w^{k,l}.$$
\end{corollary}

\subsection{Hodge theory on SKT manifolds}\label{subsec:hodge theory II}

While for parallel Hermitian manifolds we could prove that the Laplacian preserves the spaces $\W^{k,l}_+$, in and SKT manifold we can go further and show that there is an identity of Laplacians:

\begin{theorem}\label{theo:Hodge} In a positive SKT manifold,
$$\triangle_{\delta_+^N}= \triangle_{\overline{\delta_+^N}}  =\tfrac{1}{4} \triangle^H.$$
\end{theorem}
\begin{proof}
From $\D^H_+ = \delta_+^N+ \overline{\delta_+^N}$, Lemma \ref{prop:Laplacian1} and  Proposition \ref{prop:adjoint} we have
$$\tfrac{1}{4}\triangle^H = -(\D^H_+)^2 = -(\delta_+^N+\overline{\delta_+^N})^2 = -(\delta_+^N-\delta_+^{N*})^2 = \triangle_{\delta_+^N}$$
\end{proof}

\subsection{Relation to Dolbeault cohomology}\label{subsec:dolbeault}

As we saw in  Proposition \ref{prop:SKT to gc}, given a  positive SKT structure $(\G,\I)$ on $M$ one can extend  $\I_+$ to a \gcs\ $\J_1$ of complex type. In this case, in the metric splitting of $\TM$,  the structures $\J_1$ and $\J_2 = \G \J_1$ are given by 
$$\J_1 = \begin{pmatrix}
I & 0\\
0&-I^*
\end{pmatrix},
\qquad
\J_2 = \begin{pmatrix}
0 & -\omega^{-1}\\
\omega &0
\end{pmatrix},$$
where $\omega = g \circ I$ and for such generalized Hermitian structure Proposition \ref{prop:kahlerupq} gives an automorphism of the space of forms which relates the $U^{p,q}$ decomposition of forms with the usual $\wedge^{p,q}T^*M$ decomposition of forms determined by the complex structure $I$.

In this section we relate these two decomposition of forms and corresponding cohomologies. Precisely, the effect of applying the automorphism  $\Psi$  from Proposition \ref{prop:kahlerupq} to $\wedge^{\bullet}T^*_\C M$ on $d^H$ is simply to conjugate it by $\Psi$, so we get a new operator
$$\hat{d^H} = \Psi^{-1}\circ d^H  \circ \Psi.$$
Since  $d^H$ splits according to the $U^{p,q}$ into six operators, the same is true for $\hat{d^H}$ and we can define
\begin{align*}
\hat\delta_+: \Omega^{p,q}\into \Omega^{p-1,q};&\qquad \hat\deltabar_+: \Omega^{p,q}\into \Omega^{p+1,q};\\
\hat\delta_-: \Omega^{p,q}\into \Omega^{p,q+1};& \qquad \hat\deltabar_-: \Omega^{p,q}\into \Omega^{p,q+1};\\
\hat{N}:\Omega^{p,q}\into \Omega^{p-1,q-2}&\qquad \hat{\bar{N}} :\Omega^{p,q}\into \Omega^{p+1,q+2};
\end{align*}
$$\hat{d^H} = \hat\delta_+ + \hat\deltabar_+ +\hat\delta_-+ \hat\deltabar_- + \hat{N} + \hat{\bar{N}}.$$

\begin{proposition}\label{prop:relation to Dolbeault}
Let  $(\G,\I_+)$ be a positive SKT structure on a manifold with 3-form $(M,H)$, let $\J$ be the generalized complex extension of $\I_+$ given by Proposition \ref{prop:SKT to gc}.  Let
$$\del^{\omega^{-1}} =  \{\del,\omega^{-1}\} \qquad\mbox{and}\qquad \delbar^{\omega^{-1}} =  \{\delbar,\omega^{-1}\},$$
let  $\zeta_i$ be the component of $e^{\frac{i}{2}\omega^{-1}}_* \delbar\omega$ lying in $\wedge^i T \tensor \wedge^{3-i}T^*$ and let $\chi $ be the $(1,2)$-component of $[\omega^{-1},\omega^{-1}]_{SN}$, where $[\cdot,\cdot]_{SN}$ is the Schouten--Nijenhuis bracket. Then the following hold
\begin{IEEEeqnarray*}{l l}
\hat\delta_+ = \tfrac{i}{2}\delbar^{\omega^{-1}} + 2i\zeta_2&\qquad \hat\deltabar_+ = \del\\
\hat\delta_- = \delbar+ 2i\zeta_1&\qquad \hat\deltabar_- = \tfrac{i}{2}\del^{\omega^{-1}}\IEEEyesnumber\label{eq:skt dh decompositionII}\\
\hat{N} = \tfrac{1}{8}\chi&\qquad \hat{\bar{N}} = 2i\delbar{\omega}.
\end{IEEEeqnarray*}
\end{proposition}
\begin{proof}
The proof is a direct computation of the operator $\hat{d^H}$:
\begin{equation}\label{eq:dhhat}
\begin{aligned}
\hat{d^H} &= e^{-\frac{i}{2}\omega^{-1}}e^{-i\omega}(d+H\wedge) e^{i\omega}e^{\frac{i}{2}\omega^{-1}}\\
&= e^{-\frac{i}{2}\omega^{-1}}e^{-i\omega} e^{i\omega}(d+ (H + i d\omega )\wedge)e^{\frac{i}{2}\omega^{-1}}\\
&= e^{-\frac{i}{2}\omega^{-1}}de^{\frac{i}{2}\omega^{-1}} +  2i e^{-\frac{i}{2}\omega^{-1}}(\delbar \omega \wedge)e^{\frac{i}{2}\omega^{-1}}\\
\end{aligned}
\end{equation}
We compute separately each of the two operators above making up $\hat{d^H}$:
\begin{align*}
e^{-\tfrac{i}{2}\omega^{-1}}de^{\tfrac{i}{2}\omega^{-1}} & = d + [d,\tfrac{i}{2}\omega^{-1}] + \tfrac{1}{2!}[[d,\tfrac{i}{2}\omega^{-1}],\tfrac{i}{2}\omega^{-1}] +\tfrac{1}{3!} [[[d,\tfrac{i}{2}\omega^{-1}],\tfrac{i}{2}\omega^{-1}],\tfrac{i}{2}\omega^{-1}] + \cdots
\end{align*}
The first term is just $d = \del + \delbar$, while the second term is a version of the symplectic adjoint of $d$ now obtained in a nonintegrable setting, $\delta^{\omega^{-1}} = d \circ \omega^{-1} - \omega^{-1} \circ d$. Since $\omega$ is not closed, $\delta^{\omega^{-1}} $ does not square to zero. Integrability of the complex structure gives  $d = \del +\delbar$ and recalling that $\omega^{-1}$ is of type $(1,1)$ we have
$$\delta^{\omega^{-1}} = \del^{\omega^{-1}} + \delbar^{\omega^{-1}}; \del^{\omega^{-1}}: \Omega^{p,q} \into \Omega^{p,q-1};\qquad \delbar^{\omega^{-1}}:\Omega^{p,q} \into \Omega^{p-1,q}.$$
Since $\omega^{-1}$ is even, the third term in the series coincides with a multiple of the expression for the derived bracket of $\omega^{-1}$ with itself, that is, the Schouten--Nijenhuis bracket:
$$[[d,\tfrac{i}{2}\omega^{-1}],\tfrac{i}{2}\omega^{-1}]  = -\{\{\tfrac{i}{2}\omega^{-1},d\},\tfrac{i}{2}\omega^{-1} \} = \tfrac{1}{4}[\omega^{-1},\omega^{-1}]_{SN}.$$
Since $\omega^{-1}$ is of type $(1,1)$, integrability of the complex structure implies that this term this is a 3-vector lies in $\wedge^{2,1}T \oplus \wedge^{1,2}T$, so this term decomposes as
$$[\omega^{-1},\omega^{-1}]_{SN}= \chi+ \overline{\chi}, \qquad \chi \in \wedge^{1,2} TM.$$
The fourth term is the commutator $[[\omega^{-1},\omega^{-1}]_{SN},\omega^{-1}]$ which vanishes since $\omega^{-1}$ is a bivector and hence so do the remaining terms in the series. So we have established that
\begin{equation}\label{eq:expression for d}
e^{-\frac{i}{2}\omega^{-1}}d e^{\frac{i}{2}\omega^{-1}} = \del + \delbar+ \tfrac{i}{2}\del^{\omega^{-1}} + \tfrac{i}{2}\delbar^{\omega^{-1}}  + \tfrac{1}{8}(\chi+ \bar{\chi})
\end{equation}

Next we compute the second summand in \eqref{eq:dhhat}:
\begin{align*}
2ie^{-\frac{i}{2} \omega^{-1}} \circ  \delbar\omega \circ e^{\frac{i}{2} \omega^{-1}}&  = 2ie^{-\frac{i}{2} \omega^{-1}} \circ (e_*^{\frac{i}{2}\omega^{-1}}e_*^{-\frac{i}{2}\omega^{-1}}  \delbar\omega) \circ e^{\frac{i}{2} \omega^{-1}}\\
& = 2ie^{-\frac{i}{2} \omega^{-1}} \circ e^{\frac{i}{2}\omega^{-1}}\circ(e_*^{-\frac{i}{2}\omega^{-1}} \delbar\omega)\\
& = 2ie_*^{-\frac{i}{2}\omega^{-1}}\delbar\omega,
\end{align*}
where we have used \eqref{eq:spin action} in the second equality. The element $e_*^{-\frac{i}{2}\omega^{-1}}\delbar\omega \in \wedge^3 \T_\C M$ has six components:
$$e_*^{-\frac{i}{2}\omega^{-1}}\delbar\omega \in (\wedge ^{1,2}T^*) \oplus  (T^{1,0}\tensor \wedge^{1,1} T^*) \oplus (T^{0,1}\tensor \wedge^{0,2}T^*)\oplus (\wedge^{2,0}T \tensor T^{*1,0})\oplus (\wedge^{1,1}T \tensor T^{*0,1})  \oplus (\wedge^{2,1}T)$$

Since, for a 1-form $\xi$,  $e_*^{-\frac{i}{2}\omega^{-1}} \xi = \frac{i}{2}\omega^{-1}(\xi) + \xi$, we see that the 3-form component is $\delbar\omega$ and that the 3-vector component is $-\tfrac{i}{8}\omega^{-1}\delbar\omega$. We let $\zeta_1$ be the component of $e_*^{-\frac{i}{2}\omega^{-1}}\delbar\omega$ in $T\tensor \wedge^2T^*$ and $\zeta_2$ be the component in $\wedge^2 T\tensor T^*$, so that $\zeta_1:\wedge^{p,q} \into \wedge^{p,q+1}$, $\zeta_2:\wedge^{p,q}\into \wedge^{p-1,q}$, and
$$2ie^{-\frac{i}{2} \omega^{-1}} \circ  \delbar\omega \circ e^{\frac{i}{2} \omega^{-1}} = 2i\delbar\omega + 2i \zeta_1 + 2i \zeta_2 + \frac{1}{4} \omega^{-1}(\delbar\omega)$$

Putting this together with \eqref{eq:expression for d} and the fact that $\hat{d^H}$ does not have a component mapping $\Omega^{p,q}(M)$ to $\Omega^{p-2,q-1}(M)$ we conclude that $\bar{\chi}=-2\omega^{-1}(\delbar \omega)$ and obtain \eqref{eq:skt dh decompositionII}.
\end{proof}

Next, we let $\del^{i\delbar\omega} $ be the operator $ \del + i \delbar \omega\wedge$. Then we observe that even though $\del^{i\delbar \omega}$ does not preserve the degree of a form, it preserves the holomorphic degree amd the parity of the anti-holomorphic degree, i.e., 
$$\del^{ i \delbar\omega}:\Omega^{p,q}(M)\into \Omega^{p+1,q}(M),\quad p\in\Z,~q\in \Z_2.$$
Therefore its cohomology has a natural $\Z\times \Z_2$-grading. For an operator $\delta$ with $\delta^2=0$ we denote its cohomology by $H_{\delta}(M)$.

\begin{corollary}\label{cor:cohomology isomorphisms}
Let $(\G,\I_+)$ be a positive SKT structure on a compact manifold, let $\J$ be the \gc\ extension from Proposition \ref{prop:SKT to gc} and let $(g,I)$ be the corresponding Hermitian structure with Hermitian form $\omega = g \circ I$. Then for $p,q \in \Z$
$$H^{p,q}_{\del}(M) \cong H^{q-p,n-p-q}_{\bar{\delta_+}}(M),$$
and
$$H^p_{\del^{i\delbar \omega}}(M) \cong H^{n-2p}_{\bar{\delta_+^N}}(M) \cong H_{d^H}^{n-2p}(M).$$
\end{corollary}
\begin{proof}
Indeed, $\Psi$ puts  $\Omega^{p,q}(M)$ in correspondence with $\U^{q-p,n-p-q}$ and $\Omega^{p,\bullet}(M)$ with $\W^{n-2p}$ and according to Proposition \ref{prop:relation to Dolbeault}, $\Psi \circ \del = \bar{\delta_+} \circ \Psi$ and $\Psi \circ \del^{2i\omega} = \bar{\delta_+^N}\circ \Psi$. Therefore $H^{p,q}_{\del}(M) \cong H^{q-p,n-p-q}_{\bar{\delta_+}}(M)$ and  $H^{p}_{\del^{2i\delbar\omega}}(M) \cong H^{n-2p}_{\bar{\delta_+^N}}(M)$. Of course the cohomology of $\del^{2i\delbar\omega}$ is the same as the cohomology if $\del^{i\delbar \omega}$, as the automorphism
$$m:\wedge^{\bullet,\bullet}T^*M \into m:\wedge^{\bullet,\bullet}T^*M \qquad m(\alpha) = 2^{\frac{q}{2}}\alpha \qquad \mbox{ for all } \alpha \in \wedge^{p,q}T^* M,$$
relates them.
\end{proof}

Interestingly we have that $d^H = \del^{i\delbar \omega} + \delbar^{-i\del\omega}$ and the operators $ \del^{i\delbar \omega}$ and $ \delbar^{-i\del\omega}$ both square to zero (due to the SKT condition $\del\delbar\omega =0$) and graded commute. Yet, they are not well behaved \wrt\ the $\Z^2$-bigrading of forms and only on a $\Z_2^2$-grading do they behave nicely, as there we have
$$\del^{i\delbar \omega}: \Omega^{p,q}(M) \into \Omega^{p+1,q}(M) \qquad \delbar^{-i\del \omega}: \Omega^{p,q}(M) \into \Omega^{p,q+1}(M),\qquad (p,q) \in \Z_2^2. $$
There is a standard trick to  recover one $\Z$-grading, say the `$p$'-grading. We introduce a new formal variable, $\beta$, alter  the operators $\del^{i\delbar \omega}$ and  $\delbar^{-i\del \omega}$ so that they act on $E_1 = \Omega^{p,q}(M)\tensor\IP{\beta}$:
$$\del^{i\delbar \omega}( \alpha \beta^k) = (\del^{i\delbar \omega} \alpha) \beta^k\qquad \delbar^{-i\del \omega}(\alpha \beta^k) = (\delbar\alpha )\beta^k -i\del \omega\wedge \alpha \beta^{k+1}$$
and introduce a $\Z\times \Z_2$-grading on $E_1$ by declaring that for $ \Omega^{p,q}(M)\tensor \beta^k = E_1^{p-2k,q}$, $(p-2k,q) \in \Z\times \Z_2$. I.e., $\beta$ has degree $(-2,0)$. With this arrangement 
$$\del^{i\delbar \omega}: E_1^{p,q}(M) \into E_1^{p+1,q}(M) \qquad \delbar^{-i\del \omega}:E_1^{p,q}(M) \into E_1^{p,q+1}(M),\qquad (p,q) \in \Z\times \Z_2.$$

In this setting, one can produce a $\Z \times \Z_2$-graded spectral sequence whose first page is $H_{\del^{i\delbar \omega}}^{\bullet}(M)\tensor \IP{\beta}$ and the last is  $H_{d^H}^{\bullet}(M)\tensor \IP{\beta}$. Since these are isomorphic, by Corollary \ref{cor:cohomology isomorphisms}, this sequence degenerates at the first page.

\begin{theorem}\label{theo:frolicher}{\em (SKT spectral sequence)}
In an SKT manifold the decomposition $d^{H} = (\del +i\delbar\omega) + (\delbar -i\del\omega)$ gives rise to a spectral sequence which degenerates at the first page.
\end{theorem}

On a more general note, the isomorphism $\Psi$ used to prove Proposition \ref{prop:relation to Dolbeault} can be used even in the nonintegrable case to give us information about the Euler characteristic and signature of almost Hermitian manifolds. Indeed, according to Proposition \ref{prop:kahlerupq} we have
$$\Psi: \oplus_q \wedge^{p,q}T^*M\cap \Omega^{ev}(M) \into W_+^{k,0};\qquad \Psi: \oplus_q \wedge^{p,q}T^*M\cap \Omega^{od}(M) \into W_+^{k,1}.$$
And hence conjugating $\deltabar_+$ by $\Psi$ we get an operator
$$\hat{\deltabar}_+ = \Psi\circ \deltabar_+ \circ \Psi^{-1}:\Omega^{p,q}\into \Omega^{p+1,q}.$$
Since the isomorphism $\Psi$ from Proposition \ref{prop:kahlerupq} identifies $T^{*1,0}M$ with $V_+^{1,0}$, one can readily check that 
$\del$ and $\hat{\deltabar}_+$ have the same symbol and hence the same index. Of course complex conjugation swaps $\del$ and $\delbar$, allowing us to translate Theorem \ref{theo:indices} to Dolbeault cohomology terms:

\begin{corollary}\label{theo:signature}
Let $(M,I,g)$ be a compact almost Hermitian manifold, $\slashed\delbar = \delbar - \delbar^*$ and, for $q,k \in \Z_2$,
$$w^{q} = \dim(\ker \slashed \delbar:\oplus_p\Omega^{p,q}(M) \into \oplus_p\Omega^{p,q+1}(M))$$
$$h^{k} = \dim(\ker \slashed \delbar:\Omega^{k}(M) \into \Omega^{k+1}(M))$$
Then 
$$\chi(M) = h^{ev} - h^{od}  \qquad \mbox{and}\qquad \sigma(M) =w^{ev}- w^{od}.$$
In particular in a compact complex manifold if $h^{p,q} = \dim(H^{p,q}_{\delbar}(M))$, then
\begin{equation}\label{eq:signature theorem}
 \chi = \sum (-1)^{p+q}h^{p,q}\qquad\mbox{and}\qquad  \sigma = \sum_{p,q}(-1)^qh^{p,q}(M).
\end{equation}
\end{corollary}

\begin{remark}
Using Fr\"olicher's spectral sequence, Fr\"olicher  proved the identity for the Euler characteristic assuming integrability of the complex structure \cite{MR0073262}. The identity regarding the signature is an extention of the Hodge Index Theorem, which, in its modern version (see, e.g.,  \cite{MR1967689}, Theorem 6.33), states that \eqref{eq:signature theorem} holds on a compact  K\"ahler manifold and is obtained as a consequence of the development of Hodge theory of K\"ahler manifolds. The main point of Corollary \ref{theo:signature} is that these identities do not depend on existence of K\"ahler structures or on the integrability of the complex structure, so, in effect, {\it they do not represent an obstruction to the existence of any of the structures studied here}. This is to be compared with Theorem \ref{theo:frolicher} which provides a nontrivial differential--topological obstruction to the existence of SKT structures on complex manifolds.
\end{remark}

Quite separate from the theory developed so far, one can get other obstructions to the existence of SKT structures using only classical tools. For example:

\begin{theorem}\label{theo:obstruction2}
A compact SKT manifold $(M,I,g)$ for which
\begin{equation}\label{eq:restriction}
H^{2,1}_\delbar(M) = H_\delbar^{3,0}(M) = \{0\}.
\end{equation}
admits a symplectic structure.
\end{theorem}
\begin{proof}
Indeed, if there was such an SKT manifold, $(M,g,I)$, letting $H = d^c\omega$ we have that $\del\omega$ represents a class in $H^{2,1}_\delbar(M)$ which is trivial, and hence there is $\sigma \in \Omega^{2,0}(M)$ such that $ \del\omega= \delbar \sigma$. Therefore, we have
$$\delbar \del\sigma = - \del \delbar \sigma = -\del^2 \omega =0,$$
showing that  $\del\sigma$ represents a class in $H^{3,0}_\delbar(M) =\{0\}$, hence $\del \sigma =0$ and
$$\tilde{\omega} = \omega -\sigma - \bar\sigma$$
is a closed form. Further, for any $X \in T_pM\backslash \{0\}$,
$$\tilde{\omega}(X,I X) = \omega(X,IX) -i\sigma(X,X) + i\bar{\sigma}(X,X) = g(X,X) \neq 0,$$
showing that $\tilde{\omega}$ is also nondegenerate, i.e., is a symplectic structure.
\end{proof}

\begin{example}[Calabi--Eckman manifolds]\label{ex:calabi-eckman} Among the Calabi-Eckman manifolds, $S^1 \times S^1$, $S^1\times S^3$ and $S^3 \times S^3$ are known to admit SKT structures, by virtue of being compact Lie groups. In this example we show that all the remaining  Calabi--Eckman manifolds $M_{u,v} \cong S^{2u+1}\times S^{2v+1}$ do not admit SKT structures and we give two arguments for this fact.

For the first argument, we observe that the claim can be proved directly using Theorem \ref{theo:obstruction2}:  the Dolbeault cohomology of these manifolds was computed by Borel \cite{MR0202713} and, assuming $v\geq u$, it is given by
\begin{equation}\label{eq:dolbeault calabi-eckman}
H_{\delbar}^{\bullet,\bullet}(M_{u,v}) = \wedge \mbox{span}\{x^{0,1},z^{v+1,v}\} \tensor \C[y^{1,1}]/y^{u+1},
\end{equation}
where $a^{p,q}$ is a generator of bidegree $(p,q)$. Hence, for all $M_{u,v}$, with exception of the three cases known to admit SKT structures, the hypothesis of Theorem \ref{theo:obstruction2} hold but $H^2(M_{u,v}) = \{0\}$, hence these manifolds can not be symplectic.

The second argument works for $u\neq 0$ and amounts to proving that in this case, the $d^H$-cohomology of $M_{u,v}$ is not  isomorphic to the $\delbar^{i\del\omega}$-cohomology, and hence Corollary \ref{cor:cohomology isomorphisms} fails. To prove this claim, we start considering the case $v\geq u>1$. If such a manifold had an SKT structure, then $[d^c\omega] =0$ and the $d^H$-cohomology would be isomorphic to the de Rham cohomology. From \eqref{eq:dolbeault calabi-eckman}, one  sees that $H^{2,1}_\delbar(M_{u,v}) = \{0\}$ and hence $\del\omega$ would be the trivial Dolbeault class and the $\delbar^{i\del\omega}$-cohomology would be isomorphic to the usual Dolbeault cohomology. Since the Dolbeault cohomology of the Calabi-Eckman manifolds is not isomorphic to the de Rham cohomology, $M_{u,v}$ can not be SKT.

The same argument holds for $v> u=1$, except that now  $H^3(M_{v,1}) = \R$ and there are two possibilities for $[H]$: it is either zero or not. In the former case, the $d^H$-cohomology is isomorphic to the de Rham cohomology and in the latter it vanishes completely. On the other hand, the Dolbeault cohomology is still described by  \eqref{eq:dolbeault calabi-eckman} and, as before, $H^{2,1}_{\delbar}(M)$ vanishes so the $\delbar^{i\del \omega}$-cohomology is isomorphic to the $\delbar$-cohomology which is not isomorphic to either the de Rham cohomology or the trivial one.

\hfill $\blacksquare$\end{example}

\subsection{Hermitian symplectic structures}

A particular type of SKT structure for which Theorem \ref{theo:frolicher} is particularly relevant are the so called {\it Hermitian symplectic structures}. These consist of a pair $(I,\omega)$ of integrable complex structure and symplectic structure such that $\omega(X,IX)$ is positive for every nonzero vector $X$. The difference between these structures and K\"ahler structures is that here we do not require $\omega$ to be of type (1,1). Yet, we can decompose $\omega$ into its $(p,q)$ components \wrt\ the complex structure and then one readily obtains that
$$d^c\omega^{1,1} = i d(\omega^{2,0}-\omega^{0,2})$$
hence $(I,\omega^{1,1})$ is an SKT structure and $H =id  (\omega^{2,0}-\omega^{0,2})$ is exact hence the $d^H$-cohomology is isomorphic to the de Rham cohomology and $H^{1,2} = -i\del\omega^{0,2}$ is $\del$-exact hence the $\del^{-i\delbar\omega^{1,1}}$-cohomology is isomorphic to the $\del$-cohomology and Theorem \ref{theo:frolicher} gives:
\begin{corollary}\label{cor:hermitian symplectic}
In a compact Hermitian symplectic manifold the Fr\"olicher spectral sequence degenerates at the second page, i.e., the Dolbeault and the de Rham cohomologies are isomorphic as graded vector spaces. 
\end{corollary}

\section{Hodge theory beyond $U(n)$}\label{sec:beyond}

The decomposition of harmonics into their $(p,q)$-components in a K\"ahler manifold is a phenomenon that repeats itself for any other special holonomy group and, as such, is a result on Riemannian geometry. This approach is quite different from what we have done so far as, just like in the original K\"ahler identities, we relied on the underlying complex structures heavily to develop our theory. This section we show that with the appropriate setup, our results can also be extended to any other holonomy groups. Throughout this section we let $\{e_1, \cdots, e_m\} \in \Gamma(TM)$ be an orthonormal frame and  $\{e^1,\cdots, e^m\} \in \Gamma(T^*M)$  be  its dual frame. 

As before, given an oriented Riemannian manifold with closed 3-form, $(M,g,H)$ we consider the metric connection $\hat\nabla^{\pm}$ whose torsion is $\mp H$ and let $\nabla$ denote the Levi-Civita connection. Using the orthonormal frame, we can write explicit expressions for $\nabla^{\pm}$. Indeed, if we define $h_{ijk} = H(e_i,e_j,e_k)$, then,  using Einstein summation convention and omitting the symbol for the wedge product, $H$ is given by
$$ H = \tfrac{1}{3!}h_{ijk}e^i e^j e^k.$$
and we have
\begin{equation}\label{eq:relation}
\nabla^{\pm} = \nabla \mp \tfrac{1}{4}h_{ijk}e^i e^j e_k.
\end{equation}

If the holonomy of $\nabla^+$ is the Lie group $G_+$, then using the isomorphism $TM \cong V_+$, we realize its Lie algebra, $\Gg_{+}$, as a sub Lie algebra of $\frak{so}(V_+) = \wedge^2 V_+$. {\it Mutatis mutandis}, the same holds for $\nabla^-$ and we get $\Gg_- \subset \wedge^2 V_-$. Next we notice that $\Gg_+ \oplus \Gg_- \subset \wedge^2 \TM = \frak{spin}(\TM)$ and hence the elements of $\Gg_+ \oplus \Gg_- $ act on forms, thought of as spinors. Further, since $V_+$ is orthogonal to $V_-$ and $\Gg_{\pm} \subset \wedge^2V_\pm$,  the Lie algebra action of $\Gg_+$ and $\Gg_-$ on forms commute and we get an action of  $\Gg_+\oplus \Gg_-$ on forms as a direct sum of the individual actions of the Lie algebras. Now we are in condition to state the main theorem of this section, which extends Theorem \ref{theo:hodge decomposition} to general holonomy groups.

\begin{theorem}\label{theo:reduced holonomy}
Let $(M,g,H)$ be a compact Riemannian manifold endowed with a closed 3-form. If $\nabla^{\pm}$, the metric connections  with torsion $\mp H$, have holonomy in $G_\pm$, then the $d^H$-cohomology of $M$ splits according to the decomposition of forms into irreducible representations of the action of $\Gg_+\oplus \Gg_- \subset \wedge^2 V_+ \oplus \wedge^2V_-$ on forms.
\end{theorem}

\begin{lemma}\label{lem:preserves}
The connections $\nabla^\pm$ preserves the irreducible representations of $\Gg_\pm$, respectivelly.
\end{lemma}
\begin{proof}
Indeed, if we denote by $\tilde{\Gg}_+ \subset \frak{spin}(TM) \cong \wedge^2T^*M$ the bundle of endomorphisms of $TM$ defined by the connection $\nabla^+$, the condition that $\nabla^+$ has reduced holonomy implies that this bundle is preserved by parallel transport. Now $\Gg_+$ is simply the image of $\tilde{\Gg_+}$ by the parallel isomorphism $\Id + g: T^*M \into V_+$, hence the bundle $\Gg_+ \subset \wedge^2 V_+$ is also parallel and its irreducible representations are preserved by the connection. 
\end{proof}

Next we define
$$\e^i_{\pm} = e^i \pm e_i \in \Gamma(V_\pm)$$
and
$$H_{\pm} = \tfrac{1}{3!}h_{ijk}\e^i_\pm \e^j_\pm  \e^k_\pm \in \Clif^3(V_\pm) \subset \Clif^3(\TM).$$

To prove the theorem we need to  extend  to the torsion case the formulas relating the Levi--Civita connection with the exterior derivative and its adjoint:
\begin{equation}\label{eq:classical}
d = e^i \wedge \nabla_{e_i}; \qquad d^*= -\iota_{e_i}\nabla_{e_i}.
\end{equation}

To avoid considering separate cases according the whether the dimension of the manifold is even or odd, we let
$$\slashed{\mc{D}}_+ = d^H - d^{H*}\qquad \mbox{and}\qquad \slashed{\mc{D}}_- = d^H +d^{H*}.$$

\begin{lemma}\label{lem:dirac}
With the definitions above
\begin{equation}\label{eq:first order identity}
\slashed{\mc{D}}^H_+ = \e^i_+ \nabla_{e_i}^{+} + H_+;\qquad 2\slashed{\mc{D}}^H_- = \e^i_- \nabla_{e_i}^{-} + H_-.
\end{equation}
\end{lemma}
\begin{proof} 
Since the statements are similar, we will only prove the first identity. The left hand side of \eqref{eq:first order identity} is
\begin{equation}\label{eq:lefthand}
\begin{aligned}
\slashed{\mc{D}}^{H}_{+}  & = d + H - d^*- H^*\\
&= d - d^* + \tfrac{1}{3!}h_{ijk}e^i e^j e^k + \tfrac{1}{3!} h_{ijk}e_i e_j  e_k.
\end{aligned}
\end{equation}

Next, we use relation \eqref{eq:relation} to re-write the right hand side as
$$\e^i_+ \nabla_{e_i}^{+} + H_+ = e^i \wedge  \nabla_{e_i} + \iota_{e_i} \nabla_{e_i} - \tfrac{1}{2}h_{ijk}e^ie^j e_k - \tfrac{1}{2}h_{ijk}e_ie^j e_k   +H_+.$$
Due to  \eqref{eq:classical},  the  first two terms are $d- d^*$. Expanding $H_+$  we get 
$$\e^i_+\nabla_{e_i}^{+} + H_+ =  d - d^* - \tfrac{1}{2} h_{ijk}e^ie^j e_k  -\tfrac{1}{2} h_{ijk}e_ie^j e_k   +  \tfrac{1}{3!}h_{ijk}e^i e^j e^k + \tfrac{1}{2!}h_{ijk}e^i e^j  e_k +\tfrac{1}{2!}h_{ijk}e^ie_j  e_k +\tfrac{1}{3!}h_{ijk}e_i e_j  e_k.$$
which equals \eqref{eq:lefthand}.
\end{proof}

\begin{proof}[Proof of Theorem \ref{theo:reduced holonomy}]
Since $\nabla^+$ preserves the irreducible representations of $\Gg_+$, Lemma \ref{lem:dirac} implies that $\slashed{\mc{D}}^H_+$ preserves the irreducible representations of $\Gg_+ \oplus \Clif(V_-)$ and hence so does the $d^H$-Laplacian. Similarly, $\slashed{\mc{D}}^H_-$ preserves the irreducible representations of $\Clif(V_+) \oplus \Gg_-$ and hence $\Delta_{d^H}$ preserves the intersections of these representations, which are just the irreducible representations of $\Gg_+ \oplus \Gg_-$.
\end{proof}

\subsection{Integrability}

A simple consequence of \eqref{eq:classical} is that if $W \subset \wedge^k T^*M$ is a representation of the holonomy group of the Levi--Civita connection, then the exterior derivative restricted to sections of $W$ can only land in representations present in $T^*M \wedge W$. This is the Riemannian version of the claim that in a complex manifold $d:\Omega^{p,q}(M) \into \Omega^{p+1,q}(M) \oplus \Omega^{p,q+1}(M)$. Lemma \ref{lem:dirac} has similar implications.

\begin{proposition}\label{prop:dH and irreducibles}
Let $(M,g,H)$ be a compact oriented Riemannian manifold with a closed 3-form and let $\Gg_\pm \subset \wedge^2  V_\pm$ be the Lie algebras of the holonomy groups of the connections $\nabla^\pm$.  Let $W$ be a representation of $\Gg_+ \oplus \Gg_-$, then $d^H$ sends sections of $W$ into sections of representations that appear in $(\Clif^3(V_+)\oplus \Clif^3(V_-)) \cdot W$. 
\end{proposition}
\begin{proof}
It follows from Lemmata \ref{lem:preserves} and \ref{lem:dirac} that
$$\slashed{\mc{D}}^H_\pm: \Gamma(W) \into \Gamma(\Clif^3(V_\mp)\cdot W),$$
hence, $d^H = \frac{1}{2}(\slashed{\mc{D}}^H_+ + \slashed{\mc{D}}^H_-)$ has the stated property.
\end{proof}

It follows from the material in Sections \ref{sec:parallel} and \ref{subsec:integrability} that the difference between an SKT or  generalized K\"ahler structure and parallel (almost) Hermitian or (almost) bi-Hermitian structure, respectively, is that in the former two cases  if $W$ is a representation of $\Gg_+\oplus \Gg_-$ then $d^H:\Gamma(W) \into \Gamma(\T_\C M \cdot W)$, while in the latter two cases Proposition \ref{prop:dH and irreducibles} is the best one can say. This suggests that in general there is a subclass of the space of manifolds $(M,g,H)$  with reduced holonomy which may be of further interest:

\begin{definition}
We say that $(\nabla^+,\nabla^-)$ induces an {\it integrable} $G_+ \times G_-$ structure if the holonomy of $\nabla^{\pm}$ is $G_{\pm}$ and
\begin{equation}\label{eq:integrability}
d^H:\Gamma(W) \into \Gamma(\T_\C M \cdot W)
\end{equation}
for every representation $W \subset \wedge^{\bullet}T^*_\C M$  of $\Gg_+\oplus \Gg_-$.
\end{definition}

\section{Instantons over complex surfaces}\label{subsec:instantons}

Next we use the techniques developed in this section to provide an alternative description of the SKT structure on the moduli space of instantons of a bundle over a complex surface. The tool used to describe this structure are extended actions as introduced in \cite{MR2323543} and the SKT reduction theorem as presented in \cite{MR2314216}. The argument presented here follows closely the one from \cite{BCGinstantons}, so we will spare details and refer to that paper for further reading.

Let $(M,[g],I)$ be  a compact complex surface with a conformal Hermitian structure. By a result of Gauduchon \cite{MR0470920}, there is a representative $g$ of the conformal class which makes $(M,g,I)$ into an SKT manifold, i.e., the corresponding Hermitian form $\omega$ satisfies $dd^c\omega =0$. We let $H = d^c\omega$ and consider $\TM$ endowed with the $H$-Courant bracket, so that the metric
$$ \G = \begin{pmatrix}
0 & g^{-1}\\
g & 0
\end{pmatrix}$$
and the complex structure on $V_+$ induced by $I$ via the isomorphism $\pi:V_+ \into TM$ are a positive SKT structure on $\TM$.

Given a bundle $E$ over $M$ with a compact Lie group $G$ as structure group and Lie algebra $\Gg$ we let $\AA$ be the space of all $\Gg$-connections on $E$ endowed with the trivial 3-form, so that $\AA$ is an affine space isomorphic  to space of 1-forms on $M$ with values in the adjoint bundle $\Gg_E$, $ \Omega^{1}(M;\Gg_E)$. Hence at any connection $A$ we have $T_A\AA  = \Omega^1(M;\Gg_E)$ and, letting $\kappa$ be a bi-invariant metric on $G$, we can use $\kappa$ to identify $T^*_A\AA = \Omega^3(M;\Gg_E)$. Indeed, for $X \in \Omega^1(M;\Gg_E)$ and $\xi \in \Omega^3(M;\Gg_E)$ we define the natural pairing as 
$$\xi(X) = 2\int_M \kappa(X,\xi).$$
Then for $X,Y \in \Omega^1(M;\Gg_E)$ and $\xi,\eta \in \Omega^3(M;\Gg_E)$ we have
$$\IP{X+\xi,Y+\eta} = \int_M\kappa(X,\eta) + \kappa(Y,\xi)  = \int_M\kappa(X,\eta) - \kappa(\xi,Y)  =  \int_M\kappa(X+\xi, Y + \eta)_{Ch}$$
where $\kappa(\cdot,\cdot)_{Ch}$ indicates that one uses $\kappa$ to pair elements in $\Gg_E$ and the Chevalley pairing on forms to obtain a top degree form.

We denote by $\Gau$ the gauge group of $E$ and by $\tilde{\Gg}=  \Omega^{0}(M;\Gg_E)$ its Lie algebra. The infinitesimal generator corresponding to $\gamma \in \Omega^{0}(M;\Gg_E)$ at a point $A \in \AA$ is just the vector $d_A\gamma \in \Omega^1(M;\Gg_E)$ and we can extend this action to form a lifted action as in  \cite{MR2314216}:
$$\tilde{\Psi}:\Omega^0(M;\Gg_E) \into \T \AA\qquad \tilde\Psi(\gamma)|_A = d_A^H\gamma = d_A\gamma + H \wedge \gamma,$$
as long as there are no infinitesimal symmetries, i.e., as long as $d_A:\Omega^0(M;\Gg_E) \into \Omega^1(M;\Gg_E)$ has trivial kernel. Therefore, {\it from this point onwards we only consider connections for which $d_A:\Omega^0(M;\Gg_E) \into \Omega^1(M;\Gg_E)$ has trivial kernel}. If $E$ is a simple  $SU(n)$-bundle then that is the case for all ASD connections.

Next we add a moment map to this action. Our (equivariant) moment map takes values on the $\Gau$-module $\frak{h}^* = \Omega^{2}_+(M;\Gg_E)$, the space of self-dual  2-forms:
$$\mu:\AA \into \Omega^2_+(M;\Gg_E)\qquad \mu(A) = (F_A)_+,$$
where $F_A$ is the curvature of the connection $A$ and $(\cdot)_+$ indicates projection onto the space of self dual forms, $\Omega_+^{\bullet}(M;\Gg_E)$, which, for 2-forms in 4 dimensions, agrees with the usual self dual 2-forms. If we  let $\frak{a}$ be the sum $\tilde{\Gg}\oplus \frak{h}$, $\frak{a}$ becomes a Courant algebra if we endow it with the hemisemidirect product: 
$$\Cour{(\gamma_1,\lambda_1), (\gamma_2,\lambda_2)} = ([\gamma_1,\gamma_2],\gamma_1 \cdot \lambda_2).$$ 
Then  the maps $\tilde\Psi$ and $\mu$ together give rise to an extended action given infinitesimally by map of Courant algebras:
$$\Psi:\frak{a} \into \Gamma(\T \AA) \qquad \Psi(\gamma,\lambda)|_A = d_A^H\gamma + d\IP{\mu,\lambda}.$$

We notice that $\Omega^{0}(M;\Gg_E) +\Omega^2_+(M;\Gg_E)$ is isomorphic to $\Omega^{ev}_+(M;\Gg_E)$ via the map
$$\gamma+ \lambda \mapsto \gamma+\lambda + \star \gamma,\qquad \gamma \in \Omega^0(M;\Gg_E), \lambda \in \Omega^2_+(M;\Gg_E),$$
so  we can use $\frak{a} = \Omega^{ev}_+(M;\Gg_E)$ and then the extended action is given simply by
$$\Psi:\Omega^{ev}_+(M;\Gg_E) \into \Omega^{od}(M;\Gg_E)\qquad \Psi(\alpha)|_A = d_A^H \alpha. $$

Following the reduction procedure, the reduced manifold, $\mc{M}$,  is obtained by taking the quotient of $\mu^{-1}(0)$ by the action of the gauge group. In this case, $\mu^{-1}(0)$ consists of the space of anti self-dual connections and hence $\mc{M}$ is simply the moduli space of instantons.

The reduction procedure also produces a specific Courant algebroid over the reduced manifold. Namely, if we let $\KK$ be the bundle generated  by $\Psi(\frak{a})$  and $\Kperp$ its orthogonal complement \wrt\ the natural pairing, then the space of $\Gau$-invariant sections of $\Kperp/\KK$ over $\mu^{-1}(0)$ inherits a bracket and a nondegenerate pairing which make the quotient bundle $(\Kperp/\KK)/\Gau$ into a Courant algebroid over $\M$. Notice that at a specific anti self-dual connection $A$, $\KK$ is the image of $\Omega^{ev}_+(M;\Gg_E)$ by $d_A^H$. If we let $d^H_{A+}:\Omega^{od}(M;\Gg_E)\into \Omega^{ev}_+(M;\Gg_E)$ be the composition of $d_A^H$ with projection onto the \gsd-forms, then a simple integration by parts shows that $\Kperp|_A$ is the space of $d^H_{A+}$-closed forms and hence the reduced Courant algebroid is the degree one cohomology of the following elliptic complex:
\begin{equation}\label{eq:elliptic complex}
\{0\}\into \Omega^{ev}_+(M;\Gg_E)\stackrel{d_A^H}{\into}\Omega^{od}(M;\Gg_E)\stackrel{d^H_{A+}}\into \Omega^{ev}_+(M;\Gg_E)\into \{0\}.
\end{equation}

In a way, this is the double of the usual elliptic complex describing the tangent space to $M$:
\begin{equation}\label{eq:usual elliptic complex}
\{0\}\into \Omega^{0}(M;\Gg_E)\stackrel{d_A}{\into}\Omega^{1}(M;\Gg_E)\stackrel{d_{A+}}\into \Omega^{2}_+(M;\Gg_E)\into \{0\}.
\end{equation}
Earlier we added the assumption that \eqref{eq:usual elliptic complex}  above has no cohomology in degree zero. From now on we also add the assumption that this complex has no cohomology in degree two, so that $A$ corresponds to a smooth point in $\M$ and hence the dimension of $\M$  is given by the index of $d_A$, which is a topological invariant  due to the Atiyah-Singer index theorem. Further, in this case, the cohomology of \eqref{eq:elliptic complex} also concentrates in the middle term:
$$H^{od}_{d_A^H}(M;\Gg_E)= \frac{\ker(d^H_{A+}:\Omega^{od}(M;\Gg_E)\into \Omega^{ev}_+(M;\Gg_E)}{\Im(d_A^H:\Omega^{ev}_+(M;\Gg_E)\into \Omega^{od}(M;\Gg_E)}. $$

Since $M$ has an SKT structure, we can also endow $\AA$ with an SKT structure. Firstly we let $\star$ be the generalized metric: Indeed, in four dimensions $\star^2 = \Id$,  spin invariance of the Chevalley pairing means that $\star$ is orthogonal (in even dimensions) and by design we have (c.f. \eqref{eq:metric on forms})
$$\IP{\gf,\star \gf} = \int_M\kappa(\gf,\star\gf)_{Ch} > 0,\qquad\mbox{for } \gf\in \Omega^{od}(M;\Gg_E)\backslash\{0\},$$
so $\star$ is a generalized metric and  for this metric $V_+$ is the space of \gsd\ odd forms, $\Omega_+^{od}(M;\Gg_E) = (\W^2 + \W^{-2})\cap \Omega^{od}(M;\Gg_E)$,  and $V_-$ is the space of \gasd\ odd forms, $\Omega_-^{od}(M;\Gg_E) = \W^0 \cap \Omega^{od}(\M;\Gg_E)$.

For complex structure, we let $\mathbb{I} = e^{\frac{\pi \I}{2}}$ as in Lemma \ref{lem:starSKT}. Then that same lemma implies that $\mathbb{I}^2|_{V_+} = -\star = -\Id$. Since this structure is independent of the point $A \in \AA$, it is constant and hence integrable. The spaces $V_+^{1,0}$ and $V_+^{0,1}$ can be easily described: since $\I = \mathbb{I}$ on $W^{\pm 2}$ we have $V_+^{1,0} = \W^2 \cap \Omega^{od}(M;\Gg_E)$ and $V_+^{0,1} = \W^{-2}\cap \Omega^{od}(M;\Gg_E)$.

Further, it is immediate  that this SKT structure is invariant by the action of the gauge group. Now, we are in position to use the SKT reduction theorem:

\begin{theorem}[SKT Reduction Theorem \cite{MR2314216}]\label{theo:skt reduction}
Let $\Psi:\frak{a} \into \Gamma(\T \AA)$ be an extended action with moment map $\mu$ preserving an SKT structure $(\mathbb{G},\mathbb{I})$ on $\AA$ as above and let $P= \mu^{-1}(0)$. If the underlying group action on $P$ is free and proper and $0$ is a regular value of the moment map, then $\M_{red} = P/\Gau$ is smooth and the  SKT structure on $\AA$ reduces to an SKT structure on $\M_{red}$ \iff\ $\Kperp \cap V_+|_{P}$ is invariant under $\mathbb{I}$.
\end{theorem}

As we observed earlier, at a point $A \in \AA$, $\Kperp$ corresponds to the $d^H_{A+}$-closed odd forms and hence $\Kperp \cap V_+$ corresponds to the $d^H_{A+}$-closed \gsd\ odd forms, that is the \gsd, $d_A^H$-harmonic odd forms and hence according to the theorem to prove that the moduli space of instantons inherits an SKT from $M$, we must prove that the space of \gsd, $d_A^H$-harmonic odd forms is invariant under $\mathbb{I}$. To achieve this, one must develop Hodge theory for forms with coefficients. In this case, the condition  that the connection is anti self-dual, allows us to achieve the result.

Indeed, similarly to the flat case, we can split $d_A^H$ as a sum of three operators:
$$d_A^H = \delta_+^N + \bar{\delta_+^N} + \slashed{d}_{A-}^H,$$
$$\delta_+^N:\W^k \into \W^{k+2} \qquad \bar{\delta_+^N}:\W^k \into \W^{k-2} \qquad \slashed{d}_{A-}^H:\W^k \into \W^{k}.$$ 
Indeed, locally $d_A^H = d^H + A$, for some $A \in \Omega^1(M;\Gg_E)$ and hence the splitting of $d^H$ together with the splitting of $A$ into its $V_+^{1,0}$, $V_+^{0,1}$ and $V_-$ components gives the desired decomposition of $d_A^H$. Just as in Section \ref{subsec:signature}, integration by parts gives that $(\delta_+^N)^* = -\bar{\delta_+^N}$ and $\slashed{d}_{A-}^{H*} = \slashed{d}_{A-}^H$.

Now, let $\gf \in (\W^2\oplus \W^{-2}) \cap \Omega^{od}(M;\Gg_E)$. Then the $d_A^H$-Laplacian computed on $\gf$ is given by
\begin{align*}
\triangle_H \gf & = (d_A^{H*}d^H_{A+} + d_A^H d^{H*}_{A+})\gf\\
& =  (-\delta_+^N- \bar{\delta_+^N} + \slashed{d}_{A-}^H)\slashed{d}_{A-}^H + (\delta_+^N+ \bar{\delta_+^N} +\slashed{d}_{A-}^H)\slashed{d}_{A-}^H\gf\\
& =2 (\slashed{d}_{A-}^H)^2\gf.
\end{align*}
Hence the Laplacian leaves the spaces $\W^2 \cap \Omega^{od}$ and $\W^{-2} \cap \Omega^{od}$ invariant. Therefore we can decompose the space of harmonic \gsd\ odd forms into two spaces $\mc{H}^{\pm 2} = \ker(\triangle_H) \cap \W^{\pm 2}$ and $\mathbb{I}$ acts as multiplication by $i$ on $\W^2$ and by $-i$ on $\W^{-2}$, so either way it preserves the intersection $\ker(\triangle_H) \cap \W^{\pm 2}$ and hence it preserves the space of \gsd\ harmonic odd forms. According to Theorem \ref{theo:skt reduction}, this means that the SKT structure from $\AA$ reduces to an SKT structure on $\M$ so we have re-obtained the following result, originally due to L\"ubke and Teleman \cite{MR1370660}:

\begin{theorem}
Let $(M,I, [g])$ be  a compact conformal  Hermitian 4-manifold, let  $E$ be a bundle over $M$  whose structure group is compact  with Lie algebra $\Gg$ and let $\Gg_E$ be the adjoint bundle over $M$. Let $\M_{s}$ be the quotient of the space 
$$\tilde{P} = \{A \in \AA : (F_A)_+ =0; ~H^0_{d_A}(M;\Gg_E) = H^2_{d_A}(M;\Gg_E) =\{0\}\},$$
by the action of the gauge group, i.e., $\M_{s}$ is the smooth locus of the moduli space of instantons on $E$. Then  $\M_{s}$ has an SKT structure induced by the unique SKT structure on $M$ in the conformal class $[g]$.
\end{theorem}

\section{Stability}\label{sec:deformations}

\subsection{Stability of SKT structures I}\label{subsec:stability}

Similarly to the  K\"ahler case, the space of deformations of SKT structures is infinite dimensional even after taking the quotient by the action of the diffeomorphism group. Indeed, from the classical viewpoint, an SKT structure is a Hermitian structure $(g,I)$ for which $dd^c\omega =0$. For any $f \in C^{\infty}(M)$ we can consider a new Hermitian form: $\omega + \e dd^c f$. In a compact manifold if we take  $\e$ small enough this gives rise to new non diffeomorphic SKT structures on $M$ still with  3-form $d^c\omega$. 

It is therefore more natural to study the deformation problem in  the context of stability: Given a \gc\ extension of an SKT structure $(\G,\J_1)$, one is interested in which deformations of $\J_1$   give rise to deformations of the SKT structure. Since in this setup  the SKT structure and the deformation problem are on a fixed Courant algebroid, we are implicitly requiring that the class $[H]$ does not change as the structures vary. The proof of this section's main theorem is inspired on the approach used by Goto to prove the  \gk\ stability theorem  \cite{MR2669364}.

Before we state the  theorem, we must introduce the vector bundles and operators relevant  for the deformation theory of SKT manifolds. The bundle that will govern the stability problem is $\wedge^\bullet\bar{L_2}$, where $\bar{L_2} = V_+^{0,1}\oplus V_-^{1,0}$ is the $-i$ eigenspace of $\J_2$. Due to the decomposition of $\bar{L_2}$ as a direct sum, we have a bigrading on its exterior algebra, so we define
$$\wedge^{k,l}\bar{L_2} = \wedge^k V_+^{0,1}\tensor \wedge^l V_-^{1,0}.$$
The relevant differential operators are related to $\delta_\pm$ and their complex conjugates, namely we define
\begin{equation}\label{eq:delbar+def}
\begin{aligned}
\delbar_{\pm}:\Clif(\T_\C M) \into \mathrm{Diff}(\Omega^{\bullet}(M;\C));&\qquad \delbar_{\pm}\alpha = \{\deltabar_{\pm},\alpha\};\\
\del_{\pm}:\Clif(\T_\C M) \into \mathrm{Diff}(\Omega^{\bullet}(M;\C));&\qquad \del_{\pm}\alpha = \{\delta_{\pm},\alpha\};\\
\bar{\N}:\Clif(\T_\C M) \into \Clif(\T_\C M)&\qquad \bar{\N}\alpha =\{\bar{N},\alpha\};
\end{aligned}
\end{equation}
where $\mathrm{Diff}(\Omega^\bullet(M;\C))$ is the vector space of linear differential operators on forms which is itself a $\Z_2$-graded Lie algebra with graded commutator as the bracket.

The main properties of these operators are:

\begin{proposition}\label{prop:delbar_+ and del_-}
Let $(\G,\J)$ be a generalized Hermitian extension of a positive SKT structure on a manifold with 3-form $(M,H)$.  The restriction of  $\bar{\N}$ to $\wedge^{\bullet}\bar{L_2}$ vanishes and 
\begin{align}
\delbar_+:\Gamma(\wedge^{k}V_+^{0,1}\tensor  \wedge^lV_-^{1,0}\tensor \wedge^m V_-^{0,1})&\into \Gamma(\wedge^{k+1}V_+^{0,1}\tensor  \wedge^lV_-^{1,0}\tensor \wedge^m V_-^{0,1});\label{eq:delbar_+}\\
\del_-:\Gamma(\wedge^{k,l}\bar{L_2})& \into \Gamma(\wedge^{k,l+1}\bar{L_2}).\label{eq:del_-}
\end{align}
Further, the following relations and their complex conjugates hold:
$$ \delbar_+^2 = \bar{\N}^2= 0;\qquad  \{\delbar_+, \bar{\N}\}=0;$$
$$ \del_-^2 = - \{\overline{\N},\delbar_+\};\qquad \{\overline{\N},\del_-\} =0;  \qquad \{\delbar_-,\delbar_+\}=0;\qquad \{\delbar_+,\del_-\}  =- \{\delbar_-,\bar{\N}\};$$
$$ \{\del_+, \delbar_+\} + \{\del_-,\delbar_-\} + \{\N,\overline{\N}\}=0.$$
\end{proposition}

A proof of this proposition is given in the Appendix. Since $\delbar_+^2=0$ we can form the corresponding cohomology:

\begin{definition}\label{def:delbar_+ cohomology}
Let $(\G,\I)$ be a positive SKT structure on a manifold with 3-form $(M,H)$ and $\J$ be a \gc\ extension of $\I$. The {\it $\delbar_+$-cohomology in $\bar{L_2}$}, $H_{\delbar_+}^{\bullet,\bullet}(M;\bar{L_2})$, is the quotient
$$H_{\delbar_+}^{k,l}(M;\bar{L_2}) = \frac{\ker(\delbar_+:\Gamma(\wedge^{k,l}\bar{L_2})\into \Gamma(\wedge^{k+1,l}\bar{L_2}))}{\Im(\delbar_+:\Gamma(\wedge^{k-1,l}\bar{L_2})\into \Gamma(\wedge^{k,l}\bar{L_2}))} $$
\end{definition}

Now we are in position to state this section's main theorem:

\begin{theorem}\label{theo:SKT stability theorem}{\em (SKT Stability Theorem)}
Let $(\G,\J_1)$ be a generalized Hermitian extension of a positive SKT structure. Given $\J_{1\e}$, an analytic family of deformations of $\J_1$ parametrized by a small disc $D \subset \C$, $\e:D\into \Gamma(\wedge^2 \bar{L_1})$, we let 
\begin{align*}
\e^+:D& \into \Gamma(\wedge^2 V_+^{0,1})\\
\e^{\pm}:D& \into \Gamma(V_+^{0,1}\tensor V_-^{0,1})\\
\e^-:D& \into \Gamma(\wedge^2 V_-^{0,1})
\end{align*}
be the components of $\e$ and for each of them let $\e^{\bullet} = \sum_{k=1}^\infty \tfrac{1}{k!}\e^\bullet_k t^k$ be the corresponding series expansion. Then the obstructions to finding a generalized metric $\G_\e$ which makes $(\G_\e,\J_{1\e})$ into a generalized Hermitian extension of an SKT structure  lie in the space $H^{2,1}_{\delbar_+}(M;\bar{L_2})$ and the first obstruction is the class $[\del_-\e_1^{+} + \overline{\N}\e_1^{\pm}]$.
\end{theorem}

\begin{lemma}\label{lem:make it orthogonal}
There is a map $a: D \subset \wedge^2\bar{L_1} \into \wedge^2 (L_1\oplus \wedge^2\bar{L_1})_\R$, where $D$ is a small disc around the origin and $(\cdot)_\R$ denotes the real elements in the vector space, such that for each $\e \in \Gamma(\wedge^2\bar{L_1})$, $a(\e)$ is the unique element in the codomain that satisfies
$$e^{a}_*\J_1 e^{-a}_* = I_\e \J_1 I_\e^{-1}.$$
Further $a_*|_0(\e) = -(\e + \bar{\e})$.
\end{lemma}
\begin{proof}
Indeed, the space of \gcss\ of the same parity as $\J_1$  on  $T_pM$ is the homogeneous space $SO(\T_pM)/\mathrm{Stab}(\J_1) \cong SO(n,n)/U(n,n)$. Hence, composing the exponential map with the projection
$$\frak{so}(\T_p M)\into SO(\T_p M) \into SO(\T_p M)/\mathrm{Stab}(\J_1)$$
gives a local submersion.  Since the elements in $ \frak{so}(\T_p M)$ preserving $\J_1$ are those in $(L_1\tensor \bar{L_1})_\R$, and $(\wedge^2 L_1\oplus \wedge^2\bar{L_1})_{\R}$ is a complementary subspace, we have that
$$(\wedge^2 L_1\oplus \wedge^2\bar{L_1})_{\R} \into SO(\T_p M)/\mathrm{Stab}(\J_1).$$
is a local diffeomorphism. That is, for each small deformation of a \gcs\ $\J_1$ on $V$ there is a unique element $a \in  (\wedge^2 L_1\oplus \wedge^2\bar{L_1})_{\R}$ which realizes it.

For the last claim, we observe that if $\e:(-s,s)\into  \wedge^2\overline{L_1}$ is a path which passes through $0$ at time zero and $v \in L_1$, then $\dot{I_\e}|_0 (v+\bar{v}) =\dot{\e}|_0(v) + \dot{\bar{\e}}|_0(\bar{v})$. On the other hand, for $a:(-s,s)\into   \in  (\wedge^2 L_1\oplus \wedge^2\bar{L_1})_{\R}$, say $a = \alpha  + \bar{\alpha}$ with $\alpha \in \wedge^2 L_1$, we have
$$\tfrac{d}{dt}e^a_*|_0 (v+\bar{v}) =\tfrac{d}{dt}(e^a (v+\bar{v}) e^{-a})|_0 = \{a|_0,v+\bar{v}\} = -(\alpha|_0(v)+\bar{\alpha}|_0(\bar{v})),$$
hence if the family $a$ gives rise to the same deformation as $\e$, we must have $a = -(\e+\bar{\e})$.
\end{proof}

\begin{proof}[Proof of Theorem \ref{theo:SKT stability theorem}] With this setup, $e^a_*$ is an orthogonal transformation of $\T M$ and hence $e^{a}_*\J_2 e^{-a}_*$ is an \gacs\ which commutes with $\J_{1\e}$. If $b \in \R \oplus (L_1\tensor \bar{L_1})_\R \subset \Clif^2(\T M)$ the orthogonal transformation $e^b_*$ preserves $\J_1$ hence $\J_{1\e}$ and $\J_{2\e} = e^{a}_*e^{b}_*\J_2e^{-b}_*e^{-a}_*$ give rise to a generalized Hermitian structure for any choice of $b$. Our task is to identify suitable conditions on a given $a$  under which we can choose $b$ so that the pair $\J_{1\e}$ and $\J_{2\e}$ gives an SKT structure.

Since $\J_{1\e}$ is integrable, we know that
$$d^H:\U^{0,n}_\e \into \U^{1,n-1}_\e\oplus \U^{1,n-3}_\e \oplus \U^{-1,n-1}_\e \oplus \U^{-1,n-3}_\e$$
and, according to Theorem \ref{theo:genhermitian}, the SKT condition is that the $\U^{-1,n-3}_\e$ component of the map above must vanish.

Since $\J_{1\e}$ and $\J_{2\e}$ are obtained from $\J_1$ and $\J_2$ by the action of $e^a e^b \in Spin(\T M)$, the corresponding decompositions of spinors are related by
$$\U^{p,q}_\e = e^a e^b \U^{p,q},$$
so the different components of $d^H$ in the splitting induced by $(\J_{1\e},\J_{2\e})$ (c.f. Figure \ref{fig:genhermitian})
$$d^H = \delta_{+\e} + \deltabar_{+\e} + \delta_{-\e} + \deltabar_{-\e}+ N_\e + \bar{N_\e} +N_{2\e} + \bar{N_{2\e}}$$
are in correspondence with the components of $ e^{-b}e^{-a}d^H e^a e^b$ in the splitting determined by $(\J_1,\J_2)$ (see Figure \ref{fig:obstruction}).
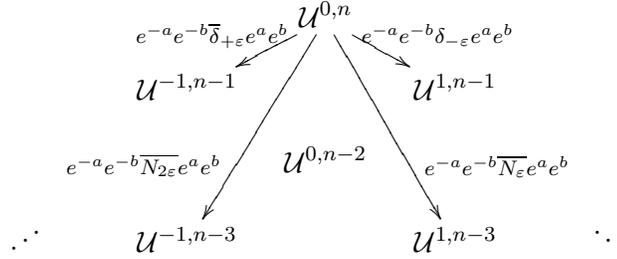
\begin{figure}[h!!]
\begin{center}
$$\xymatrix@R=12pt@C=12pt{
& & &\ar[ld]_{e^{-a}e^{-b}\deltabar_{+\e}e^a e^b} \ar[lddd]_(0.7){e^{-a}e^{-b}\bar{N_{2\e}}e^a e^b} \U^{0,n}\ar[rd]^{e^{-a}e^{-b}\delta_{-\e} e^a e^b} \ar[rddd]^(0.7){e^{-a}e^{-b}\bar{N_\e}e^a e^b}& &&\\
&& \U^{-1,n-1}              &               &\U^{1,n-1}&& \\
&&& \U^{0,n-2}              &&&\\
\Ddots&&\U^{-1,n-3}                & &\U^{1,n-3}              & &\ddots\\
}
$$
\caption{Nontrivial components of $e^{-b}e^{-a}d^He^a e^b$ and their relation to the components of $d^H$ \wrt\ the decomposition of forms determined by $(\J_{1\e},\J_{2\e})$.}\label{fig:obstruction}
\end{center}
\end{figure}

Since $a \in (\wedge^2 L_1\oplus \wedge^2\bar{L_1})_{\R}$, there are $\alpha^+\in \Gamma(\wedge^2V_+^{0,1})$, $\alpha^{\pm} \in \Gamma(V_+^{0,1}\tensor V_-^{0,1})$ and $\alpha^-\in \Gamma(\wedge^2V_-^{0,1})$ such that
\begin{equation}\label{eq:a decomposition}
a = \alpha^+ + \alpha^\pm +\alpha^- +\bar{\alpha^+} + \bar{\alpha^\pm} +\bar{\alpha^-}.
\end{equation}
Similarly, any  $b \in \Gamma(\R \oplus L_1\tensor \bar{L_1})_\R$ can be decomposed in four  components according to the splitting
$$\R \oplus (L_1\tensor \bar{L_1})_\R = \R \oplus (V_+^{1,0}\tensor V_+^{0,1})_\R \oplus (V_-^{1,0}\tensor V_-^{0,1})_\R \oplus ((V_+^{0,1}\tensor V_-^{1,0})\oplus (V_+^{1,0}\tensor V_-^{0,1})_\R)$$
however, as we will see later, only the last of these components is relevant for the deformation problem. So we let $b$ is of the form
\begin{equation}\label{eq:b decomposition}
b = \beta + \beta^{\pm} + \bar{\beta^{\pm}}
\end{equation}
with  $\beta^{\pm} \in \Gamma(V_+^{0,1}\tensor V_-^{1,0})$ and $\beta \in  \Gamma(\R \oplus (V_+^{1,0}\tensor V_+^{0,1})_\R \oplus (V_-^{1,0}\tensor V_-^{0,1})_\R)$.

Since the deformation family $\e(t)$ is analytic on $t$, so is $a(t)$  and we can try to solve the condition that the $\U^{-1,n-3}$ component of $e^{-b}e^{-a}d^He^ae^b \psi $ vanishes by a power series argument. In the sequence given an analytic function $f(t)$, where $f$ can be $\alpha^+$, $\alpha^\pm$ and so on,  we denote denote the coefficients of its power series expansion by $f_k$, so that $f = \sum \tfrac{1}{k!}f_k t^k$. 

The first nontrivial condition on $a$ imposed by  the integrability of $\J_1$ is that $[d^H,a_1] $ must have no components in $\U^{n-3,\bullet}$. In terms of the decomposition of $a$ above, this implies
$$[\deltabar_+,\alpha^+_1]=0\qquad  [\deltabar_+ ,\alpha^{\pm}_1] + [\deltabar_-, \alpha^+_1] =0$$
or equivalently
\begin{equation}\label{eq:integrability1}
\delbar_+\alpha^+_1 =0\qquad \delbar_+ \alpha^{\pm}_1 + \delbar_-\alpha^+_1=0.
\end{equation}

The first term in the series expansion of $e^{-b}e^{-a}d^He^ae^b \psi$ is 
$$[d^H,b_1+ a_1]\psi.$$
Splitting $a$ and $b$ into their components, we see that the requirement that the $\U^{-1,n-3}$ component vanishes is equivalent to
$$([\deltabar_+,\beta^{\pm}_1] +  [\delta_-,\alpha^+_1] + [\overline{N},\alpha^{\pm}_1])\psi =0$$ 
Isolating $\beta^{\pm}$  and using that $\psi$ is nonzero, we can rewrite  the condition above as
$$\delbar_+\beta^{\pm}_1 =  -\del_-\alpha^+_1 - \bar{\N}\alpha^{\pm}_1.$$ 
And one obvious necessary condition for this to have a solution is that the right hand side must be $\delbar_+$-closed. That is indeed the case, since 
$$ \delbar_+(\del_-\alpha^+_1 + \overline{\N}\alpha^{\pm}_1) =  - \del_-\delbar_+ \alpha^+_1  - \delbar_-\bar{\N}\alpha^+_1 - \bar{\N}\delbar_-\alpha^+_1 -\overline{\N}\delbar_+\alpha^{\pm}_1 =0$$
where for the first equality we used the commutation rule for $\delbar_+$ and $\del_-$ from Proposition \ref{prop:delbar_+ and del_-} and in the second equality we used equations \eqref{eq:integrability1} and the fact that $\N \alpha^+_1=0$ as $\alpha^+_1 \in \bar{L_2}$. If this $\delbar_+$-closed form is in fact $\delbar_+$-exact, we can use  $G_{\delbar_+}$, the Green operator for $\delbar_+$, to find a suitable $\beta_1^{\pm}$:
$$\beta_1^{\pm} =  -\delbar_+^* G_{\delbar_+} (\del_-\alpha^+_1 + \bar{\N}\alpha^{\pm}_1) $$
 
Finally, we recall that according to Lemma \ref{lem:make it orthogonal}, $a_1 = -(\e_1 + \bar{\e_1})$, so the components of $a_1$ are determined by the different components of $\e_1$ hence  the equation $\beta^{\pm}_1$ must solve is
$$\delbar_+\beta^{\pm}_1=  \del_-\e^+_1 + \overline{\N}\e^{\pm}_1,$$ 
so the first obstruction class is $[\del_-\e^+_1 + \overline{\N}\e^{\pm}_1 ] \in H_{\delbar_+}^{2,1}(M;\bar{L_2})$.

Now we move to the general case, which is proved by induction. We assume that we have chosen $a_i$ and $b_i$ for $i<k$ such that the  component of $e^{-b}e^{-a}d^He^a e^b$ mapping $\U^{0,n}$ into $\U^{-1,n-3}$ vanishes to order $k-1$. Then the vanishing of the order $k$ component of this map is the condition
$$[\beta_k^{\pm},\deltabar_+] + [\alpha_k^+,\delta_-] + [\alpha_k^{\pm},\bar{N}] + F_k(a_1,\cdots, a_{k-1},b_1,\cdots,b_{k-1})= 0,$$
where $F_k$ is some  function. This condition is equivalent to
$$\bar{\del_+}\beta_k^{\pm} =-( \del_- \alpha_k^++ \bar{\mc{N}}\alpha_k^{\pm}+ F_k(a_1,\cdots, a_{k-1},b_1,\cdots,b_{k-1})).$$
And if there is $\beta_k^{\pm}$ which solves this equation, we can proceed to the next step. If for all possible choices of $b_i$ for $i<k$ which guarantee the vanishing of the lower order terms this equation has no solution, then the deformation is obstructed.

Next we show that $ \del_- \alpha_k^++ \bar{\mc{N}}\alpha_k^{\pm}+ F_k(a_1,\cdots, a_{k-1},b_1,\cdots,b_{k-1})$ is $\delbar_+$-closed so we can conclude that the obstruction space is $H^{2,1}(M;\bar{L_2})$. Since $\delbar_+ \beta^{\pm}_k$ is clearly $\delbar_+$-closed, we must prove that
$$\delbar_+(\delbar_+\beta_k^{\pm} +  \del_- \alpha_k^++ \bar{\mc{N}}\alpha_k^{\pm}+ F_k(a_1,\cdots, a_{k-1},b_1,\cdots,b_{k-1}))=0$$
but the term in parenthesis is precisely the order $k$ component of $e^{-b}e^{-a}\bar{N_{2\e}}e^a e^b$ (see Figure \ref{fig:obstruction}) and by assumption, the lower order terms of $e^{-b}e^{-a}\bar{N_{2\e}}e^a e^b$ vanish, therefore we must prove that
\begin{equation}\label{eq:middle step}
\delbar_+ e^{-b}e^{-a}\bar{N_{2\e}}e^a e^b = \{\deltabar_+,e^{-b}e^{-a}\bar{N_{2\e}}e^a e^b\}
\end{equation}
vanishes to order $k$. But since $\deltabar_+=e^{-b}e^{-a}\deltabar_{+\e}e^{b}e^{a} +O(t)$, we see that \eqref{eq:middle step} vanishes to order $k$ \iff\ 
$$\{e^{-b}e^{-a}\deltabar_{+\e}e^{a}e^{b},e^{-b}e^{-a}\bar{N_{2\e}}e^a e^b\}$$
vanishes to order $k$. But for this operator, we have
$$\{e^{-b}e^{-a}\deltabar_{+\e}e^{a}e^{b},e^{-b}e^{-a}\bar{N_{2\e}}e^a e^b\} = e^{-b}e^{-a}\{\delta_{+\e},\bar{N_{2\e}}\}e^a e^b = 0,$$
 where in the last equality we used that, in a generalized Hermitian manifold, the condition $(d^H)^2=0$ gives, among other things, that $\{\deltabar_+,\bar{N_2}\}=0$.
 
 Finally, if all obstructions vanish, standard elliptic estimates show that the sequence constructed above converges. See for example \cite{MR2669364} for the proof of convergence for the analogous problem in the \gk\ case. 
\end{proof}

Of course, the same arguments with the obvious changes give deformation results for negative SKT structures. Similarly, Propositions \ref{prop:kahlerupq} and  \ref{prop:relation to Dolbeault} can be used to translate our theorem to classical SKT structures:

\begin{theorem}
Let $(M,g,I)$ be a compact SKT manifold with $d^c \omega =H$. Given $I_\e$, an analytic family of deformations of the complex structure parametrized by a small disc $D \subset \C$, the obstructions to finding a metric $g_\e$ which makes $(g_\e,I_{\e})$ into an SKT structure   whose 3-form represents the class $[H]$ lie in the Dolbeault cohomology space $H^{1,2}_{\delbar}(M)$.
\end{theorem}
\begin{proof}
Under the isomorphism $\Psi$ of Proposition \ref{prop:kahlerupq}, $\bar{L_2}$ is identified with $T^*_\C M$. More precisely, we have the following identifications
$$V_+^{0,1} \stackrel{\cong}{\into} T^{*1,0}M;\qquad V_-^{1,0}  \stackrel{\cong}{\into} T^{*0,1}M.$$
Hence $\wedge^{k,l}\bar{L_2} = \wedge^{k,l}T^*M$ and due to Proposition \ref{prop:relation to Dolbeault}, $\deltabar_+$ is identified with $\del$  so we have
$$H^{k,l}_{\delbar_+}(M;\bar{L_2})  = H^{k,l}_\del(M) = H^{l,k}_\delbar(M).$$
In particular we see that the obstruction space is the Dolbeault cohomology $H^{1,2}_\delbar(M)$.
\end{proof}

\begin{example}
Let $H$ be a Hopf surface. As a manifold, $H$ is diffeomorphic to $S^1 \times S^3$. Since for every compact complex surface the Fr\"olicher spectral sequence degenerates at the first page \cite{MR0239114}, the Dolbeault cohomology is isomorphic to the de Rham cohomology as graded vector spaces. Further the identities for the Euler characteristic and signature
$$0 =\chi = \sum (-1)^{p+q} h^{p,q}\qquad \mbox{and}\qquad 0 = \sigma =\sum (-1)^q h^{p,q}$$
give us $h^{0,0} = h^{0,1}=h^{2,1} = h^{2,2} =1$ and the remaining Hodge numbers vanish. Since $h^{1,2}=0$ we conclude that SKT structures on the Hopf surface are stable. Of course, due to Gauduchon's result on the existence of SKT structures on any compact complex surface \cite{MR0470920}, this conclusion is not particularly strong.

Due to K\"unneth's formula for the Dolbeault cohomology, if $M_n$ is a product of $n$ Hopf surfaces (with the product complex structure) we have that $h^{1,2}(M_n) = 0$ and hence again the SKT deformation problem is unobstructed and this claim {\it does not} follow from Gauduchon's theorem. Since we can deform the product complex structure into one that is not a product, this produces nontrivial examples of SKT structures on $\times_n (S^1\times S^3)$.

Finally, the Dolbeault cohomology of the Lie group $\mathrm{Spin}(4) \cong S^3 \times S^3$ with a left invariant complex structure compatible with a bi-invariant metric is
$$H_\delbar(\mathrm{Spin}(4)) = \wedge \mathrm{span}\{x^{0,1}, z^{2,1}\} \tensor \C[y^{1,1}]/y^2,$$
where $a^{p,q}$ is a generator of bidegree $(p,q)$ \cite{MR0202713}, hence $H^{1,2}(M)$ is 1-dimensional, generated by $x^{0,1}\cup y^{1,1}$. Using the Kunneth formula, we have that $(S^1\times S^1) \times (S^3 \times S^3)$ with the product structure has $h^{1,2} =3$, showing that for this structure, the SKT deformation problem is potentially obstructed. The problem is in fact obstructed and one way to understand this is by observing that a left invariant complex structure on $S^3 \times S^3$ compatible with a bi-invariant metric can be deformed into one that is not compatible with any bi-invariant metric and such structure ceases to be SKT.
\hfill$\blacksquare$
\end{example}

\begin{example}[Rank two, simple Lie groups]
While the proof that  a left invariant complex structure on $S^3 \times S^3$ which is not orthogonal \wrt\ a bi-invariant metric can not be made into an SKT structure seems to be just a computation, for rank two, compact, simple Lie groups  the nature of this  obstruction is a little clearer. 

Firstly, we observe that if $I$ is indeed orthogonal \wrt\ the bi-invariant metric on a simple Lie group, then the 3-form $H$ is the Cartan 3-form
$$H(X,Y,Z) = \IP{[X,Y],Z}.$$
For such a form, the $d^H$-cohomology of the Lie group vanishes \cite{cavalcanti-2004}.

On the other hand, if $G$ has rank 2, the Dolbeault cohomology of $G$ depends on whether the complex structure is orthogonal \wrt\ the Killing form or not:

\begin{theorem}[Pittie \cite{MR994129}, Proposition 4.5]
Let $G$ be a rank two, compact, simple Lie group of complex dimension $n+1$. Then the Dolbeault cohomology of a left invariant complex structure $I$ on $G$ is given, as a graded vector space, by
$$H^{\bullet}_\delbar(G) = \wedge^{\bullet} \mathrm{span}(x^{0,1},z^{2,1})\tensor \C[y^{1,1}]/(y^{1,1})^n$$
if $I$ is compatible with the bi-invariant metric and by
$$H^{\bullet}_\delbar(G) = \wedge^{\bullet} \mathrm{span}(x^{0,1},z^{n,n-1})\tensor \C[y^{1,1}]/(y^{1,1})^2$$
otherwise, where $a^{p,q}$ is a generator of bidegree $(p,q)$.
\end{theorem}
In the case when $I$ is compatible with the metric, then $\del \omega$ represents the class $z^{2,1}$ above and it is immediate that the $\delbar^{i\del\omega}$-cohomology vanishes, as it should, since in an SKT manifold this cohomology is isomorphic to the $d^H$-cohomology due to Corollary \ref{cor:cohomology isomorphisms}.

On the other hand, for a complex structure which not compatible with the metric, $H^{2,1}_\delbar(G) =H^{3,0}_\delbar(M)=\{0\}$ and since $G$ is simple, $H^2(G) =\{0\}$ hence the complex structure can not be extended to an SKT one, due to Theorem \ref{theo:obstruction2}. Similarly to Example, \ref{ex:calabi-eckman}, one can also arrive at the conclusion also using Corollary \ref{cor:cohomology isomorphisms}. Indeed, if the complex structure were part of an SKT structure then  $\del \omega$ would  be $\delbar$-exact and the $\delbar^{i\del\omega}$-cohomology would be isomorphic to the Dolbeault cohomology. But since $H^{0,1}_\delbar(G) = \C$ we see that the Dolbeault cohomology is not isomorphic, as a graded vector space, to either the usual de Rham cohomology of $G$ or the $d^H$ cohomology of $G$. Hence, regardless of the background 3-form chosen, this complex structure is not part of an SKT structure.

Of course, one can deform a complex structure compatible with the metric into one that is not (by changing the choice of complex structure on the Lie algebra of the maximal torus) and the argument above shows that in this case not only do we have $H^{2,1}_\delbar(G) \neq \{0\}$, but also that the obstruction map is also nontrivial. 
\hfill$\blacksquare$
\end{example}

\subsection{Stability of generalized K\"ahler structures}\label{subsec:gkdeformations}

Similarly to the SKT deformation problem studied in Section \ref{subsec:stability}, a common way to study the question of deformations of a (generalized) K\"ahler  structure  $(\J_1,\J_2)$ is in the context of stability:  For which deformations of $\J_1$ is there a corresponding deformation of $\J_2$ such that the pair of deformed structures is still a (generalized) K\"ahler structure?

In its classical setting, this question was successfully answered by Kodaira \cite{MR0115189}, where $\J_1$ is taken to be the complex structure and in the generalized setting, for analytic deformations,  by Goto \cite{MR2669364} under the additional hypothesis that the canonical bundle of $\J_2$ has a globally defined nowhere vanishing closed section. In both cases, there are no obstructions to the problem. However, even in the K\"ahler setting, if  we let $\J_1$ be the symplectic structure and $\J_2$ be the complex,  the problem does have obstructions \cite{MR2198780}.

As we have seen in Section \ref{sec:intrinsic torsion}, a \gks\ is at the same time a positive and a negative SKT structure. Hence, when studying the stability problem of a \gks\ one can expect three possible outcomes:
\begin{enumerate}
\item The deformation is unobstructed and one can complete the deformed $\J_{1}$ to a \gk\ structure;
\item The deformation is half obstructed and one can complete the deformed $\J_1$ to a positive or negative SKT structure;
\item The deformation is fully obstructed and $\J_1$ can not be completed to either a  positive or a negative SKT structure.  
\end{enumerate}

\begin{theorem}[Generalized K\"ahler Stability Theorem]\label{theo:gkstability}
Let $(M,\J_1,\J_2)$ be a compact \gkm. Given $\J_{1\e}$, an analytic family of deformations of $\J_1$ parametrized by a small disc $D \subset \C$, $\e:D\into \Gamma(\wedge^2 \bar{L_1})$, we let 
\begin{align*}
\e^+:D& \into \Gamma(\wedge^2 V_+^{0,1})\\
\e^{\pm}:D& \into \Gamma(V_+^{0,1}\tensor V_-^{0,1})\\
\e^-:D& \into \Gamma(\wedge^2 V_-^{0,1})
\end{align*}
be the components of $\e$ and for each of them let $\e^{\bullet} = \sum_{k=1}^\infty \tfrac{1}{k!}\e^\bullet_k t^k$ be the corresponding series expansion. Then the obstructions to finding a \gcs\ $\J_{2\e}$ which makes $(\J_{1\e},\J_{2\e})$ into a \gks\  lie in the space
$$\mc{H}^{(2,1)+(1,2)} = \frac{\ker(d_{L_2}:\Gamma(\wedge^{2,1}\bar{L_2}\oplus \wedge^{1,2}\bar{L_2}))}{\Im(d_{L_2}:\Gamma(\wedge^{1,1}\bar{L_2}))},$$
and the first obstruction class is $[\delbar_+\bar{\e_1^-} +\del_-\e_1^+] $.

\end{theorem}
\begin{proof}
According  to Lemma \ref{lem:make it orthogonal}, there is an element $a \in \Gamma(\wedge^2 L_1\oplus \wedge^2 \bar{L_1})_{\R}$ such that
$$e^{-a}_* \J_1e^a_* = e^{-\e}_*\J_1 e^{-\e}$$
and, as in Theorem \ref{theo:SKT stability theorem}, one can compose the orthogonal transformation $e^a_*$ with an orthogonal tranformation $e^b_*$ which preserves $\J_1$ and still obtain the same deformation, that is, for any $b \in \Gamma(L_1\tensor \bar{L_1})_{\R}$,  $\J_{1\e}$ and $\J_{2\e}$ are given by
$$ \J_{1\e} =  e^{-a}_* \J_1e^a_* \qquad \mbox{and}\qquad  \J_{2\e} = e^{-a}e^{-b}\J_2 e^b e^a.$$
As before,  our quest is to find $b$ for which $\J_{2\e}$ is integrable and once again we do so using a power series argument. The canonical bundle of $\J_{2\e}$ is given by $e^a e^b U^{0,n}$ and integrability is equivalent to
$$e^{-b}e^{-a}d^He^{a}e^b: \U^{0,n}\into \U^{n-1,1}\oplus \U^{n-1,-1}.$$
If we decompose $a$ and $b$ into their $V_{\pm}^{1,0}$ and $V_{\pm}^{0,1}$ components as in \eqref{eq:a decomposition} and \eqref{eq:b decomposition} and write the corresponding components in power series, then the first condition we must solve is the linear vanishing of the components of $\bar{N_{\J_2\e}}$,  the Nijenhuis tensor of $\J_{2\e}$
$$e^{-b}e^{-a}d^H e^a e^b:\U^{0,n}\into \U^{1,n-3}\oplus \U^{-1,n-3}.$$
The linear term which maps  $\U^{0,n}$ into $\U^{-1,n-3} \oplus \U^{1,n-3}$ is
$$d_{L_2}\beta_1^{\pm} + \del_-\alpha_1^+ +\delbar_+\bar{\alpha_1^-} $$
So  we see that the vanishing of the linear part of the Nijenhuis tensor is equivalent  the following condition
\begin{equation}\label{eq:obstruction}
d_{L_2}\beta_1^{\pm} = -(\delbar_+\bar{\alpha_1^-} +\del_-\alpha_1^+).
\end{equation}
Integrability of $\J_{1\e}$ means that the $\U^{\pm3,n-3}$ components of  $e^{-a}e^{-b}d^He^ae^b|_{\U^{0,n}}$ vanish and the linear part of these components  is $\delbar_+\alpha_1^+:\U^{0,n}\into \U^{-3,n-3}$ and $\del_-\bar{\alpha_1^-}:  \U^{0,n}\into \U^{3,n-3}$, hence $\alpha_1^+$ is $\delbar_+$-closed and $\bar{\alpha_1^-}$ is $\del_-$-closed, showing that $\delbar_+\bar{\alpha_1^-} +\del_-\alpha_1^+$ is $d_{L_2}$-closed and hence represents a class in $\mc{H}^{(2,1)+(1,2)}$. If this class  vanishes, we can find $\beta^{\pm}_1$ which solves \eqref{eq:obstruction}. That is the first obstruction class is 
$$[\delbar_+\bar{\e_1^-} +\del_-\e_1^+]  \in \mc{H}^{(2,1)+(1,2)}.$$

Next we assume that we have chosen $\beta^{\pm}_j$ for $j<k$  such that the Nijenhuis tensor of $\J_2$ vanishes to order $k-1$ and we must choose $\beta^{\pm}_k$ which makes the Nijenhuis tensor of $\J_2$ vanish to order $k$. Expanding the operator $e^{-b}e^{-a}d^He^ae^b$ to order $k$ we see that the degree $k$ part of the Nijenhuis  tensor is:
\begin{equation}\label{eq:stepk1}
(\bar{N_{\J_2\e}})_k = d_{L_2}\beta_k^{\pm} + \del_-\alpha_k^+ + \delbar_+\bar{\alpha_k^-}+ F_k(a_1,\cdots, a_{k-1},b_1,\cdots, b_{k-1}),
\end{equation}
for some function $F_k$ which  takes values in $\Gamma(\wedge^{2,1}\bar{L_2}\oplus \wedge^{1,2}\bar{L_2})$, and the vanishing of $\bar{N_{\J_2\e}}$ to order $k$ is equivalent to the requirement that the term above vanishes, which can only be achived if 
$$d_{L_2}\beta_k^{\pm} = - (\del_-\alpha_k^+ + \delbar_+\bar{\alpha_k^-}+ F_k(a_1,\cdots, a_{k-1},b_1,\cdots, b_{k-1})).$$
Next we  prove that the right hand side above is $d_{L_2}$-closed, which allows us to conclude that  the obstruction to find $\beta_k^{\pm}$ lies in $\mc{H}^{(2,1)+(1,2)}$.

We let
$$b_{<k+1} = \sum_{j<k+1}\tfrac{1}{j!}(\beta^{\pm}_j +\bar{\beta^{\pm}}_j)t^j.$$
Since $e^{-a}\J_1 e^{a} = e^{-a}e^{-b_{<{k+1}}}\J_1 e^{b_{<k+1}}e^a$ is integrable for any choice of $b_k$, we have that the operator
$$e^{-b_{<k+1}}e^{-a} d^H e^{a} e^{b_{<k+1}}|_{\U^{0,n}}$$
has two components,  say $d^H|_{\U^{n,0}} = \delbar_{\J_2 \e} +  \bar{N_{\J_2\e}}$,
$$\delbar_{\J_2 k} = \delbar_{\J_2} + \bar{v}(t), \qquad \bar{v(t)} \in \Gamma(\bar{L_2}),~ \bar{v}(0)=0;$$ 
$$\bar{N_{\J_2\e}}  = (\bar{N_{\J_2 \e}})_k t^k+  O(t^{k+1}) \in \Gamma(\wedge^{1,2}\bar{L_2}).$$
The condition $(d^H)^2 =0$ implies $\{\delbar_{\J_2 \e},\bar{N_{\J_2\e}} \}=0$, hence, in particular, the degree $k$ part of this operator vanishes, i.e.,
$$ 0= \{\delbar_{\J_2},(\bar{N_{\J_2\e}})_k \}= d_{L_2}(\bar{N_{\J_2\e}})_k.$$
\end{proof}

\begin{remark}
There are natural maps
$$\pi_{2,1}:\mc{H}^{(2,1)+(1,2)} \into H_{\delbar_+}^{2,1}(M;\bar{L_2}) \qquad{and}\qquad  \pi_{1,2}:\mc{H}^{(2,1)+(1,2)} \into H_{\del_-}^{1,2}(M;\bar{L_2}).$$
given by taking the appropriate components of a class in $\mc{H}^{(2,1)+(1,2)}$. It is possible, however, that a nontrivial class in $\mc{H}^{(2,1)+(1,2)}$ is mapped into the trivial class by both projections. In terms of deformations, this means that there may be no obstruction to deforming a \gks\ as either a positive or a negative SKT structure, while the \gk\ deformation problem itself is obstructed. A simple example of this phenomenon happens if we regard a K\"ahler structure  $(I,\omega)$ as a generalized K\"ahler structure with $\J_1$ determined by the symplectic structure and $\J_2$ by the complex structure: On the one hand not all deformations $\omega_t$ of $\omega$ can be accompanied by a deformation of $I$ to keep the pair generalized K\"ahler. On the other hand, keeping $I$ fixed makes the pair  $(\omega_t,I)$ into a Hermitian symplectic structure as positivity of the tensor $\omega(I\cdot,\cdot)$ is an open condition. We can regard this Hermitian symplectic structure as either a positive or a negative SKT structure, hence the obstructions to deformations of the SKT structures vanish.
\end{remark}

\subsection{Stability of SKT structures II}

An inherent part of the setup used for the study of the stability problem of SKT structures in Section \ref{subsec:stability}, was that the cohomology class of 3-form was fixed. Depending on the reader's upbringing, that hypothesis may seem unnatural and one might want to study the problem of stability of a classical SKT structure $(g,I)$ without the requirement that $[d^c\omega]$ is fixed. As we will see next, this is actually an easier problem as it has perfectly acceptable solutions which use only classically available tools.

\begin{definition} Let $(M,I)$ be a complex manifold. The {\it $(p,q)$-Aeppli cohomology} of $M$ is 
$$H^{p,q}_A(M) = \frac{\mathrm{ker}(\del\delbar:\Omega^{p,q}(M)\into \Omega^{p+1,q+1}(M))}{\mathrm{Im}(\del:\Omega^{p-1,q}(M)\into \Omega^{p,q}(M) )+ \mathrm{Im}(\delbar:\Omega^{p,q-1}(M)\into \Omega^{p,q}(M))},$$
and the {\it Aeppli numbers} are $h^{p,q}_A = \mathrm{dim}(H^{p,q}_A(M))$.
\end{definition}

\begin{theorem}
Let $(M,g,I)$ be an SKT structure and  for $\e \in (-\delta, \delta)$, let  $I_\e$ be a smooth family of complex structures with $I_0 = I$. If the Aeppli number $h^{1,1}_A$ is constant on $\e$, then there is a family of metrics $g_\e$ such that $(g_\e,I_\e)$ is an SKT structure for small values of $\e$. In particular, if $h^{1,1}_A$ is constant  in  a \nhood\ of $I$ in the moduli space of complex structures, then every nearby complex structure is part of and SKT structure.
\end{theorem}
\begin{proof}
If $h^{1,1}_A$ is constant in $\e$, we can form the family of  Green operators $G_\e$ and projection onto harmonics, $\mathscr{H}_\e$, for the elliptic complex
$$\Omega^{1,0}(M)\oplus \Omega^{0,1}\stackrel{\delbar+ \del}{\into}\Omega^{1,1}(M)\stackrel{\del\delbar}{\into}\Omega^{2,2}(M)$$ 
Since $h^{1,1}_A$ is constant, the form $\omega_\e = (\mathscr{H}_\e + (\del+ \delbar)(\del+ \delbar)^*G_\e)  \pi^{1,1}_\e\omega$ varies smoothly on $\e$, where $\pi^{1,1}_\e$ is the projection from $\Omega^2(M)$ onto $\Omega^{1,1}(M)$ \wrt\ the complex structure $I_\e$. The operator $(\mc{H}_\e + (\del+ \delbar)(\del+ \delbar)^*G_\e) $ is just the projection onto $\del\delbar$-closed forms. We see immediately that $\omega_0 = \omega$, that $\omega_\e$ is of type $(1,1)$ and $\del\delbar$-closed. Since $\omega(\cdot,I\cdot)$ is positive definite, the same is true for $\omega_\e(\cdot,I_\e\cdot)$ as long as $I_\e$ is not too far from $I$.
\end{proof}

\begin{example}[The Iwasawa manifold]
The Iwasawa manifold $M$  is the quotient of the 3-dimensional complex Heiseberg group $\mathbb{H}$  by the left action of the subgroup, $\Lambda$, of matrices whose entries are Gaussian integers: $M = \mathbb{H}/\Lambda$. The question of existence of SKT structures on the underlying differentiable manifold as well as the stability of such structures was studied by  Fino, Parton and Salamon \cite{MR2059435}: while the standard complex structure is not part of an SKT structure there are other complex structures on $M$ which are. Further, there are deformations of the SKT complex structures which cease to be SKT.

From our point of view, this is to be expected and amounts to the observation that for the SKT structures, $h^{1,1}_A(M) =7$ while for generic complex structures $h^{1,1}_A(M)= 6$ \cite{angella2013}. It is also worth noticing that for the families of SKT structures provided by  Fino {\it et al} the cohomology class of the 3-form changes along the family \cite{baarsmaphd}.\hfill$\blacksquare$

\end{example}

\section{Appendix - SKT differentials}\label{sec:Lie algebroids}

Throughout this section we let $(M,\G,\J)$ be a manifold with a generalized complex extension of an SKT structure and let  $\delbar_\pm$, $\del_\pm$, $\mc{N}$ and $\bar{\mc{N}}$ be the operators defined in \eqref{eq:delbar+def}.

\vskip6pt
\noindent
{\bf Proposition \ref{prop:delbar_+ and del_-}.}
The restriction of  $\bar{\N}$ to $\wedge^{\bullet}\bar{L_2}$ vanishes and 
\begin{align}
\delbar_+:\Gamma(\wedge^{k}V_+^{0,1}\tensor  \wedge^lV_-^{1,0}\tensor \wedge^m V_-^{0,1})&\into \Gamma(\wedge^{k+1}V_+^{0,1}\tensor  \wedge^lV_-^{1,0}\tensor \wedge^m V_-^{0,1});\tag{\ref{eq:delbar_+}}\\
\del_-:\Gamma(\wedge^{k,l}\bar{L_2})& \into \Gamma(\wedge^{k,l+1}\bar{L_2}).\tag{\ref{eq:del_-}}
\end{align}
Further, the following relations and their complex conjugates hold:
$$ \delbar_+^2 = \bar{\N}^2= 0;\qquad  \{\delbar_+, \bar{\N}\}=0;$$
$$ \del_-^2 = - \{\overline{\N},\delbar_+\};\qquad \{\overline{\N},\del_-\} =0;  \qquad \{\delbar_-,\delbar_+\}=0;\qquad \{\delbar_+,\del_-\}  =- \{\delbar_-,\bar{\N}\};$$
$$ \{\del_+, \delbar_+\} + \{\del_-,\delbar_-\} + \{\N,\overline{\N}\}=0.$$
\vskip6pt
\begin{proof}
The relations between the operators follow from the corresponding relations for $\deltabar_+$, $\delta_-$ and $\bar{N}$ and the fact that the vector space of operators on forms are a $\Z_2$-graded Lie algebra with Lie bracket given by the graded commutator of operators. If we let $\alpha \in \Clif(\T_\C M)$, using the graded Jacobi identity we have
$$\{d^H,\{d^H,\alpha\}\}  = \{\{d^H,d^H\},\alpha\} -\{d^H,\{d^H,\alpha\}\}=  -\{d^H,\{d^H,\alpha\}\} =0,$$
where we have used that $\{d^H,d^H\} = 2(d^H)^2 =0$. Hence, splitting $d^H$ into its six components, taking $\alpha \in \Gamma(\wedge^jV_+^{1,0} \tensor \wedge^{k}V_+^{0,1}\tensor  \wedge^lV_-^{1,0}\tensor \wedge^m V_-^{0,1})$, applying this operator to elements in $\U^{p,q}$ and imposing the vanishing of the different components of the result we get the relations stated in the theorem.

Since $\bar{N} \in \Gamma(\wedge^3\bar{L_2})$ and $\bar{L_2}$ is isotropic, $\{\bar{N},\alpha\}$ vanishes for every $\alpha \in \Gamma(\wedge^{\bullet}\bar{L_2}).$ 

Next we will prove that \eqref{eq:delbar_+}  holds and will ommit the proof of \eqref{eq:del_-} as it follows the same principles. Firstly, we prove that, when restricted to the domains indicated in \eqref{eq:delbar_+}, $\delbar_+$  actually takes values in the Clifford algebra of $\TM$. Indeed, for $\alpha \in \Clif(\T M)$, $\gf \in \U^{p,q}$ and $f\in C^{\infty}(M)$  we  decompose $df$ into its $V_+^{1,0}$, $V_+^{0,1}$, $V_-^{1,0}$ and $V_-^{0,1}$ components, say, $\del_+ f$,  $\delbar_+ f$, $\del_- f$  and $\delbar_- f$, respectively and  compute
\begin{equation}\label{eq:basic computation}
\begin{aligned}
\{d^H,\alpha\} (f\gf) & = d^H \alpha f \gf -(-1)^{|\alpha|} \alpha d^H (f\gf)\\
& =  df\wedge \alpha  \gf + f d^H\alpha \cdot \gf -(-1)^{|\alpha|} \alpha (df \wedge \gf + f d^H\gf)\\
& =  (\del_- f + \delbar_- f +\del_+ f + \delbar_+ f) \cdot \alpha  \gf + f d^H\alpha  \gf \\ &-(-1)^{|\alpha|} \alpha ((\del_- f + \delbar_- f + \del_+ f + \delbar_+f)\cdot \gf + fd^H\gf).\\
& =  f \{d^H\alpha\} \gf  + (\del_- f + \delbar_- f +\del_+ f + \delbar_+ f) \cdot \alpha  \gf \\&- (-1)^{|\alpha|}\alpha ((\del_- f + \delbar_- f + \del_+ f + \delbar_+f)\cdot \gf).\\
\end{aligned}
\end{equation}
If $\alpha \in \Gamma(\wedge^{k}V_+^{0,1}\tensor  \wedge^lV_-^{1,0}\tensor \wedge^m V_-^{0,1})$, then $\delbar_+ f$ graded commutes with $\alpha$ and hence if we compose $\{d^H,\alpha\}$ with the projection onto $\U^{p-k+l-m-1,q-k-l+m-1}$ we have
$$\pi_{\U^{p-k+l-m-1,q-k-l+m-1}}\circ \{d^H, \alpha\}: \U^{p,q} \into \U^{p-k+l-m-1,q-k-l+m-1}$$
is $C^{\infty}$-linear, that is, it is a tensor. Further, using the decomposition $d^H$ into its six components, as in Theorem \ref{theo:genhermitian}, we see that the tensor above corresponds to:
$$\delbar_+\alpha = \{\deltabar_+,\alpha\}: \U^{p,q} \into \U^{p-k+l-m-1,q-k-l+m-1}$$
Since the action of $\Clif(\T M)$ establishes an isomorphism between $\Clif(\TM)$ and  the space of endomorphisms of $\wedge^{\bullet}T^*M$, the tensor $\delbar_+\alpha$ corresponds to an element in $\Clif(\TM)$:
$$\delbar_+:\Gamma(\wedge^{k}V_+^{0,1}\tensor  \wedge^lV_-^{1,0}\tensor \wedge^m V_-^{0,1}) \into \Clif(\T_\C M)$$

Next we observe that a Leibniz rule holds. Indeed, not only is the space of graded linear differential operators on $\Omega^{\bullet}(M)$ a graded Lie algebra with the graded commutator as a bracket, but also this bracket satisfies a Leibniz rule, \wrt\ the usual composition of operators. In our case, say, for $\alpha,\beta \in \Clif(\T_\C M)$, this translates to
$$\{\deltabar_+,\alpha\circ \beta\} = \{\deltabar_+,\alpha\}\circ \beta + (-1)^{|\alpha|+1}\alpha\circ \{\deltabar_+,\beta\}.$$
Or, in terms of $\delbar_+$, the Leibniz rule holds:
$$\delbar_+(\alpha\circ \beta) = \delbar_+(\alpha)\circ \beta + (-1)^{|\alpha|+1}\alpha\circ \delbar_+(\beta).$$
Therefore, in order to prove, for example, \eqref{eq:delbar_+} it is enough to establish the result when the argument of $\delbar_+$ is a section of any of the bundles $V_+^{0,1}$, $V_-^{1,0}$ or $V_-^{0,1}$.

For sections of $\bar{L_1} = V_+^{0,1}\oplus V_-^{0,1}$, we notice that $L_1$ is a Lie algebroid whose dual space is naturally isomorphic to $\bar{L_1}$ hence we have a Lie algebroid differential
$$d_{L_1}: \Gamma(\bar{L_1}) \into \Gamma(\wedge^2 \bar{L_1})$$
and the decomposition $\bar{L_1} = V_+^{0,1}\oplus V_-^{0,1}$ gives a corresponding decomposition of $d_{L_1}$ in two components. Then one can readily see that $\delbar_+$ is nothing but the component of $d_{L_1}$ mapping which increases the $V_+^{0,1}$-degree:
$$\delbar_+:\Gamma(V_+^{0,1}) \into \Gamma( \wedge^2V_+^{0,1}),\qquad \delbar_+:\Gamma(V_+^{0,1}) \into \Gamma( V_+^{0,1}\wedge  V_-^{0,1}).$$ 

Therefore to finish the proof we only need to show that
\begin{equation}\label{eq:representation1}
\delbar_+:\Gamma(V_-^{1,0})\into \Gamma(V_+^{0,1}\tensor V_-^{1,0}).
\end{equation}

\begin{lemma}\label{lem:representation1}
The vector bundle $V_-^{1,0}$ is a representation of the Lie algebroid $V_+^{1,0}$ if endowed with the connection
$$\nabla:\Gamma(V_+^{1,0}\tensor V_-^{1,0})\into \Gamma(V_-^{1,0}); \qquad \nabla_{v} w = \pi_{V_-^{1,0}}(\Cour{v,w}),$$
where  $v \in \Gamma(V_+^{1,0})$ and $ w \in \Gamma(V_-^{1,0})$. That is, for $v, v_1,v_2\in \Gamma(V_+^{1,0})$, $w \in \Gamma(V_-^{1,0})$ and $f \in C^{\infty}(M)$ we have
$$\nabla_{fv} w = f \nabla _vw,\qquad \nabla_v fw = (\mc{L}_{\pi_T(v)} f ) w + f \nabla_v w$$
$$\nabla_{v_1}\nabla_{v_2}-\nabla_{v_2}\nabla_{v_1} - \nabla_{\Cour{v_1,v_2}} =0.$$
\end{lemma}
\begin{proof}
The first two identities follow from basic properties of the Courant bracket:
$$\nabla_{fv} w = \pi_{V_-^{1,0}}(\Cour{fv,w}) =  \pi_{V_-^{1,0}}(f\Cour{v,w} -( \mc{L}_{\pi_T(w)}f) v) = \pi_{V_-^{1,0}}(f\Cour{v,w} ) = f \nabla _v w.$$
$$\nabla_v fw =  \pi_{V_-^{1,0}}(\Cour{v,fw}) =  \pi_{V_-^{1,0}}(f\Cour{v,w} +( \mc{L}_{\pi_T(v)}f) w) =  f \nabla_v w + (\mc{L}_{\pi_T(v)} f ) w.$$

Now we prove flatness. Let us first consider the term $\nabla_{v_1} \nabla_{v_2} w$. We decompose $\Cour{v_2,w}$ into two components, according to the decomposition $L_1 = V_+^{1,0} \oplus V_-^{1,0}$ and $\nabla_{v_2} w$ is the component in $V_-^{1,0}$:
$$\Cour{v_2,w} = u_+^{1,0} + \nabla_{v_2} w,\qquad u_+^{1,0} \in \Gamma(V_+^{1,0}).$$
Since $V_+^{1,0}$ is involutive, we see that
$$\Cour{v_1, u_+^{1,0}} \in  V_+^{1,0}.$$
Hence the $V_-^{1,0}$ component of $\Cour{v_1,\nabla_{v_2} w}$ is the $V_-^{1,0}$ component of $\Cour{v_1,\Cour{v_2, w}}$, that is
$$\nabla_{v_1} \nabla_{v_2} w = \pi_{V_-^{1,0}}(\Cour{v_1,\Cour{v_2, w}}).$$
Therefore we have
$$\nabla_{v_1}\nabla_{v_2}w-\nabla_{v_2}\nabla_{v_1}w - \nabla_{\Cour{v_1,v_2}}w  = \pi_{V_-^{1,0}}(\Cour{v_1,\Cour{v_2, w}} - \Cour{v_2,\Cour{v_1, w}} - \Cour{\Cour{v_1,v_2}, w}) =0,$$
where in the last equality we used the fact that the Courant bracket satisfies the Jacobi identity.
\end{proof}
Since $\nabla$ is a connection, we think of it dually as an operator
$$\nabla:\Gamma(V_-^{1,0})\into\Gamma(V_+^{0,1}\tensor V_-^{1,0}).$$
The following lemma establishes \eqref{eq:representation1} and hence finishes the proof of the proposition.

\begin{lemma}\label{lem:nabla and delbar_+}
For $\alpha \in \Gamma(V_-^{1,0})$ we have
$$\delbar_+\alpha = \nabla \alpha.$$
\end{lemma}
\begin{proof}
Indeed, we check that when paired (i.e., graded commuted) with an element $v\in \Gamma(V_+^{1,0})$ these two operators give the same result. Unwinding the definition of $\delbar_+$,  we have that for $\gf \in \U^{p,q}$ $\{v,\delbar_+ \alpha\} \gf$ is the $\U^{p+1,q-1}$ component of $\{v,\{d^H, \alpha\}\} \gf$.
\begin{align*}
\{v,\delbar_+ \alpha\} \gf &= \pi_{U^{p+1,q-1}}(\{v,\{d^H, \alpha\}\} \gf)\\
&=\pi_{U^{p+1,q-1}}(\Cour{v, \alpha} \gf)\\
&=\pi_{V_-^{1,0}}\Cour{v,\alpha} \gf\\
&=\{v,\nabla\alpha\} \gf,
\end{align*}
where in the second equality we used that the Courant bracket is the derived bracket associated to $d^H$.

Since $\delbar_+ \alpha$ and $\nabla \alpha$ are both elements in $\Clif(\T_\C M)$ we have just established that their difference lies in the annihilator of $V_+^{1,0}$ under graded commutator, that is
$$  \nabla \alpha-\delbar_+ \alpha \in \Gamma(\wedge^{\bullet}V_+^{1,0}\tensor\wedge^{\bullet}V_-).$$
And we have that when acting in $\U^{p,q}$
$$\nabla \alpha - \delbar_+ \alpha : \U^{p,q} \into \U^{p,q-2} \subset \W^{p+q-2}$$
as individually both $\delbar_+ \alpha$ and $\nabla \alpha$ have this property. On the other hand, elements in $\wedge^{l}V_+^{1,0}\tensor\wedge^{m}V_-$  map $\U^{p,q}\subset \W^{p+q}$ into $\W^{p+q+l}$, so the difference $\nabla \alpha - \delbar_+ \alpha $ vanishes and the result follows.
\end{proof}
\noqed\end{proof}

\bibliographystyle{amsplain}
\bibliography{references}

\end{document}